\documentclass[a4paper, oneside]{article} 
\usepackage[a4paper, left=17mm, right=17mm, top=25mm, bottom=25mm]{geometry}
\usepackage{amsfonts}
\usepackage{amsmath}
\usepackage{graphicx}
\usepackage{multicol}
\usepackage{mathtools}
\usepackage{dsfont}
\usepackage{xparse}
\usepackage{caption}
\usepackage{algorithm}
\usepackage{algpseudocode}
\usepackage{color}
\usepackage{listings}
\usepackage{xcolor}
\usepackage{enumitem}
\usepackage{hyperref}

\newcommand{\hold}{\mathcal{D}}
\newcommand{\Sb}{S}
\newcommand{\Om}{\Omega}
\newcommand{\VV}{\theta}
\newcommand{\R}{\mathds{R}}
\newcommand{\ddo}{:}
\newcommand{\singset}{L}

\newcommand{\fopt}{\texttt{FormOpt} }
\newcommand{\fx}{\texttt{FEniCSx}}

\newcommand{\id}{\rm{id}}

\newtheorem{definition}{Definition}
\newtheorem{proposition}{Proposition}

\hypersetup{
    colorlinks=true,
    linkcolor=blue,
    citecolor=blue,
    urlcolor=blue
}

\newcommand{\code}[1]{\texttt{#1}}

\title{FormOpt: A FEniCSx toolbox for level set-based shape optimization supporting parallel computing}
\author{
	Josu\'{e} D. D\'{i}az-Avalos\thanks{Email: \texttt{josue.diazavalos@uni-due.de}} \qquad
	Antoine Laurain\thanks{Email: \texttt{antoine.laurain@uni-due.de}} \\
	\\
	Faculty of Mathematics, University of Duisburg-Essen
}
\date{\phantom{\today}}

\begin{document}

\maketitle

\begin{abstract}
	This article presents the toolbox \fopt for two- and three-dimensional shape optimization with parallel computing capabilities, built on the \fx\,  software framework. 
	We introduce fundamental concepts of shape sensitivity analysis and their numerical applications, mainly for educational purposes, while also emphasizing computational efficiency via parallelism for practitioners.
	We adopt an optimize-then-discretize strategy based on the distributed shape derivative and its tensor representation, following the approach of \cite{MR3843884} and extending it in several directions.
The numerical shape modeling relies on  a level set method, whose evolution is driven by a descent direction computed from the shape derivative.
Geometric constraints are treated accurately through a Proximal-Perturbed Lagrangian approach.
\fopt leverages the powerful features of  \fx, particularly its support for weak formulations of partial differential equations, diverse finite element types, and scalable parallelism.
The implementation supports three different parallel computing modes: data parallelism, task parallelism, and a mixed mode. Data parallelism exploits \fx's mesh partitioning features, and we implement a task parallelism mode which is useful for problems governed by a set of partial differential equations with varying parameters. The mixed mode conveniently combines both strategies to achieve efficient utilization of computational resources.
\end{abstract}

\begin{multicols}{2}

\section{Introduction}

Shape and topology optimization \cite{MR2731611,SokZol92} has built over the years a solid theoretical basis and has found its way in many industrial applications for optimal design thanks to the fast development of computational capacities. 
The topic has found natural applications first  in structural and fluid mechanics, then has been developed for various other partial differential equations (PDEs) constraints. 

Different approaches have been developed for shape and topology optimization. 
The solid isotropic microstructure with penalty (SIMP) method \cite{MR2008524}, based on a material density approach,  is the most popular in structural mechanics and is the basis of many educational codes since the work \cite{Sigmund2014}. 
It has been extended to different frameworks, beyond structural optimization.
Another widely used approach is the level-set paradigm, with different variations of the level-set method introduced by Osher and Sethian in \cite{MR965860}. 
Other important methods include phase-field  \cite{MR2108636} and topological derivative-based approaches \cite{MR2235384}. 
Reference works for the level set approach in structural optimization are \cite{MR1911658,MR2033390}. 
Given the vast literature on shape and topology optimization, we provide only a few illustrative references here and refer the reader to review articles for a more comprehensive overview. 

The range of software available for shape and topology optimization is rapidly expanding.
Along commercial softwares, numerous open-source tools have been developed.
For topology optimization in structural mechanics, the trend of educational codes began in 2001 with Sigmund's 99-line \texttt{MATLAB} code \cite{Sigmund2014}, followed by an improved version in \cite{Andreassen2010}. Another pioneering contribution is the \texttt{FreeFEM} implementation \cite{MR2252743}.
Since then, codes have been released for various platforms, including \texttt{NGsolve} \cite{MR4233697}, \texttt{FEniCS} \cite{fenicsx, BLAUTH2021100646,MR3843884} and \texttt{Firedrake}~\cite{Paganini2021}. 
An open-source \texttt{Python} implementation of the Null Space Optimizer is also available \cite{MR4686093}. For a comparison of algorithms, we refer to \cite{Gain2013}.
 
A significant portion of the literature and educational implementations in structural and topology optimization is based on the SIMP method \cite{Sigmund2014}.
Another widely used class of methods is the level-set approach, which exists in several variants, as reviewed in \cite{vanDijk2013}. Some formulations employ a Heaviside or smoothed Heaviside function to represent the domain~\cite{MR1951408}, while others are built upon the concept of the topological derivative \cite{MR2235384,Novotny2013,MR4916810}.
A further line of level set-based methods relies on shape sensitivity analysis in the infinite-dimensional setting, using shape derivatives \cite{MR2731611,SokZol92}. Educational codes following this approach include \cite{Challis2009,MR2459656}. These methods typically employ the boundary representation of the shape derivative, known as the Hadamard form, which must then be extended into the domain for use within a level-set framework~\cite{MR2225309}.
In the present work, we base our approach on the so-called {\it distributed} or {\it weak shape derivative} \cite{LAURAIN2020328,MR3535238}.
This formulation offers several advantages including higher accuracy \cite{MR2642680,MR3348199},  weaker regularity requirements \cite{LAURAIN2020328}, and a natural framework for domain extension. An educational code for structural mechanics based on this approach was first proposed in \cite{MR3843884}. The present work represents a natural continuation and extension of \cite{MR3843884}, incorporating additional and modernized features.
Our overall objective is to provide a framework for learning the fundamental principles of shape sensitivity analysis and their numerical implementation, with an educational focus. 
The \fx\ platform \cite{fenicsx} employs the Unified Form Language (UFL) to represent weak formulations of PDEs, offering an intuitive notation that closely mirrors the underlying mathematics.
At the same time, we place strong emphasis on computational efficiency by leveraging the parallel capabilities of \fx, enabling fast and practical testing of shape optimization problems.
A distinguishing feature of our approach is its emphasis on the optimize-then-discretize paradigm: the shape derivative is derived in the infinite-dimensional setting, and a discrete version is then used in the numerical algorithm.

With the steady increase in computational power, parallel computing has become not only commonplace but often essential for tackling large-scale problems. 
This trend is also visible in the topology optimization community, where several recent educational codes incorporate parallel capabilities. 
Examples include a \texttt{Julia}-based toolbox \cite{Wegert2025}, a parallel \texttt{FreeFEM} framework~\cite{MR4916810}, a PETSc  framework \cite{Aage2014}, and a \fx\ implementation for 2D and 3D topology optimization \cite{Jia2024}. The latest version of \texttt{cashocs}, built on FEniCS, has also introduced MPI-based parallelism \cite{BLAUTH2021100646,BLAUTH2023101577}.
In this work, we propose a natural extension of the distributed shape derivative approach introduced in \cite{MR3843884} by incorporating parallel computing capabilities.
Unlike many educational papers on topology optimization, which focus primarily on structural mechanics, our aim is to provide a broader perspective by addressing shape optimization problems subject to different types of PDE constraints.
A particularly relevant class of problems considered here are inverse problems, where task parallelism plays a key role, since these problems typically involve families of PDE constraints with varying parameters \cite{MR4261115}. Standard applications in  structural mechanics and linear elasticity are also included as a particular case.

We extend the work of \cite{MR3843884} in several directions. 
First, in the present framework we rely exclusively on finite elements, whereas \cite{MR3843884} combined finite elements for solving the PDE constraints with finite differences for solving the Hamilton–Jacobi equation governing domain evolution. This unified finite element approach enables more flexible geometries and makes more efficient use of \fx\ resources.
Second, while \cite{MR3843884} focused exclusively on linear elasticity, we address more general PDE constraints and illustrate the framework through a broader set of examples. 
We also provide a systematic way of handling volume constraints without resorting to penalization in the cost functional.
Another major improvement is the integration of parallel computing, which allows for substantial speedups in numerical experiments. 
Three paradigms of parallelism are supported: data parallelism, distributing computations across multiple processes; task parallelism, enabling for example the simultaneous solution of the same PDE with different sources and mixed parallelism, combining both approaches.
An additional objective of this paper is to offer an accessible introduction to parallelization and its application to PDE-constrained shape optimization.
The foundation of our framework remains the tensor representation of the shape derivative, as in~\cite{MR3843884}. Practically, this means that when considering a new problem, the main task is to specify the appropriate PDE and geometry, and to analytically compute the corresponding tensors of the distributed shape derivative.
The \fopt files are available for download at 
\href{https://github.com/JD26/FormOpt}{\texttt{github.com/JD26/FormOpt}}, 
along with several tutorials.

The structure of the paper is as follows. 
In Section~\ref{sec:shape-opt}, the model shape optimization with constraints is presented, and the basic notions of distributed and boundary expressions of the shape derivative, as well as the  tensor representation
are introduced. 
In Section~\ref{sec:model_prob}, an inverse problem for linear elasticity is presented, which allows us to explain general principles of shape derivative computation through a concrete example.
In Section~\ref{sec:numerics}, the numerical implementation is discussed in details, and various snippets of the code are used for reference. In particular, the discretization of the level set equation and the parallelism paradigms are explained.
In Section~\ref{sec:results}, numerical results are given for several model problems, including benchmarks of compliance minimization for linear elasticity, and performance tests to assess the parallel scalability are performed.

\section{Model problem and shape derivatives}\label{sec:shape-opt}

Let $\mathcal{D} \subset \mathds{R}^d$ be an open, bounded domain with boundary
$\Gamma \coloneqq \partial \mathcal{D}$.
Throughout the paper, $\mathcal{D}$ is fixed, while $\Om \subset \overline{\mathcal{D}}$ represents the variable domain subject to optimization.
The spaces $H^1(\mathcal{D})$ and $H^1_0(\mathcal{D})$ denote the usual Sobolev spaces.
The notation $\chi_{\Omega}:\overline{\mathcal{D}} \rightarrow \mathds{R}$ stands for
the indicator function of $\Omega$.
The notation $n$ is used for both the outward unit normal vector to $\hold$ and to $\Om$. 
Let $I$ denote the $d$-dimensional identity matrix and $\left|\cdot\right|_{\infty}$ the infinity norm.
Given a vector-valued function $w \in H^1(\mathcal{D})^d$,
$Dw$ denotes its Jacobian matrix.
We use $\boldsymbol{0}$ for the zero vector of $\mathds{R}^d$ and 
$\boldsymbol{1}$ for the all-ones vector of $\mathds{R}^{N_c}$,
where $N_c$ is the number of constraints.

We consider the following generic shape optimization problem with constraints:
\begin{equation} \label{eq:problem0}
	\min_{\Omega\subset\mathcal{D}} \mathcal{J}(u,\Omega)
	\quad \text{subject to} \quad C(\Omega) = \boldsymbol{1}, \, e(u,\Om)=0,
\end{equation}
where $\mathcal{J}: \mathbb{H}\times \mathds{P}\to\R$ is a cost functional, 
$C: \mathds{P}\to\R^{N_c}$ is a vector-valued constraint function,
$e: \mathbb{H}\times \mathds{P}\to\mathbb{H}^*$ is a PDE constraint for $u$,
$\mathds{P}$ is a set of subsets of $\hold$,
$\mathbb{H}$ is an appropriate function space, and $\mathbb{H}^*$ is its dual space.
The additional constraint $C(\Omega) = \boldsymbol{1}$ is typically a geometric constraint, such as a volume or perimeter constraint.
Assuming $u(\Om)$ is the unique solution of $e(u,\Om)=0$,  we introduce the reduced cost functional $J(\Om):= \mathcal{J}(u(\Om),\Om)$ and we consider the problem:
\begin{equation} \label{eq:problem}
	\min_{\Omega\subset\mathcal{D}} J(\Omega)
	\quad \text{subject to} \quad C(\Omega) = \boldsymbol{1}.
\end{equation}

For numerical purposes, we need a notion of derivative with respect to the shape $\Om$.
We thus consider a diffeomorphism $T_t:\overline{\hold}\to\overline{\hold}$, $t\in [0,\tau]$ and the associated parameterized shape $\Om_t := T_t(\Om)$.
The transformation must be at least bi-Lipschitz, as we will use a change of variables in integrals later to compute the shape derivative.
We denote its time derivative at $t=0$ as $\VV: = \frac{\partial}{\partial t} T_t|_{t=0}$.
Conversely, one can also start with a given $\VV\in\Theta^k(\hold)$, where
$$ \Theta^k(\hold) := \{\VV\in C^k_c(\R^d,\R^d)| \VV\cdot n_{|\partial\hold\setminus \singset} = 0 \mbox{ and } \VV_{|\singset} = 0\}$$
and $\singset\subset\partial\hold$ is the set of points where the normal $n$ is not defined, i.e., the set of singular points of $\partial\hold$.
Then, one can build $T_t$ satisfying $\frac{\partial}{\partial t} T_t|_{t=0} = \VV$ using a flow, as in the {\it speed method}
\cite{MR2731611,SokZol92} or using a {\it perturbation of identity} \cite{MR2512810}.
The choice of the method is irrelevant for the computation of the first-order shape derivative.
In any case, once $T_t$ is available and assuming  $\VV\in\Theta^k(\hold)$, we can define the shape derivative of $J$ as follows.

\begin{definition}[Shape derivative]\label{def1}
	Let $J : \mathds{P} \rightarrow \R$ be a shape function.
	\begin{itemize}
	\item[(i)] The Eulerian semiderivative of $J$ at $\Omega$ in direction $\theta \in \Theta^k(\hold)$, when the limit exists,
	is defined by
	\begin{equation}
	dJ(\Omega ;\VV):= \lim_{t \searrow 0}\frac{J(\Omega_t)-J(\Omega)}{t}.
	\end{equation}
	\item[(ii)] $J$ is  \textit{shape differentiable} at $\Omega$ if it has a Eulerian semiderivative at $\Omega$ for all $\theta \in \Theta^k(\hold)$ and the mapping
	\begin{align*}
	dJ(\Omega): \Theta^k(\hold) &  \to \R,\; \VV     \mapsto dJ(\Omega ;\VV)
	\end{align*}
	is linear and continuous, in which case $dJ(\Omega)$ is called the \textit{shape derivative} at $\Omega$.
	\end{itemize}
\end{definition}
	When both the cost functional and the PDE constraints are expressed as bulk integrals, the corresponding shape derivative can likewise be formulated as a bulk integral, and often admits the following canonical form.
\begin{definition}[Tensor representation]\label{def:tensor2}
	Let $\Om\in \mathds{P}$ be open.
	A shape differentiable function $J : \mathds{P} \rightarrow \R$  admits  a tensor representation of order $1$ if
	there exist tensors $\Sb_0\in L^1(\hold,\R^{d})$  and $\Sb_1\in L^1(\hold,\R^{d\times d})$, such that
	\begin{align}
	\label{ea:volume_form2}
	dJ(\Omega ;\VV) = \int_\hold \Sb_1 \ddo D\VV +  \Sb_0\cdot \VV,
	\end{align}
	for all  $\VV\in \Theta^k(\hold)$. 
\end{definition}
	Expression \eqref{ea:volume_form2} is called distributed, volumetric, domain, or weak expression of the shape derivative.
	Tensor representations of arbitrary order can be defined in a similar way, see \cite{MR3535238}.
The order of the representation is essentially the maximal order of the differential operators appearing either in the variational formulation of the PDE constraint or in the cost functional, see \cite{MR4531392,LL2024} for examples of tensor representations of order two. 
Tensor representations can also be formulated when the cost function and PDE constraints contain boundary terms, see \cite{MR3535238}.
In this work, we focus exclusively on representations of order one and bulk integrals in order to keep the presentation concise.

	Under natural regularity assumptions, the shape derivative only depends on the restriction of the normal component $\VV\cdot n$ to the interface $\partial\Om$.  
	This fundamental result is known as the Hadamard-Zol\'esio structure theorem; see \cite[pp. 480-481]{MR2731611}.  
Applying the divergence theorem, this structure follows immediatly from the tensor representation \eqref{ea:volume_form2}, as shown in the following proposition.

	\begin{proposition}[Hadamard form]\label{tensor_relations}
	Let $\Om\in \mathds{P}$ and assume $\partial \Omega$ is $C^2$. 
	Suppose that $dJ(\Omega)$ has the tensor representation \eqref{ea:volume_form2}.
	If $\Sb_1^+\in W^{1,1}(\Om,\R^{d\times d})$ and $\Sb_1^-\in W^{1,1}(\hold\setminus\overline \Omega,\R^{d\times d})$, 	then we obtain the so-called {\it Hadamard form} or {\it boundary expression} of the shape derivative: 
	\begin{equation}
	\label{eq:general_boundary_exp}
	dJ(\Omega;\theta) = \int_{\partial \Omega} G\, \theta\cdot n,
	\end{equation}
	with
	$
	G :=  [(\Sb_1^+-\Sb_1^-)n]\cdot n,
	$
	where $+$ and $-$ denote the traces on $\partial\Om$ of the restrictions of $\Sb_1$ to $\Omega$ and $\hold\setminus\overline \Omega$, respectively. 
	\end{proposition}
	See \cite[Proposition~1]{LAURAIN2020328} and \cite[Proposition~4.3]{MR3535238} for proofs of Proposition \ref{tensor_relations} in more general settings.
	The Hadamard formula \eqref{eq:general_boundary_exp} is the basis of most shape derivative–based numerical methods \cite{MR2252743,MR2033390,vanDijk2013}. It provides a natural way to express the shape derivative and derive a descent direction for optimization algorithms via the shape gradient $G$, and it integrates well with the original level-set method formulation \cite{MR965860} based on the Hamilton–Jacobi equation.
	However, this approach has certain drawbacks, such as stronger regularity requirements on $\Sb_1$ in Proposition~\ref{tensor_relations} compared to Definition~\ref{def:tensor2}. 
	By contrast, the distributed shape derivative provides a natural framework for domain extension, and several recent studies have shown that it achieves higher accuracy \cite{MR2642680,MR3348199}.
In this work, we pursue the distributed expression-based approach described in \cite{MR3843884} using \eqref{ea:volume_form2}.
	A Lagrangian approach is employed to compute the tensor representation \eqref{ea:volume_form2}.
	
Throughout the paper, we will use the following notation for the shape derivatives of the cost and constraint functionals:
	\begin{equation}\label{dJ1}
		dJ(\Omega;\VV) = \int_{\mathcal{D}} S_0^J \cdot \theta + S_1^J : D\theta, 
	\end{equation}
	\begin{equation}\label{dC1}
		dC(\Omega;\VV) = \int_{\mathcal{D}} S_0^C \cdot \theta + S_1^C : D\theta .
	\end{equation}

\section{Shape derivative for a model problem}
\label{sec:model_prob}

In order to illustrate a concrete example of the calculation of a tensor representation \eqref{ea:volume_form2}, we consider a standard inverse problem in linear elasticity, where the objective is to reconstruct the shape of an inclusion with different elastic material parameters, from boundary and/or domain measurements.
The uniqueness and stability of identifying a rigid inclusion in an elastic body from a boundary measurement have been investigated in \cite{Morassi2009}, while the case of cavities was addressed in \cite{MR1772729}.
Shape sensitivity analysis for a related problem was studied in \cite{MEFTAHI20151554}, but with a different cost functional and using strong, rather than weak, formulations of the PDEs. The emphasis in \cite{MEFTAHI20151554} is on the Hadamard form \eqref{eq:general_boundary_exp} of the shape derivative, rather than on the distributed formulation \eqref{ea:volume_form2}.
Further examples of tensor representation of the shape derivative for other PDEs are provided in Section~\ref{sec:results}.

\textbf{Inverse problem.}
Let $\Gamma_0$, $\Gamma_1$ be subsets of $\Gamma:=\partial\hold$ such that $\Gamma=\Gamma_0\cup\Gamma_1$. 
Let $H_{\Gamma_0}(\mathcal{D})^d$ be the space of functions
\begin{equation*}
	H_{\Gamma_0}(\mathcal{D})^d\coloneqq\left\{ u\in H(\mathcal{D})^{d}\mid u|_{\Gamma_0}=\boldsymbol{0}\right\}. 
\end{equation*}
The strain and stress tensors are given by, respectively,
\begin{equation*}
	\epsilon(w) \coloneqq \frac{Dw+{Dw}^{\top}}{2},\quad
	\sigma(w) \coloneqq \lambda\,\mathrm{tr}(\epsilon(w))I + 2\mu\,\epsilon(w),
\end{equation*}
where $\mu$ and $\lambda$ are the Lam\'e parameters.
Let $f\in H^{1/2}(\Gamma_1)^d$ and $g\in H^{1/2}(\Gamma_1)^d$ be vector-valued functions defined on $\Gamma_1$.
Let $\alpha$ and $\beta$ be positive weights.
\begin{figure}[H]
	\centering
\includegraphics[width=0.7\linewidth]{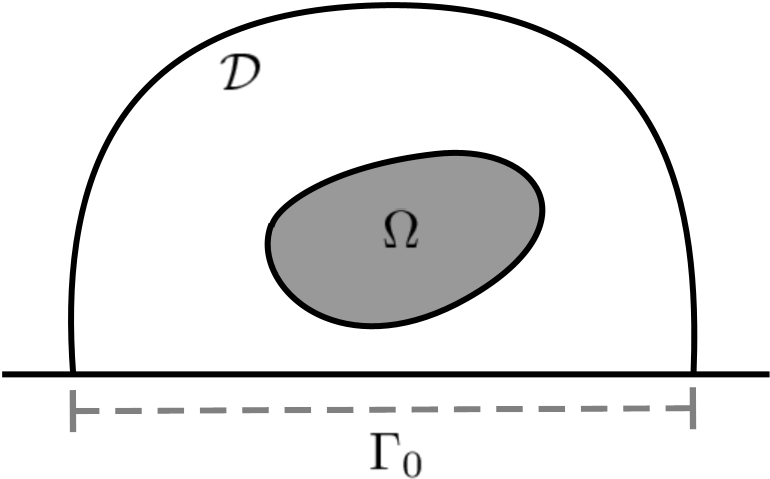}
\captionof{figure}{
	Geometric configuration for the inverse problem in elasticity.
}
\label{fig:elast_dom}
\end{figure}
We follow a Kohn-Vogelius-type approach \cite{KohnVogelius1984}, considering the shape optimization problem:
\begin{equation} \label{eq:cost}
	\min_{\Omega \subset \mathcal{D}}
	J(\Omega) \coloneqq
	\frac{\alpha}{2}\int_{\mathcal{D}}|u-y|^{2} + 
	\frac{\beta}{2}\int_{\Gamma_{1}}|u-y|^{2}
\end{equation}
where $u$ and $y$ are the solutions to the transmission problems
\begin{equation} \label{eq:f-prob}
	\begin{aligned}
    -\mathrm{div} A_{\Omega} \sigma(u) &= 0 && \text{in } \Omega \text{ and } \mathcal{D}\setminus\overline{\Omega} \\
                                     u &= 0 && \text{on } \Gamma_{0} \\
                 A_{\Omega} \sigma(u)n &= f && \text{on } \Gamma_{1}   \\
		     (A_{\Omega} \sigma(u)n)^+ &= (A_{\Omega} \sigma(u)n)^- && \text{on } \partial\Omega\\      
				                   u^+ &= u^- && \text{on }  \partial\Omega      				              
	\end{aligned}	
\end{equation}
and
\begin{equation} \label{eq:g-prob}
	\begin{aligned}
    -\mathrm{div} A_{\Omega} \sigma(y) &= 0 && \text{in } \Om  \text{ and } \mathcal{D}\setminus\overline \Omega \\
                                     y &= 0 && \text{on } \Gamma_{0} \\
                                     y &= g && \text{on } \Gamma_{1}   \\
			 (A_{\Omega} \sigma(y)n)^+ &= (A_{\Omega} \sigma(y)n)^- && \text{on } \partial\Omega\\      
								   y^+ &= y^- && \text{on } \partial\Om       
	\end{aligned}
\end{equation}
respectively, and $A_{\Omega}:\overline{\mathcal{D}} \rightarrow \mathds{R}$ is given by
\begin{equation}\label{eq:AOmega}
A_{\Omega} = \kappa\chi_{\Omega} + \chi_{\overline{\mathcal{D}}\setminus\Omega}
\quad \text{with }\kappa>0.
\end{equation}
See Figure \ref{fig:elast_dom} for an illustration of the geometry. In \eqref{eq:cost}, observe that $u$ and $y$ depend on $\Omega$ due to
the presence of $A_\Omega$ in \eqref{eq:f-prob},\eqref{eq:g-prob}.
We present the model for a single measurement for simplicity. 
For several measurements, one can simply sum over the measurements in \eqref{eq:cost}.

\textbf{The state and adjoint problems.}
Now we write the weak formulations of \eqref{eq:f-prob} and \eqref{eq:g-prob},
and of their corresponding adjoint problems.
The weak formulation of~\eqref{eq:f-prob} reads:
Find $u \in  H^1_{\Gamma_0}(\mathcal{D})^d$ such that
\begin{equation} \label{eq:w-f}
	\int_{\mathcal{D}} A_{\Omega}\sigma(u):\epsilon(w) =
	\int_{\Gamma_1} f \cdot w
	\quad
	\forall \, w \in  H^1_{\Gamma_0}(\mathcal{D})^d.
\end{equation}
It is convenient to lift the non-homogeneous Dirichlet condition $y=g$ on $\Gamma_1$ in \eqref{eq:g-prob}.
To this end, we assume that there exists an extension $g\in H^1(\hold)^d$, using the same notation for simplicity.
This allows us to rewrite \eqref{eq:g-prob} using the equivalent weak formulation:
Find $v \in  H^1_{0}(\mathcal{D})^d$ such that
\begin{equation} \label{eq:w-g}
	\int_{\mathcal{D}} A_{\Omega}\sigma(v + g):\epsilon(w) = 0
	\quad
	\forall \, w \in  H^1_{0}(\mathcal{D})^d,
\end{equation}
and ${v + g} =y$ solves \eqref{eq:g-prob}.
Considering the cost functional~\eqref{eq:cost}
and the weak formulations, we write
the Lagrangian functional
as follows:
\begin{equation} \label{eq:elastic_lag}
	\mathcal{L}(\Omega, \varphi, \psi, \rho, \varrho) =
	\mathcal{L}_{\mathrm{cost}} +
	\mathcal{L}_{\mathrm{state}[f]} +
	\mathcal{L}_{\mathrm{state}[g]},
\end{equation}
with
\begin{align*}
	\mathcal{L}_{\mathrm{cost}}    &= \frac{\alpha}{2}\int_{\mathcal{D}} \left|\varphi-\psi - g\right|^{2} + \frac{\beta}{2}\int_{\Gamma_{1}} \left|\varphi - g \right|^{2}, \\
	\mathcal{L}_{\mathrm{state}[f]}   &= \int_{\mathcal{D}} A_{\Omega} \sigma(\varphi) : \epsilon(\rho) - \int_{\Gamma_{1}} f \cdot \rho, \\
	\mathcal{L}_{\mathrm{state}[g]} &= \int_{\mathcal{D}} A_{\Omega} \sigma(\psi + g) : \epsilon(\varrho),
\end{align*}
where $\Omega\subset\mathcal{D}$ and $\varphi, \psi, \rho, \varrho$ are functions in $H^1_{\Gamma_0}(\mathcal{D})^d$.

By differentiating $\mathcal{L}$ with respect to $\varphi$ and $\psi$,
we obtain the adjoint equations:
Find $p \in  H^1_{\Gamma_0}(\mathcal{D})^d$ and $ q\in H^1_0(\mathcal{D})^d$ such that
\begin{align}
\label{dL1}	\partial_{\varphi}\mathcal{L}(u,v,p,q)(r)&=0 & \forall \, r \in  H^1_{\Gamma_0}(\mathcal{D})^d,\\
\label{dL2}	\partial_{\psi}\mathcal{L}(u,v,p,q)(r)&=0 & \forall \, r \in  H^1_0(\mathcal{D})^d.
\end{align}
From \eqref{dL1} we obtain the first adjoint problem:
Find $p \in  H^1_{\Gamma_0}(\mathcal{D})^d$ such that
\begin{align*} \label{eq:f-adj}
\int_{\mathcal{D}} A_{\Omega}\sigma(p):\epsilon(r) & = -\alpha \int_{\mathcal{D}} (u-v-g) \cdot r \\
&\quad -\beta \int_{\Gamma_1} (u - g) \cdot r \quad
\forall \, r \in  H^1_{\Gamma_0}(\mathcal{D})^d.
\end{align*}
From \eqref{dL2} we get the second adjoint problem:
Find $q \in  H^1_{0}(\mathcal{D})^d$ such that
\begin{equation*} \label{eq:g-adj}
	\int_{\mathcal{D}} A_{\Omega}\sigma(q):\epsilon(r) =
	\alpha \int_{\mathcal{D}} (u-v-g) \cdot r \quad \forall \, r \in  H^1_0(\mathcal{D})^d.
\end{equation*}
Observe that both the state and the adjoint problems are linear.

\textbf{Shape derivative components.}
Let $T_t:\overline{\hold}\to\overline{\hold}$, $t\in [0,\tau]$ be a diffeomorphism as described in Section~\ref{sec:shape-opt} and $\Om_t := T_t(\Om)$.
Let $\id$ denote the identity mapping.
We assume in addition that $T_t|_{\Gamma_1}=\id$ on $\Gamma_1$, where the measurements are taken.
Note that the Lagrangian $\mathcal{L}$ in \eqref{eq:elastic_lag} depends on $\Om_t$  through $A_{\Om_t}$.
The standard procedure to compute the shape derivative is to first perform a change of variables $x\mapsto T_t(x)$ in the integrals of $\mathcal{L}$, in order to transform $A_{\Om_t}$ into $A_{\Om}$.
As this change of variables induces a composition of the variables by $T_t$, it is natural to reparameterize $\mathcal{L}$ by pre-composing the functions with $T_t^{-1}$; see \cite{MR4531392} for a more detailed explanation. 
To this end, we introduce the so-called {\it shape-Lagrangian} as
\begin{align*}
	& \mathcal{G}(t,\varphi,\psi,\rho,\varrho) \\
	& \coloneqq \mathcal{L}\left(\Omega_{t},\varphi\circ T_{t}^{-1},\psi\circ T_{t}^{-1},\rho\circ T_{t}^{-1},\varrho\circ T_{t}^{-1}\right)
\end{align*}
for all $t\in[0, t_1]$, $\varphi,\psi$ in $H^1_{\Gamma_0}(\mathcal{D})^d$, and $\rho,\varrho$ in $H^1_0(\mathcal{D})^d$.
For $t\in[0, t_1]$ and $w \in H^1(\mathcal{D})^d$, let
\begin{equation*}
	E(t,w)\coloneqq\frac{Dw\,DT_{t}^{-1}+(Dw\,DT_{t}^{-1})^{\top}}{2}.
\end{equation*}
It holds that $E(t,w) = \epsilon\left(w\circ T_{t}^{-1}\right)\circ T_{t}$.
In particular, $E(0, w) = \epsilon(w)$.
On the other hand, recall that $\sigma = \mathrm{C}\epsilon$,
where $\mathrm{C}$ is the $4$-index stiffness tensor associated to
$\lambda$ and $\mu$. In terms of Kronecker symbols, it reads
\begin{equation}
	\mathrm{C}_{ijkl} \coloneqq \lambda \delta_{ij} \delta_{kl} + \mu(\delta_{ik}\delta_{jl} + \delta_{il}\delta_{jk}).
\end{equation}
Applying the change of variable $x \mapsto T_t(x)$ in $\mathcal{G}(t,\cdot)$ and using the fact that $T_t|_{\Gamma_1}=\id$ on $\Gamma_1$, we obtain
\begin{equation}
	\mathcal{G}(t,\varphi,\psi,\rho,\varrho) =
	\mathcal{G}_{\mathrm{cost}} +
	\mathcal{G}_{\mathrm{state}[f]} +
	\mathcal{G}_{\mathrm{state}[g]}
\end{equation}
with
\begin{align*}
	\mathcal{G}_{\mathrm{cost}}     &= \frac{\alpha}{2}\int_{\mathcal{D}}\left|\varphi-\psi-g\circ T_{t}\right|^{2}\xi(t)+\frac{\beta}{2}\int_{\Gamma_{1}}\left|\varphi-g\right|^{2}, \\
	\mathcal{G}_{\mathrm{state}[f]} &= \int_{\mathcal{D}}A_{\Omega}\mathrm{C}E(t,\varphi):E(t,\rho)\xi(t)-\int_{\Gamma_{1}}f\cdot\rho, \\
	\mathcal{G}_{\mathrm{state}[g]} &= \int_{\mathcal{D}}A_{\Omega}\mathrm{C}E(t,\psi + g\circ T_{t}):E(t,\varrho)\xi(t),
\end{align*}
where $\xi(t):=\det DT_t$.
Then, one can show that the shape derivative is given as the partial derivative of the shape-Lagrangian $\mathcal{G}$ with respect to $t$, using for instance the results of \cite[Chapter~10, Section~5.4]{MR2731611} or the averaged adjoint method \cite{MR3535238}.
This yields
\begin{equation}
	dJ(\Omega; \theta) = \partial_{t}\mathcal{G}(0, u, v, p, q).
\end{equation}
Using the formulas $\frac{\partial}{\partial t} DT_{t}^{-1}|_{t=0} =-D\theta$,   
\begin{align*}
	\xi'(0) & = \operatorname{div}(\theta) = I : D\theta,\\
	\mathrm{C}M:N &= M\vcentcolon\mathrm{C}N \quad \forall\,M,N\in\mathds{R}^{d\times d},\\
	\partial_t E(0,w) &= -\frac{Dw\,D\theta + (Dw\,D\theta)^\top}{2},
\end{align*}
and rearranging the various terms using tensor calculus (see for instance \cite[Lemma 2.2]{MR4531392}), we obtain
\begin{equation}\label{eq:SDelast_inv}
	dJ(\Omega; \theta) = \int_{\mathcal{D}} S_0^J \cdot \theta + S_1^J : D\theta 
\end{equation}
with $S_0^J$ and $S_1^J$ given by
\begin{align*}
	S^{J}_{0} & := -\alpha{Dg}^{\top}(u - v - g) +
	A_{\Omega}\sum_{i,j=1}^{d}{\sigma (q)}_{i,j}\nabla(\epsilon(g)_{i,j}),\\
	S^{J}_{1} & := \frac{\alpha}{2}{|u-v-g|}^2 I\\
			  & \quad + \left(A_{\Omega}\sigma(u):\epsilon(p) + A_{\Omega}\sigma(v+g):\epsilon(q)\right)I\\
			  & \quad - A_{\Omega} \left({Du}^{\top}\sigma(p) + {Dp}^{\top}\sigma(u)\right)\\
			  & \quad - A_{\Omega} \left({Dv}^{\top}\sigma(q) + {Dq}^{\top}\sigma(v+g)\right),
\end{align*}
respectively.

\textbf{Hadamard form.} Even though our numerical algorithm is based on the distributed expression~\eqref{eq:SDelast_inv}, it is instructive to compute the Hadamard form of the shape derivative, following Proposition~\ref{tensor_relations}.
Using \eqref{eq:general_boundary_exp} and using the fact that the jump of $|u-v-g|^2$ on $\partial\Om$ vanishes, we get
\begin{equation}
	dJ(\Omega;\theta) = \int_{\partial \Omega} G\, \theta\cdot n,
	\end{equation}
	with
	$G :=  \mathfrak{G}^+- \mathfrak{G}^-$ and 
\begin{align*}
	\mathfrak{G} & := A_{\Omega}\sigma(u):\epsilon(p) + A_{\Omega}\sigma(v+g):\epsilon(q)\\
			  & \quad - A_{\Omega} \left({Du}^{\top}\sigma(p) + {Dp}^{\top}\sigma(u)\right)n\cdot n\\
			  & \quad - A_{\Omega} \left({Dv}^{\top}\sigma(q) + {Dq}^{\top}\sigma(v+g)\right)n\cdot n \text{ on }\partial\Om.
\end{align*}

\section{Main algorithm}\label{sec:algo}

In this section we describe the main steps of the algorithm, which can be summarized as follows. 
We employ the level set method, introduced by Osher and Sethian~\cite{MR965860}, to update $\Om$. 
In this method, the shape $\Om$ is represented as the subzero level set of a function $\phi$:
\begin{equation} \label{eq:subdomain}
	\Omega \coloneqq \left\{x\in\mathcal{D}\mid \phi(x)<0\right\},
	\quad
	\phi:\overline{\mathcal{D}} \rightarrow \mathds{R}.
\end{equation}
For a given $\phi$, we start by solving the state and adjoint equations, which are used  to construct the shape derivative components of  \eqref{ea:volume_form2}, which were first computed analytically.
Next, a descent direction $\VV: \overline{\mathcal{D}}\to\R^d$ is computed by solving a so-called {\it velocity equation}:
Find $\theta\in\mathbb{H}$ such that
\begin{equation}\label{eq:velocity_eq0} 
	B(\theta, \xi) = - dJ(\Omega ;\xi)\qquad \forall \; \xi \in \mathbb{H},
\end{equation}	
where $dJ(\Omega ;\xi)$ is the distributed expression \eqref{ea:volume_form2} and $B:\mathbb{H}\times \mathbb{H} \rightarrow \mathds{R}$ is a user-defined positive definite bilinear form, with $\mathbb{H}$ an appropriate function space.
The solution is indeed a descent direction, as  $dJ(\Omega ;\VV) = - B(\VV,\VV)<0
$.
The implementation is described in Subsection \ref{sub:velocity}.

Further, $\VV$ is used to update
the level set function by solving the following transport-like equation on a short time interval:
\begin{equation} \label{eq:transport0}
	\begin{aligned}
	\partial_t \phi+\theta\cdot\nabla\phi &= h^2\Delta\phi && x\in \mathcal{D},\, t>0, \\
						\partial_{n} \phi &= 0 && x\in \partial\mathcal{D},\, t>0, \\
							   \phi(0, x) &= \phi^{i}(x) && x\in \overline{\mathcal{D}},
	\end{aligned}
	\end{equation}
	where $\phi^{i}$ is the current iteration.
	The implementation is described in Subsection \ref{sub:level}.

	These are the key steps of the algorithm, summarized in Algorithm~\ref{eq:algo_lv}.
The numerical implementation also involves a reinitialization, a Lagrangian approach for handling constraints, and a stopping criterion, which
will be explained in the following subsections.

\begin{algorithm}[H] 
\caption{Level set algorithm} \label{eq:algo_lv}
\small
\begin{algorithmic}[1]
\State \textbf{Initialization:} Choose initial $\phi^0:\overline{\mathcal{D}} \rightarrow \mathds{R}$.
\For{$i = 0, 1, 2, \dots$}
	\State Solve \textbf{state} problems with $\Omega^i$.
	\Statex \hspace{3em} Let $U$ be the solutions.
    \State Solve  \textbf{adjoint} problems with $\Omega^i$ and $U$.
	\Statex \hspace{3em} Let $P$ be the solutions.
	\State Use $\Omega^i$ and $U$ to compute \textbf{cost} value $J(\Om^i)$.
	\State Use $\Omega^i$ and $U$ to compute \textbf{constraint} error.
	\State Check stopping criterion based on cost value and constraint error.
	\State Use $\Omega^i$, $U$, $P$ to build derivative components:
	\Statex \hspace{1.2em} $S_0^{J}$, $S_1^{J}$ (cost) and $S_0^{C}$, $S_1^{C}$ (constraint).
	\State Solve \textbf{velocity} equation \eqref{eq:velocity_eq0} with $S_0^{J}$, $S_1^{J}$, $S_0^{C}$, $S_1^{C}$.
	\Statex \hspace{3em} Let $\theta:\overline{\mathcal{D}} \rightarrow \mathds{R}^d$ be the solution.
	\State Solve \textbf{transport} equation \eqref{eq:transport0} with $\theta$.
	\Statex \hspace{3em} Let $\phi^{i+1}:\overline{\mathcal{D}} \rightarrow \mathds{R}$ be the solution.
	\State Check stopping criterion.
\EndFor
\end{algorithmic}
\end{algorithm}

\section{Numerical implementation}\label{sec:numerics}

The Python module \code{formopt.py} was developed to apply
the level set method to shape optimization problems.
It implements the level set algorithm (see Algorithm \ref{eq:algo_lv}) using
the finite element method to solve the weak formulations of all problems involved.
The shape optimization problem is defined as a model that includes
all relevant equations and the distributed shape derivatives,
while sub-problems are addressed using parallel computing methodologies.

This section covers four aspects:
the main classes and functions provided by \code{formopt.py},
the construction of the model problem,
the three parallelism paradigms that were implemented,
and the numerical methods used to solve the sub-problems.

\subsection{Toolbox structure}
\label{sub:structure}

The main classes included in \code{formopt.py} are the following:
\begin{itemize}[left=0.0em]
\item \code{Model}:
This is an abstract class that serves as a base for user-defined classes
(i.e., subclasses of \code{Model}).
It must contain
the weak formulations of the state and adjoint equations,
the cost functional, the constraints,
the distributed shape derivative components, and
the bilinear form used to compute the velocity $\theta$, all written exclusively using Unified Form Language (UFL) and native Python functions.
To achieve this, the following methods must be overridden:
\begin{lstlisting}
	pde(level_set_func)
	adjoint(level_set_func, states)
	cost(level_set_func, states)
	constraint(level_set_func, states)
	derivative(level_set_func, states, adjoints)
	bilinear_form(velocity_func, test_func)
\end{lstlisting}
\code{Model} class provides its subclasses with
methods to create the initial level set function $\phi^0$
and to run Algorithm \ref{eq:algo_lv} using
different parallelism modes:
\begin{lstlisting}
	create_initial_level(centers, radii)
	runDP(...) # data parallelism
	runTP(...) # task parallelism
	runMP(...) # mixed parallelism
\end{lstlisting}
These methods must \textbf{not} be overridden.

\item \code{Velocity}:
This class builds and solves the velocity equation \eqref{eq:velocity_eq0}. 
Since the bilinear form $B$ remains unchanged, the \code{\_\_init\_\_} method
compiles the left-hand side of the corresponding linear system during initialization.
Thus, at each iteration, only the right-hand side is updated before solving
the linear system (see the \code{run} method of this class). 

\item \code{Level}:
This class provides the level set function at each iteration
by solving the diffusive version of the transport equation \eqref{eq:transport0}. 
The linear system is precompiled during initialization
(see the \code{\_\_init\_\_} method). 
Calling the \code{run} method updates the level set function.
Before starting the time iterations, the right-hand side is compiled,
and only the left-hand side is updated during these iterations.

\item \code{InitialLevel}:
This class creates the initial level set function $\phi^0$
from a set of centers and radii.
Its method \code{func} returns a callable function representing $\phi^0$.

\item \code{PPL}:
This class implements the Lagrangian method developed in \cite{MR4589221}
to handle the constraints; see Section~\ref{sub:lagrangian}.

\item \code{AdapTime}:
This class implements an adaptive time-stepping method to
estimate the number of steps and the final time to
solve the transport equation.

\item \code{Reinit}:
This class implements the reinitialization of the level set function
by approximating the distance function associated to $\Omega^i$.
\end{itemize}

The following are some helper functions:
\begin{itemize}[left=0.0em]
\item \code{create\_domain\_2d\_<DP|TP|MP>}:
Here, \code{DP|TP|MP} refers to the parallelism paradigm.
This function uses Gmsh to create 2D polygonal domains 
from a set of vertices, along with marked boundaries.
Interior holes and curves (formed by line segments) are supported.

\item \code{create\_space}:
This is a wrapper function for creating a \fx\ function space.

\item \code{homogeneous\_dirichlet}:
This function creates homogeneous Dirichlet boundary conditions
on faces of dimension $d-1$.

\item \code{dir\_extension\_from}:
This function computes the Dirichlet extension of the solutions
corresponding to a set of partial differential equations.
\end{itemize}

The main functions in \code{formopt.py}
implement Algorithm \ref{eq:algo_lv} using
different parallelism modes:
\begin{itemize}[left=0.0em]
\item \code{runDP} function for data parallelism,
\item \code{runTP} function for task parallelism, 
\item \code{runMP} function for mixed parallelism.
\end{itemize}

The parameters of these functions, their types, and default values are:
\begin{lstlisting}
niter: int = 100,
dfactor: float = 1e-2,
lv_time: Tuple[float, float] = (1e-3, 1e-1),
lv_iter: Tuple[int, int] = (8, 16),
smooth: bool = False,
reinit_step: int | bool = False,
reinit_pars: Tuple[int, float] = (8, 1e-1),
start_to_check: int = 30,
ctrn_tol: float = 1e-2,
lgrn_tol: float = 1e-2,
cost_tol: float = 1e-2,
prev: int = 10,
random_pars: Tuple[int, float] = (1, 0.05)
\end{lstlisting}
The \code{runMP} function has an additional parameter:
\begin{lstlisting}
sub_comm: MPI.Comm # Without default value
\end{lstlisting}
It is a MPI communicator obtained by splitting the processes in groups.
See Subsection \ref{sub:parallel} for more details.

The \code{niter} parameter is the number of iterations.
The parameters \code{dfactor}, \code{lv\_time}, \code{lv\_iter}, and \code{smooth}
are related to the level set equation.
See Subsection \ref{sub:level} and ``Adaptive time-stepping'' in Subsection \ref{sub:other} for more details.
The parameters \code{reinit\_step} and \code{reinit\_pars}
are related to the reinitialization of the level set function.
See ``Reinitialization'' in Subsection~\ref{sub:other} for more details.
The parameters \code{start\_to\_check}, \code{ctrn\_tol}, \code{lgrn\_tol}, \code{cost\_tol}, and \code{prev}, are related
to the errors and tolerances to decide when to stop the algorithm. See ``Stopping criterion'' in Subsection \ref{sub:other}.
Finally, the \code{random\_pars} parameter adds randomness to the integration of the level set equation.
See the last part of ``Adaptive time-stepping'' in Subsection \ref{sub:other} for more details.

We point out that the user-defined subclasses inherit these functions from the \code{Model} class.

In addition to \code{formopt.py}, we have implemented several models in \code{models.py}:
\begin{itemize}[left=0.0em]
\item \code{Compliance}: compliance minimization with one load;
\item \code{CompliancePlus}: similar to \code{Compliance}, but with multiple loads;
\item \code{InverseElasticity}: inverse problem in linear elasticity;
\item \code{Heat}: heat conduction problem;
\item \code{HeatPlus}: same as \code{Heat}, but with multiple sinks and sources;
\item \code{Logistic}: resource distribution for a population governed by the logistic equation.
\end{itemize}

Tests using these models are provided in \code{test.py}.
For each test, there is a corresponding folder in \code{results/},
where the results are saved. 
A tutorial covering a few of these tests can be found in \code{code/Tutorial.ipynb}.
The files \code{load.py} and \code{plots.py} contain auxiliary code
for visualizing the results. 
We recommend ParaView \cite{10.5555/2789330}
to view the results which are saved in XDMF format. 
For further details,
detailed documentation is available at
\href{https://JD26.github.io/FormOpt/}{\texttt{JD26.github.io/FormOpt/}}.

Several classes and functions are used internally, and the user does not need
to understand how to use them. Below, we describe the stages the user must follow,
along with the classes and functions that must be used:
\begin{enumerate}[left=0.0em]
\item
Creation of a model class (a subclass of \code{Model})
containing the problem equations. 
\item
Definition of the domain $\mathcal{D}$, function space, and boundary conditions.
To facilitate this step, we provide the function
\code{create\_domain\_2d\_<DP|TP|MP>} and several helper functions
for setting Dirichlet and Neumann boundary conditions. 
\item
Initialization and execution of the model
by calling the method \code{run<DP|TP|MP>}.
\end{enumerate}
Only in the last two stages does the user need to consider the parallelism paradigm.
Since the problem model is a subclass of \code{Model},
no parallelism considerations are needed during its creation.

\subsection{Model construction}
\label{sub:model}

We begin by creating a subclass of \code{Model}
with the following required attributes and methods:
\begin{lstlisting}
from formopt import Model

class MyModel(Model):

    def __init__(self, dim, domain, space, path):
        self.dim = dim
        self.domain = domain
        self.space = space
        self.path = path

	def pde(self, phi):
		pass
	
	def adjoint(self, phi, U):
		pass

	def cost(self, phi, U):
		pass
	
    def constraint(self, phi, U):
		pass

	def derivative(self, phi, U, P):
	 	pass
	
	def bilinear_form(self, th, xi):
		pass
\end{lstlisting}
The attributes \code{dim}, \code{domain}, \code{space}, and \code{path}
correspond, respectively, to the dimension $d$ of the problem domain $\mathcal{D}$,
the mesh that represents $\mathcal{D}$,
the function space for the solutions of the state and adjoint equations
(we assume the same space for both),
and the path where the results will be saved. 

The function arguments \code{phi}, \code{U}, \code{P} represent
the level set function, the state solutions, and the adjoint solutions, respectively.
The \code{pde} method defines the partial differential equations of the problem.
It must return a list with elements of the form \code{(wk, bc)},
where \code{wk} is the weak formulation of the equations,
and \code{bc} is a list of Dirichlet boundary conditions.
The \code{adjoint} method defines the adjoint equations of the problem and
must return a list with the same structure as that of \code{pde}. 
The \code{cost} method must return the cost functional $J(\Omega)$.
The \code{constraint} method must return a list with the components of $C(\Omega)$.
The \code{derivative} method must return two tuples,
each containing the derivative components $S_0$ and $S_1$,
corresponding to the cost functional and the constraints.
Finally, the \code{bilinear\_form} method must return
the chosen bilinear form $B$, and \code{th},  \code{xi}
represent the arguments $\theta, \xi$ in $B(\theta, \xi)$.

We adopt the convention that all components of this class are implemented exclusively using
UFL and native Python functions. For example, the following UFL
functions were used across all our models (see \code{models.py}):
\begin{lstlisting}
from ufl import (
    TrialFunction, TestFunction,
    FacetNormal, Identity, Measure,
    SpatialCoordinate, Coefficient,
    conditional, indices, as_vector,
    inner, outer, grad, sym, dot,
    lt, pi, cos, sqrt, nabla_div
)	
\end{lstlisting}

Details about the input and outputs of the required methods
can be found in the documentation of the \code{Model} class.

\subsection{Parallelization}\label{sub:parallel}

The main principle of parallel computing is to break a large problem into smaller tasks that can be solved independently and simultaneously using multiple processors, typically CPU cores. As the processing power of a single CPU has nearly reached its physical limits, and as data sizes increase and simulations become more detailed, computing in parallel has become essential to overcome these limitations.
The main challenges in parallel computing are how to divide the problem effectively and how to manage communication and coordination between processors. These must be handled efficiently to ensure that parallelization leads to a significant reduction in computation time.

Our module supports three parallelism paradigms:
\textit{data parallelism} (\code{DP}),
\textit{task parallelism} (\code{TP}), and a combination of the two,
known as \textit{mixed parallelism} (\code{MP}).
\textit{Data parallelism} refers to solving the problem on a mesh
that is distributed across multiple processes.
One uses domain decomposition methods to partition the domain and solve the PDE locally in each part \cite{MR2104179}. 
All the weak formulations are solved on this distributed mesh.
In the \textit{task parallelism} paradigm, each state equation is solved
in a separate process, and the same applies to each adjoint equation.
In this case, each process must have its own copy of the mesh.
The velocity field and the level set function are computed
in the first process (identified as \code{rank = 0}).
The \textit{mixed parallelism} paradigm combines the previous two by
solving each state equation on a mesh distributed across a group
of subprocesses. The adjoint equations are solved in the same 
way, while the velocity field and the level set function are computed
in the first group of subprocesses (identified as \code{color = 0}).

We implement parallelism using the Message Passing Interface (MPI). 
Communication between processes is handled through a {\it communicator}, which must be specified at the start of the program. The communicator is used in particular to broadcast data from one process to all other processes or to gather values distributed across multiple processes into a single process.
We use the default communicator \texttt{MPI.Comm\_World}.
The model problem (i.e., a subclass of \code{Model})
and its initialization are independent of the parallelism paradigm.
However, when creating the domain, the appropriate communicator must be passed.
For all paradigms, we begin by accessing the MPI communicator (\code{comm}),
the number of processes (\code{size}), and the current process (\code{rank}):
\begin{lstlisting}
from mpi4py import MPI
comm = MPI.COMM_WORLD
size = comm.size
rank = comm.rank
\end{lstlisting}

To create the domain for data parallelism, use \code{comm},
the main communicator that connect all the processes.
For 2D domains, we provide the function
\code{create\_domain\_2d\_DP}, which internally uses \code{comm}.
This function can be imported from the module.

Before starting task parallelism, we check that the number of processes matches
the number of state equations. We use the variable
\code{task\_nbr} to store this number.
After this verification, call the ``self'' communicator
(used for communication within a single process):
\begin{lstlisting}
if size != task_nbr:
    return
comm_self = MPI.COMM_SELF
\end{lstlisting}
Then use \code{comm\_self} to create the domain in each process.
Alternatively, one can import the function \code{create\_domain\_2d\_TP}
to create 2D domains. It is not necessary to pass \code{comm\_self} to this function.
Note that if adjoint equations are present,
their number is assumed to be equal to the number of state equations.
This assumption is particularly important in task parallelism,
where one process is assigned to each state equation,
and after they are solved,
the same applies to the adjoint equations (if present).

In the mixed parallelism paradigm (MP),
processes are divided into groups of equal size.
Thus, the total number of processes (\code{size}) must be divisible
by the number of groups (\code{nbr\_groups}), where \code{nbr\_groups}
corresponds to the number of state equations (and adjoint equations, if present).
Once this condition is verified, we assign a \code{color} to each process:
processes in the same group share the same \code{color}.
We expect each group to contain more than one process;
otherwise, the configuration corresponds to task parallelism.
For simplicity,
we do not consider the case where \code{size} is not divisible by \code{nbr\_groups},
thus avoiding process groups of different sizes.

The \code{Split} function returns the sub-communicator associated to each group:
\begin{lstlisting}
if size%nbr_groups != 0:
	return
color = rank * nbr_groups // size
sub_comm = comm.Split(color, rank)
\end{lstlisting}
Using \code{sub\_comm}, one copy of the same domain must be created
in each group of processes.
Alternatively, one can use the function \code{create\_domain\_2d\_MP} to do this,
passing \code{sub\_comm} as an argument.

\subsection{Lagrangian method}\label{sub:lagrangian}

To handle the constraints $C(\Omega)=\boldsymbol{1}$,
we implement the Lagrangian-based first-order method developed in \cite{MR4589221}.
We begin by formulating the Proximal-Perturbed Lagrangian
\begin{equation} \label{eq:lagrangian}
\begin{aligned}
L(\Omega, z, \lambda, \mu) \coloneqq\,
& J(\Omega) + \langle \lambda, C(\Omega) -  \boldsymbol{1} \rangle \\
&+ \langle \mu, z \rangle + \frac{\alpha}{2} \vert z \vert^2 + \frac{\beta}{2} \vert \lambda - \mu \vert^2,
\end{aligned}
\end{equation}
where $z$ is a perturbation variable,
$\lambda$ and $\mu$ are Lagrange multipliers,
$\alpha$ is a penalty parameter, and
$\beta$ is a proximal parameter.
Here, $\langle \cdot, \cdot \rangle$ and $\vert \cdot \vert$
are the standard inner product and norm of $\mathds{R}^{N_c}$.
We then follow \cite[Algorithm~1]{MR4589221}:
Given the parameters 
$r\in(0,1)$,
$\alpha\in(1,\infty)$, $\beta\in(0,1)$, 
$\delta^{0}\in(0,1]$, and initial values
$(z^{0},\lambda^{0},\mu^{0})$,
the iterations are defined by
\begin{equation} \label{eq:lag_method}
\begin{aligned}
\mu^{i} & = \mu^{i-1}+\delta^{i-1}\tfrac{\lambda^{i-1}-\mu^{i-1}}{\left\Vert \lambda^{i-1}-\mu^{i-1} \right\Vert^{2} + 1}\\
\lambda^{i} & = \mu^{i} + \tfrac{\alpha}{1+\alpha \beta} (C(\Omega^{i})- \boldsymbol{1})\\
z^{i} & = \tfrac{1}{\alpha}(\lambda^{i}-\mu^{i})\\
\delta^{i} & =r\delta^{i-1}
\end{aligned}
\end{equation}
for $i=1,2,\dots$, where $\Omega^{i}$ is given by \eqref{eq:subdomain}.
Note that, in our implementation,
the update of the level set function via the transport equation \eqref{eq:transport}
replaces the gradient descent step in \cite[Algorithm~1]{MR4589221},
which aims to minimize $L$ with respect to its first primal variable (here $\Omega$).

The Lagrange multiplier $\lambda^i$ is used to construct
the shape derivative components
$S_0$ and $S_1$ given by
\begin{equation} \label{eq:S0andS1}
S_0 \coloneqq S_0^{J} + \lambda^{i} S_0^{C}, \quad S_1 \coloneqq S_1^{J} + \lambda^{i} S_1^{C}, 	
\end{equation}
see \eqref{dJ1},\eqref{dC1}.
In the initialization of the \code{PPL} class,
we set the parameters to the same values used in \cite{MR4589221}:
$r = 0.999$, $\delta^{0} = 0.5$,
$\alpha = 2000$, $\beta = 0.5$,
$z^{0} = 0$, $\lambda^{0} = 0$, and $\mu^{0} = 0$.
These values are also used in all our tests.

Our module recognizes the number of constraints by calling
the method \code{constraints} of the user-defined class.
If there are no constraints (i.e., \code{constraints} returns an empty list),
then the Lagrangian method is not applied, and we simply set
$S_0 = S_0^J$ and $S_1 = S_1^J$.

\subsection{Velocity}\label{sub:velocity}

In view of \eqref{eq:velocity_eq0}, the velocity field $\theta$ is obtained by solving the following weak formulation:
Find $\theta\in\mathbb{H}$ such that
\begin{equation} \label{eq:velocity}
	B (\theta, \xi) =
	- \int_{\mathcal{D}} S_0\cdot\xi + S_1 : D\xi \quad\forall \; \xi \in \mathbb{H},
\end{equation}
where $S_0$ and $S_1$ are defined in \eqref{eq:S0andS1}, and
the Hilbert space $\mathbb{H}$ is either $H^1(\mathcal{D})^d$ or $H_0^1(\mathcal{D})^d$.

The \code{bilinear\_form} method must return 
the UFL expression of bilinear form $B$, along with
a Boolean value indicating whether homogeneous Dirichlet boundary conditions
should be imposed:
\code{False} if $\mathbb{H} = H^1(\mathcal{D})^d$ and
\code{True} if $\mathbb{H} = H_0^1(\mathcal{D})^d$.
In the case $\mathbb{H} = H^1(\mathcal{D})^d$,
we recommend adding to $B$ a term of the form
\begin{equation}\label{eq:bd_plt}
	{10}^4 \int_{\partial \mathcal{D}}  (\theta \cdot n) (\xi \cdot n)
\end{equation}
to penalize the normal component of $\theta$ on the boundary.
This enforces the velocity field $\theta$ to be
(almost) tangential along the boundary.
For instance:
\begin{lstlisting}
def bilinear_form(self, th, xi):
	nv = FacetNormal(self.domain)
	B = 0.1*dot(th, xi)*self.dx
	B += inner(grad(th), grad(xi))*self.dx
	B += 1e4*dot(th, nv)*dot(xi, nv)*self.ds

	return B, False	
\end{lstlisting}
corresponds to the bilinear form
\begin{equation*}
B(\theta, \xi) =  {10}^{-1} \int_{\mathcal{D}} \theta \cdot \xi + 
	\int_{\mathcal{D}} D\theta : D\xi +
	{10}^4 \int_{\partial \mathcal{D}}  (\theta \cdot n) (\xi \cdot n),	
\end{equation*}
indicating that $\mathbb{H} = H^1(\mathcal{D})^d$.

Recall that the \code{Velocity} class sets up and solve \eqref{eq:velocity} internally.
The outputs of the \code{derivative} and \code{bilinear\_form} methods
are sufficient to configure this part. 

\subsection{Transport equation}\label{sub:level}

In the standard level set method \cite{MR965860}, the level set function is updated by solving
a Hamilton-Jacobi equation on a short time interval, using a descent direction $\theta\cdot n$
in the normal direction,  derived from the boundary expression.
In the distributed shape derivative-based level set method,
one rather solves a linear transport equation, as $\theta$ is available
in $\hold$ rather than $\theta\cdot n$ on $\partial\Om$, see~\cite{MR3535238}.  
We update the level set function by solving the following diffusion--transport problem:
\begin{equation} \label{eq:transport}
	\begin{aligned}
	\partial_t \phi+\theta\cdot\nabla\phi &= h^2\Delta\phi && x\in \mathcal{D},\, t>0, \\
						\partial_{n} \phi &= 0 && x\in \partial\mathcal{D},\, t>0, \\
							   \phi(0, x) &= \phi^{i}(x) && x\in \overline{\mathcal{D}},
	\end{aligned}
	\end{equation}
where $\theta$ is the velocity field obtained from \eqref{eq:velocity}
and $h$ is the mesh diameter.
Note that a diffusion term was added to
the linear transport equation \eqref{eq:transport}.
Indeed, during the shape optimization process, small artificial interfaces
can sometimes persist between regions that are expected to merge.
These thin ``ghost boundaries'' prevent the complete merging of nearby holes,
leading to non-smooth intermediate geometries.
The additional diffusion smooths the level set evolution
and removes spurious residual interfaces.
The diffusion term in \eqref{eq:transport} can be added or removed
by setting \code{True} or \code{False}, respectively, to the parameter \code{smooth},
which is passed as an argument to the \code{run<DP|TP|MP>} functions.
By default, \code{smooth=False}.
We recommend keeping \code{smooth=False} for heat conduction problems
where the goal is to obtain tree-like material morphologies.
In contrast, for compliance minimization problems,
we suggest setting \code{smooth=True} to suppress
the formation of spurious residual interfaces.

The mesh diameter $h$ is defined as
\begin{equation}
	h \coloneqq
	\begin{cases}
		4\frac{\left|\mathcal{D}\right|}{N_T\sqrt{3}} & \text{if } d = 2,\\
		\left(6\sqrt{2}\frac{\left|\mathcal{D}\right|}{N_T}\right)^{2/3} & \text{if } d = 3,	
	\end{cases}
\end{equation}
where $\left|\mathcal{D}\right|$ denotes the area ($d=2$) or volume ($d=3$)
of the domain $\mathcal{D}$, and $N_T$ denotes the number of triangles ($d=2$) 
or tetrahedra ($d=3$).
This expression provides an approximation of the maximum diameter of the finite elements,
assuming that the elements are identical equilateral triangles or regular tetrahedra.

We apply the Petrov-Galerkin method \cite{MR1640142} to \eqref{eq:transport},
along with the Crank-Nicolson method to the time derivative.
This results in iteratively solving the following weak formulation:
Find $\phi(t+\delta t,\cdot) \in H^1(\mathcal{D})$ such that
\begin{align} \label{eq:weak_transport}
	& \int_{\mathcal{D}} \phi(t+\delta t,\cdot)\hat{\psi} + \delta t \, \mathcal{F}(\phi(t+\delta t,\cdot), \hat{\psi})  \nonumber \\
	& = \int_{\mathcal{D}} \phi(t,\cdot) \hat{\psi} - \delta t \, \mathcal{F}(\phi(t,\cdot), \hat{\psi})
	\quad \forall \, \psi\in H^1(\mathcal{D}),
\end{align}
where $\delta t$ is the time step,
$\mathcal{F}(u, v) = \frac{(\theta \cdot \nabla u) v}{2} + h^2\frac{\nabla u \cdot \nabla v}{2}$,
and $\hat{\psi}=\psi+\tau \theta\cdot\nabla\psi$, with 
\begin{equation} \label{eq:tau_lv}
	\tau(x) = \frac{1}{2}\left(\frac{1}{\delta t^2}  + \frac{|\theta(x)|^2}{h^2}\right)^{-1/2}.
\end{equation}
This choice of $\tau$ is motivated by standard practice for conservation laws,
where $\tau$ serves as a stabilization parameter that controls numerical dissipation along characteristic directions.
We solve \eqref{eq:weak_transport} up to a final time $t_{\text{end}}$ and
define the new level set function as $\phi^{i+1}(x) \coloneqq \phi(t_{\text{end}},x)$.

The Neumann boundary condition in \eqref{eq:transport} serves to eliminate
the boundary integral term of the weak form of the diffusion term.
Moreover, no inflow boundary condition is considered in \eqref{eq:transport}
since we assume that the velocity field $\theta$ is either $\theta = 0$ 
or $\theta \cdot n \approx 0$ on $\Gamma$.

The implementation and resolution of \eqref{eq:weak_transport} are handled
as in the \code{Level} class. 

\subsection{Other implementation aspects}\label{sub:other}

\textbf{Reinitialization.} 
The goal of reinitialization in the standard level set method \cite{MR965860}, see also \cite{MR3843884}, is to prevent the level set function from degenerating over successive iterations.
This means that one strives to preserve the property $|\nabla\phi|\approx 1$, at least in the vicinity of the interface.
We follow the same approach and regularly reinitialize the level set function $\phi = \phi(t_{\text{end}}, \cdot)$
(obtained from \eqref{eq:transport})
by solving the diffusive Eikonal equation with pseudo-time derivative:
\begin{equation} \label{eq:Eikonal}
\begin{aligned}
    \partial_t \varphi - h^2 \Delta \varphi &= S(\phi)(1 - |\nabla \varphi|)
    && x\in \mathcal{D},\, t>0, \\
    \partial_n \varphi &= 0
    && x\in \partial\mathcal{D},\, t>0, \\
	\varphi(0, x) &= \phi(x)
	&& x\in \overline{\mathcal{D}},
\end{aligned}
\end{equation}
where $S(\phi):\overline{\mathcal{D}}\rightarrow (-1, 1)$ is given by 
\begin{equation*}
	S(\phi) = \frac{\phi}{\sqrt{\phi^2 + h^2}}\quad (h= \text{mesh diameter}).
\end{equation*}
The function $S(\phi)$ approximates the signed distance function associated
with the set $\{x\in\overline{\mathcal{D}} \mid \phi(x)<0\}$.
We write the first equation in \eqref{eq:Eikonal}
as a Hamilton-Jacobi equation:
\begin{equation} \label{eq:HJ}
	\partial_t \varphi + H(\nabla \varphi) =
	S(\phi) + h^2 \nabla^2\varphi \quad x\in \mathcal{D},\, t>0,
\end{equation}
where $H(p) \coloneqq S(\phi) |p|$ is the Hamiltonian.
Note that $H$ is homogeneous of degree $1$: $H(\lambda p) = \lambda H(p)$
for all positive $\lambda$. By Euler's homogeneous function theorem,
$H(p) = \nabla H (p) \cdot p$, where
\begin{equation}\label{eq:grad_H}
	\nabla H (p) = S(\phi)\frac{p}{|p|}.
\end{equation}
These observations lead us to consider again
the Petrov-Galerkin method developed in \cite{MR1640142}. 
Discretizing the time derivative with the two-step Adams-Bashforth method,
we obtain the following iterative scheme:
Find $\varphi(t + \delta t, \cdot) \in H^1(\mathcal{D})$ such that
\begin{align} \label{eq:weak_HJ}
	& \int_{\mathcal{D}}  \varphi(t + \delta t, \cdot) \hat{\psi} = \int_{\mathcal{D}} (\varphi(t, \cdot) + {\delta t} \, S(\phi)) \hat{\psi} + {}\nonumber \\
	& \int_{\mathcal{D}} {\delta t}\,\mathcal{G}(\varphi(t, \cdot), \varphi(t - \delta t, \cdot), \hat{\psi}) \quad \forall \, \psi\in H^1(\mathcal{D}),
\end{align}
where $\mathcal{G}(u, v, w) = \tfrac{H(\nabla v)-3H(\nabla u)}{2} w + h^2 \tfrac{\nabla v - 3\nabla u}{2} \cdot \nabla w$
and $\hat{\psi}=\psi+\tau\nabla H(\nabla \varphi(t, \cdot)) \cdot\nabla\psi $, with 
\begin{equation}
	\tau(x) = \frac{1}{2}\left(\frac{1}{\delta t^2}  + \frac{|\nabla H(\nabla \varphi(t, x))|^2}{h^2}\right)^{-1/2}.
\end{equation}

Note that $|\nabla H(\nabla \varphi(t, \cdot))| = \left|S(\phi)\right|$ (set $p=\nabla \varphi$ in \eqref{eq:grad_H}),
so $\tau$ is time-independent, just as in \eqref{eq:tau_lv}.
The Adams-Bashforth method provides an explicit scheme
to the nonlinear equation \eqref{eq:Eikonal},
with second-order accuracy in time.
The diffusion term prevents the propagation of
small instabilities due to $\nabla H(\nabla \varphi (t, \cdot))$.
Although $\varphi$ approximates the distance function
associated to $\{x\in\overline{\mathcal{D}} \mid \phi(x)<0\}$, and
therefore ideally satisfies $\left|\nabla \varphi\right| \approx 1$,
this may not hold during the first iterations of \eqref{eq:weak_HJ}
due to the initial condition.

In practice, we observed that computing
the first iterate $\varphi(\delta t, \cdot)$
using the explicit Euler method is sufficient to start the scheme.
Finally, we point out that, unlike in \cite{MR1640142}, inflow fluxes were not considered.

Recall that the \code{Reinit} class implements the reinitialization.
Its constructor creates a precompiled solver for \eqref{eq:weak_HJ},
and the method \code{run} performs the iterations.
Two user-defined parameters are available in \code{run<DP|TP|MP>}
to configure the reinitialization:
\begin{itemize}[left=0.0em]
	\item \code{reinit\_step}:
	Either \code{False} or a positive integer.
	If \code{False}, no reinitialization is performed. Otherwise
	the reinitialization is performed whenever $\code{reinit\_step\%i=0}$, where
	\code{i} denotes the current iteration. By default, \code{reinit\_step=False}.
	\item \code{reinit\_pars}:
	Tuple with the number of iterations and the final time for \eqref{eq:weak_HJ}.
	By default, \code{reinit\_pars=(8,1e-1)}.
\end{itemize}

\vspace*{0.5cm}
\noindent\textbf{Adaptive time-stepping.} 
The final time and the number of steps/iterations in \eqref{eq:weak_transport}
are chosen as follows:
\begin{align}
    t_{\text{end}}  & = \tfrac{\code{dfactor}}{B(\theta, \theta)}\quad\text{and} \label{eq:tend} \\
    \hat{s} & =
	(s_{\max} - s_{\min}) \left( \tfrac{t_{\text{end}} - t_{\min}}{t_{\max} - t_{\min}} \right)^{1/6}
	+ s_{\min}, \label{eq:steps}
\end{align}
respectively, where
$\theta$ is the velocity field solution to \eqref{eq:velocity}.
The other parameters are user-defined:
\begin{itemize}[left=0.0em]
	\item \code{dfactor}:
	Positive float to scale the inverse of $B(\theta, \theta)$.
	We recommend choosing $\code{dfactor}\leq1$.
	By default, \code{dfactor=1e-2}.
	\item $\code{lv\_time} = (t_{\min}, t_{\max})$:
	Tuple with the minimum and maximum times allowed.
	By default, \code{lv\_time=(1e-3,1e-1)}.
	\item $\code{lv\_iter} = (s_{\min}, s_{\max})$:
	Tuple with the minimum and maximum number of steps allowed.
	By default, \code{lv\_iter=(8,16)}.
\end{itemize}
These parameters are configured in \code{run<DP|TP|MP>}.

Note that in \eqref{eq:tend} the final time $t_{\text{end}}$
varies inversely with the magnitude of $\theta$,
which carries information of the shape derivative components.
In \eqref{eq:steps}, we just apply a concave increasing function
to $t_{\text{end}}$ to find the number of steps.

To asses the robustness of the adaptive time-stepping,
randomness can be introduced by specifying
a seed number $N_\mathrm{seed}$ and a noise level $\vartheta$ in the argument
$\code{random\_pars} = (N_\mathrm{seed}, \vartheta)$:
$t_{\text{end}}$ is scaled by a random factor drawn from the uniform distribution
$\mathcal{U}(1 - \vartheta, 1 + \vartheta)$
and $\hat{s}$ is replaced by a sample drawn from the normal distribution
$\mathcal{N}(\hat{s}, (\vartheta \hat{s})^2)$.
By default, \code{random\_pars=(1, 0.05)}.

A lower bound is imposed on $B(\theta, \theta)$
to prevent division by zero in \eqref{eq:tend}.
Before returning $t_{\text{end}}$ and $\hat{s}$, we make sure that
$t_{\min}\leq t_{\text{end}} \leq t_{\max}$, $s_{\min}\leq \hat{s} \leq s_{\max}$,
and $\hat{s}$ is rounded to an integer.
We recall that the \code{AdapTime} class 
performs all these procedures.

\vspace*{0.5cm}
\noindent \textbf{Stopping criterion.}
Starting from iteration
\[
i>\code{start\_to\_check},
\]
we apply the following stopping criteria:
\begin{itemize}[left=0.0em]
	\item (problem without constraints) we check the relative error of the cost functional $J$ given by
	\[
		\left|(J(\Omega^i)-J(\Omega^j))_{j\in I(i)}\right|_{\infty}<\code{cost\_tol}\cdot |J(\Omega^i)|,
	\]
	\item (problem with constraints) we check the relative error of the Lagrangian functional $L$ given by
	\begin{align*}
		&\left|(L(\Omega^i, z^i, \lambda^i, \mu^i)-L(\Omega^j, z^j,\lambda^j,\mu^j))_{j\in I(i)}\right|_\infty\\
		&<\code{lgrn\_tol}\cdot|L(\Omega^i,z^i,\lambda^i,\mu^i)|
	\end{align*}
	and the constrain error
	\[
		|C(\Omega^i)-\textbf{1}|_\infty<\code{ctrn\_tol},
	\]
\end{itemize}
where $I(i) = \{i-\code{prev}, \dots,i-1\}$.
The parameters \code{cost\_tol}, \code{lgrn\_tol}, and \code{ctrn\_tol} 
are error tolerances, and \code{prev} is
the number of previous values considered.
The default parameter values are:
\begin{lstlisting}
start_to_check=30, cost_tol=1e-2,
lgrn_tol=1e-2, ctrn_tol=1e-2, prev=10	
\end{lstlisting}

\vspace*{0.5cm}
\noindent \textbf{Initial guess.}
The initial level set function $\phi^0$ is constructed
using two arrays:  \code{centers} and \code{radii}.
The array \code{centers} contains coordinates that
represent centers of balls included in $\mathcal{D}$,
and the array \code{radii} contains their corresponding radii.
Then, we define $\phi^0: \mathcal{\overline{D}} \rightarrow \mathds{R}$ as follows:
\begin{equation} \label{eq:phi0}
	\phi^0(x)= \code{factor} \cdot \max_k \phi_k(x),	
\end{equation}
where each $\phi_k :\mathcal{\overline{D}} \rightarrow \mathds{R}$ is given by
\begin{equation}
\phi_k(x) = \code{radii}[k] - \left|x - \code{centers}[k]\right|_{\code{ord}}. 
\end{equation}
Here, $\code{factor}$ is non-zero float number and \code{ord}
is the order of the norm $\left|\cdot\right|_{\code{ord}}$.
By default, \code{factor=1.0} and \code{ord=2} (Euclidean norm).
Note that the set
$\left\{x\in\overline{\mathcal{D}}\mid \code{factor} \cdot \phi_k(x) < 0\right\}$
determines either the complement of a closed ball
or an open ball, if \code{factor} is positive or negative, respectively.

From \eqref{eq:phi0}, and applying level set operations,
we can see that the initial domain
$\Omega^0=  \left\{x\in\overline{\mathcal{D}}\mid \phi^0(x) < 0\right\}$ is given by
\[
\Omega^0 = 
\begin{cases}
  \overline{\mathcal{D}} \setminus \bigcup_k \overline{\mathbb{B}(k)} & \text{if } \code{factor} > 0, \\
  \bigcup_k \mathbb{B}(k) & \text{if } \code{factor} < 0,
\end{cases}
\]
where $\mathbb{B}(k)$ is the $\code{ord}$-norm open ball
centered at $\code{centers}[k]$ with radius $\code{radii}[k]$.
Thus, $\Omega^0$ can be either
$\overline{\mathcal{D}}$ with ball shaped holes
or the union of balls included in $\overline{\mathcal{D}}$.

For instance,
after instantiating an object of a subclass of $\code{Model}$,
we can invoke the method \code{create\_initial\_level} to pass
the arrays that define the initial level set function $\phi^0$:
\begin{lstlisting}
md = MyModel(...)
md.create_initial_level(centers, radii)	
\end{lstlisting}
Internally, the $\code{InitialLevel}$ class
provides a callable method that uses the
\code{centers} and \code{radii} arrays to construct $\phi^0$.

\vspace*{0.5cm}
\noindent \textbf{Dirichlet extension.}
In order to compute distributed shape derivatives, one sometimes need to extend functions defined on the boundary to the entire domain $\hold$. 
This is precisely the case with the boundary measurements $g$ in Section~\ref{sec:model_prob}. 
For this purpose, we consider the Dirichlet extension of a function $u\in H^{-1/2}(\Gamma_{\text{sub}})^{d}$,
where $\Gamma_{\text{sub}}\subseteq \Gamma$ is a subset of the boundary, defined as the solution to the following problem:
Find $\eta\in H^1(\mathcal{D})^{d}$ such that
\begin{equation} \label{eq:dir}
\begin{aligned}
	\int_{\mathcal{D}} D\eta : D\zeta &= 0 &&\quad \forall\, \zeta \in H^1_{\Gamma_{\text{sub}}}(\mathcal{D})^{d}, \\
	\eta &= u &&\quad \text{on } \Gamma_{\text{sub}}.
\end{aligned}
\end{equation}
The $\code{dirichlet\_extension}$ function 
sets up and solve \eqref{eq:dir} for a list
of functions in $H^{-1/2}(\Gamma_{\text{sub}})^{d}$, and returns
a list with the corresponding Dirichlet extensions.

For instance, by solving \eqref{eq:dir}, we generate the data displacement $g$ used in \eqref{eq:w-g}.
Since the data displacement must correspond to a force application,
we provide the $\code{dir\_extension\_from}$ function, which
solves~\eqref{eq:w-f} and then extend its solution by calling
the $\code{dirichlet\_extension}$ function.

\vspace*{0.5cm}
\noindent \textbf{Penalized subdomains.}
In some applications, specific subregions of $\overline{\mathcal{D}}$ must remain fixed. 
To handle this, we have implemented the \code{Subdomain} class.
This class constructs an indicator function from a list of inequalities that
define a subdomain $\Om_0$ of $\overline{\mathcal{D}}$. We can then use this function
to add a penalty term of the form
\begin{equation}\label{eq:sub_dom}
	{10}^4 \int_{\mathcal{D}} (\theta \cdot \xi) \chi_{\Omega_0}
\end{equation}
to the bilinear form $B$.
Solving \eqref{eq:velocity} with this additional term enforces $\theta \approx 0$ within $\Omega_0$, effectively keeping the subdomain $\Omega_0$ unperturbed during the shape optimization process.
For instance, the subdomain
\begin{equation}
	\Omega_0 = \left\{ (x,y)\in \overline{\mathcal{D}} \mid 1.95 < x,\, 0.42 < y < 0.58   \right\}
\end{equation}
is modeled using the function
\begin{lstlisting}
@dib.region_of(domain)
def sub_domain(x):
    ineqs = [x[0] - 1.95, x[1] - 0.42, 0.58 - x[1]]
    return ineqs
md.sb = sub_domain.expression()
\end{lstlisting}
Each element in \code{ineqs} represents an inequality greater than zero.
It is required to include the decorator \code{region\_of} with the
current domain as its argument. This decorator is a wrapper for the \code{Subdomain} class,
turning \code{sub\_domain} into an indicator via the \code{expression} method.
Once defined, the penalty term \eqref{eq:sub_dom}
can be added to the bilinear form by including the
following line in the \code{bilinear\_form} method:
\begin{lstlisting}
def bilinear_form(self, th, xi):
	#...
	B += 1e4*self.sb*dot(th, xi)*self.dx
\end{lstlisting}

\vspace*{0.5cm}
\noindent \textbf{Computational setup.}
All numerical experiments were carried out in an \code{Anaconda} environment using \code{Python~3.11.10},
together with \code{FEniCSx~0.9.0} for the finite element computations and
\code{Gmsh~4.12.2} for mesh generation. Other important modules are
\code{NumPy~1.26.3}, \code{SciPy~1.12.0}, \code{Matplotlib~3.8.3}, and \code{MPI4Py~4.0.3}.
An installation manual is provided in the \code{README} file of the \code{GitHub} repository
(\href{https://github.com/JD26/FormOpt}{\texttt{github.com/JD26/FormOpt}}).

Most of the tests were executed on a laptop; a server was used only when explicitly indicated.
Their specifications are as follows:
\begin{itemize}[left=0.0em]
\item \textbf{Laptop:} Ubuntu 22.04.5 LTS, equipped with an Intel Core i9-13900H
processor (14 physical cores, 20 threads) and 62 GB of RAM.
\item \textbf{Server:} Debian 12 (bookworm) equipped with two AMD EPYC 9684X
processors (192 physical cores in total, 1 thread per core) and 1.5 TB of RAM.
\end{itemize}

\section{Numerical results}\label{sec:results}

\subsection{Inverse elasticity}\label{sec:num_inverse_elasticity}

We consider the inverse elasticity problem described in Section~\ref{sec:model_prob},
where the distributed shape derivative has been computed.
To account for multiple force-displacement pairs $\{(f_k, g_k)\}_k$,
we extend the cost functional \eqref{eq:cost} as follows:
\begin{equation}
	J(\Omega)\coloneqq
	\frac{\alpha}{2}\sum_{k}\int_{\mathcal{D}}|u_{k}-y_{k}|^{2}+
	\frac{\beta}{2}\sum_{k}\int_{\Gamma_{1}}|u_{k}-y_{k}|^{2},
\end{equation}
where $u_k$ is the solution to \eqref{eq:f-prob} with $f=f_k$,
and $y_k$ is the solution to \eqref{eq:g-prob} with $g=g_k$.

We start by writing the model class for this problem,
which we call \code{InverseElasticity}:
\begin{lstlisting}
from formopt import Model

class InverseElasticity(Model):
	def __init__(self, dim, domain, space, path):
		pass
\end{lstlisting}
In addition to the mandatory parameters
\code{dim}, \code{domain}, \code{space}, and \code{path},
we include the following parameters to the initializer:
\begin{lstlisting}
def __init__(self, dim, domain, space, forces, ds_forces, ds1, dirbc_partial, dirbc_total, path):
\end{lstlisting}

The parameters \code{forces} and \code{ds\_forces} are lists
containing the forces and their corresponding surfaces of application, respectively.
The parameter \code{ds1} corresponds to the surface $\Gamma_1$.
Although in \eqref{eq:f-prob} the forces are applied over the entire $\Gamma_1$,
here we assume $f\equiv\textbf{0}$ on a subset of $\Gamma_1$.
Therefore, the \code{ds\_forces} list contains only those surfaces
that are subsets of $\Gamma_1$ where $f$ is not zero.

The parameters \code{dirbc\_partial} and \code{dirbc\_total} are
the homogeneous Dirichlet conditions imposed on $\Gamma_0$ and $\Gamma$,
respectively.

In \code{\_\_init\_\_} are defined the variables and functions
that will be used to write the problem equations.
For instance, the functions $\sigma$, $\epsilon$, and $A_\Omega$
are defined as Python lambda functions:
\begin{lstlisting}
self.sigma = lambda w: (
	lm*nabla_div(w)*self.Id + 2.0*mu*self.epsilon(w)
)
self.epsilon = lambda w: sym(grad(w))
self.A = lambda w: conditional(lt(w, 0.0), 10.0, 1.0)
\end{lstlisting}
Here, \code{self.A} is an attribute of the \code{InverseElasticity} class
and represents $A_\Omega$ with $\kappa = 10$, see \eqref{eq:AOmega}.
The argument passed to \code{self.A} will be the level set function:
if negative, \code{self.A} returns $10$, else $1$.

The weak formulations for the force application \eqref{eq:w-f} and
the observed displacement \eqref{eq:w-g},
along with their boundary conditions, are returned
by the \code{f\_prob} and \code{g\_prob} methods:
\begin{lstlisting}
def f_prob(self, u, w, phi, f, df):
	# f-problem
	su = self.sigma(u)
	ew = self.epsilon(w)
	W = self.A(phi)*inner(su, ew)*self.dx
	W -= dot(f, w)*df
	return (W, self.bcF)

def g_prob(self, v, w, phi, g):
	# g-problem
	sv = self.sigma(v)
	sg = self.sigma(g)
	ew = self.epsilon(w)
	W = self.A(phi)*inner(sv, ew)*self.dx
	W += self.A(phi)*inner(sg, ew)*self.dx
	return (W, self.bcG)
\end{lstlisting}

Using these methods we write the required \code{pde} method.
It returns a list (\code{F+G}) with the weak formulation
and boundary condition corresponding
to each applied force and boundary displacement:
\begin{lstlisting}
def pde(self, phi):
	a = TrialFunction(self.space)
	b = TestFunction(self.space)
	zipf = zip(self.fs, self.dfs)
	F = [self.f_prob(a,b,phi,f,df) for f, df in zipf]
	G = [self.g_prob(a,b,phi,g) for g in self.gs]
	return F + G
\end{lstlisting}

Similarly, we write the \code{adj\_f\_prob} and \code{adj\_g\_prob} methods
to get the adjoint problems, and then we call them to collect all the adjoint equations
in the required \code{adjoint} method:  
\begin{lstlisting}
def adjoint(self, phi, U):
	a = TrialFunction(self.space)
	b = TestFunction(self.space)
	AF, AG = [], []
	for u,v,g in zip(U[:self.N],U[self.N:],self.gs):
		AF += [self.adj_f_prob(a,b,phi,u,v,g)]  
		AG += [self.adj_g_prob(a,b,phi,u,v,g)]
	return AF + AG
\end{lstlisting}
Observe that in the \code{adjoint} method we iterate over \code{U}:
the sub-lists \code{U[:self.N]} and \code{U[self.N:]} 
contain the solutions of \eqref{eq:f-prob} and \eqref{eq:g-prob}, respectively.

The cost functional \eqref{eq:cost} is returned by the \code{cost} method
(the weights \code{self.alpha} and \code{self.beta} are defined in \code{\_\_init\_\_}):
\begin{lstlisting}
def cost(self, phi, U):
	uvg = zip(U[:self.N], U[self.N:], self.gs)
	ug = zip(U[:self.N], self.gs)
	Ja = [dot(u-v-g,u-v-g)*self.dx for u,v,g in uvg]
	Jb = [dot(u-g,u-g)*self.ds1 for u,g in ug]
	return (self.alpha*sum(Ja)+self.beta*sum(Jb))/2.0
\end{lstlisting}

Since there are no constraints in this problem,
the \code{constraint} method returns an empty list:
\begin{lstlisting}
def constraint(self, phi, U):
	return []
\end{lstlisting}

In order to write the \code{derivative} method,
we first write the derivative components $S_0$ and $S_1$ separately:
\begin{lstlisting}
def S0(self, u, v, q, g, phi):
	i, j, k = indices(3)
	sq = self.sigma(q)
	eg = self.epsilon(g)
	s0a = grad(g).T*(u - v - g)
	s0b = as_vector(sq[i,j]*(grad(eg))[i,j,k], (k))
	return -self.alpha*s0a + self.A(phi)*s0b	
\end{lstlisting}
\begin{lstlisting}
def S1(self, u, v, p, q, g, phi):
	su = self.sigma(u)
	sp = self.sigma(p)
	sq = self.sigma(q)
	svg = self.sigma(v + g)
	ep = self.epsilon(p)
	eq = self.epsilon(q)
	uvg = u - v - g
	s1i = inner(su, ep) + inner(svg, eq)
	s1i = self.alpha*dot(uvg,uvg)/2.0+self.A(phi)*s1i
	s1j = grad(u).T*sp + grad(p).T*su 
	s1j += grad(v).T*sq + grad(q).T*svg
	return s1i*self.Id - self.A(phi)*s1j
\end{lstlisting}

The \code{derivative} method collect the derivative components
for each pair of force and displacement
by evaluating the \code{S0} and \code{S1} methods
at state and adjoint solutions,
then the sums \code{sum(S0)} and \code{sum(S1)} are returned: 
\begin{lstlisting}
def derivative(self, phi, U, P):
	fu = U[:self.N]
	gv = U[self.N:]
	fp = P[:self.N]
	gq = P[self.N:]
	uvqg = zip(fu, gv, gq, self.gs)
	uvpqg = zip(fu, gv, fp, gq, self.gs)
	S0 = [
		self.S0(u, v, q, g, phi)
		for u, v, q, g in uvqg
	]
	S1 = [
		self.S1(u, v, p, q, g, phi)
		for u, v, p, q, g in uvpqg
	]
	return (sum(S0), []), (sum(S1), [])
\end{lstlisting}
Empty lists are returned in the place corresponding
to the derivative components of the constraints.

Finally, we write the bilinear form used to compute
the velocity field $\theta$ by solving \eqref{eq:velocity}:
\begin{lstlisting}
def bilinear_form(self, th, xi):
    biform = 0.1*dot(th, xi)*self.dx
	biform += inner(grad(th), grad(xi))*self.dx    
	return biform, True
\end{lstlisting}
By returning \code{True} we are specifying homogeneous Dirichlet boundary condition,
that is, $\mathbb{H}=H^1_0(\Omega)$.

We are now ready to write the numerical test. 
Results using all parallelism modes are reported,
but here we describe only the test using mixed parallelism.
We start by calling the MPI communicator:
\begin{lstlisting}
comm = MPI.COMM_WORLD
rank = comm.rank
size = comm.size
\end{lstlisting}

We will consider two examples with three and eight pairs of
applied forces with their corresponding displacement observations.
Thus, there are six (resp. sixteen) state problems.
Setting \code{nbr\_groups = 6} (resp. \code{16}),
and considering \code{size=12} processes (resp. \code{32}),
the number of processes \code{size} is a multiple of \code{nbr\_groups}, and each group contains two processes.
Then, the sub-communicator can be created:
\begin{lstlisting}
color = rank * nbr_groups // size
sub_comm = comm.Split(color, rank)
\end{lstlisting}

Both examples have similar structure; the main differences lie in
the number of inclusions and the number of force--displacement pairs. 
We now describe the examples, alternating between them.

For clarity and good coding practice, we first specify
the output path for saving results,
the domain dimension ($d=2$),
the rank (i.e., the number of components of the state/adjoint solutions),
and the mesh-size parameter which will be used by \code{Gmsh} to create the mesh:
\begin{lstlisting}
test_path = Path("../results/t08/")
dim = 2
rank_dim = 2
mesh_size = 0.015
\end{lstlisting}
The above path corresponds to the first example.

The pairs of applied forces and observed displacements $\{(f_k,g_k)\}_k$
are generated on a finer mesh than that used in the examples,
by applying forces on $\Gamma_1$ and recording the corresponding displacements.
Each $g_k$ is then extended to the entire domain by solving \eqref{eq:dir}.
This procedure is performed at the beginning of the example codes
(see the \code{test\_6} and \code{test\_34} functions in \code{test.py}).

In both examples, the domain $\mathcal{D}$ is a semi-ellipse truncated at the bottom.
Below we define the vertices and the boundary parts of $\mathcal{D}$ for the second example:
\begin{minipage}{\linewidth}
\begin{lstlisting}
npts = 80
part = npts // 8

vertices = np.column_stack(semi_ellipse(0.75,0.5,0.15,npts))

# 1 dirichlet boundary and 8 neumann boundaries
dir_idx, dir_mkr = [npts], 1
neu_idxA, neu_mkrA = np.arange(1, part+1), 2
neu_idxB, neu_mkrB = np.arange(part+1, 2*part+1), 3
neu_idxC, neu_mkrC = np.arange(2*part+1, 3*part+1), 4
neu_idxD, neu_mkrD = np.arange(3*part+1, 4*part+1), 5
neu_idxE, neu_mkrE = np.arange(4*part+1, 5*part+1), 6
neu_idxF, neu_mkrF = np.arange(5*part+1, 6*part+1), 7
neu_idxG, neu_mkrG = np.arange(6*part+1, 7*part+1), 8
neu_idxH, neu_mkrH = np.arange(7*part+1, npts), 9

boundary_parts = [
	(dir_idx, dir_mkr, "dir"),
	(neu_idxA, neu_mkrA, "neuA"),
	(neu_idxB, neu_mkrB, "neuB"),
	(neu_idxC, neu_mkrC, "neuC"),
	(neu_idxD, neu_mkrD, "neuD"),
	(neu_idxE, neu_mkrE, "neuE"),
	(neu_idxF, neu_mkrF, "neuF"),
	(neu_idxG, neu_mkrG, "neuG"),
	(neu_idxH, neu_mkrH, "neuH"),
]
\end{lstlisting}
\end{minipage}

These variables and the sub-communicator are passed to
the \code{create\_domain\_2d\_MP} function:
\begin{lstlisting}
output = dib.create_domain_2d_MP(
	sub_comm, color, vertices, boundary_parts,
	mesh_size, path=test_path, plot=True
)
domain, nbr_tri, boundary_tags = output	
\end{lstlisting}
Set \code{plot=False} to hide the mesh plotted by \code{Gmsh}.

In both examples, there is only one subset of $\Gamma$ where
we impose zero displacement. In addition, it is necessary to consider the 
boundary condition of problem~\eqref{eq:w-g}.  
Then we create the space of functions and the Dirichlet boundary conditions
using predefined functions from our module:
\begin{lstlisting}
space = dib.create_space(domain, "CG", rank_dim)
dirbc_partial = dib.homogeneous_dirichlet(
	domain, space, boundary_tags, [dir_mkr], rank_dim
)
dirbc_total =
dib.homogeneus_boundary(domain, space, dim, rank_dim)
\end{lstlisting}

Boundary measures are created to impose the Neumann conditions, i.e.,
the forces applied along the boundary $\Gamma_1$:
\begin{lstlisting}
ds_parts = dib.marked_ds(
	domain, boundary_tags,
	[bR_mkr, neu_mkrA, neu_mkrB, neu_mkrC, bL_mkr]
)
ds_forces = [ds_parts[1], ds_parts[2], ds_parts[3]] 
ds1 = sum(ds_parts[1:], start = ds_parts[0])
\end{lstlisting}
The code above corresponds to the first example, where three forces are applied
on three different parts of $\Gamma_1$.
The \code{ds\_forces} list contains the boundary parts where each force is applied, 
and \code{ds1} groups these parts together.
The markers \code{bR\_mkr} and \code{bL\_mkr} correspond to segments of $\Gamma_1$ where
no forces are applied;
they are located at the right-bottom and left-bottom parts of $\Gamma_1$, respectively. 
The markers \code{neu\_mkrA}, \code{neu\_mkrB}, and \code{neu\_mkrC}
correspond to segments where the forces are applied.
In the second example, the applied forces cover the entire boundary $\Gamma_1$.

\begin{figure}[H]
\centering
		\includegraphics[width=0.48\linewidth]{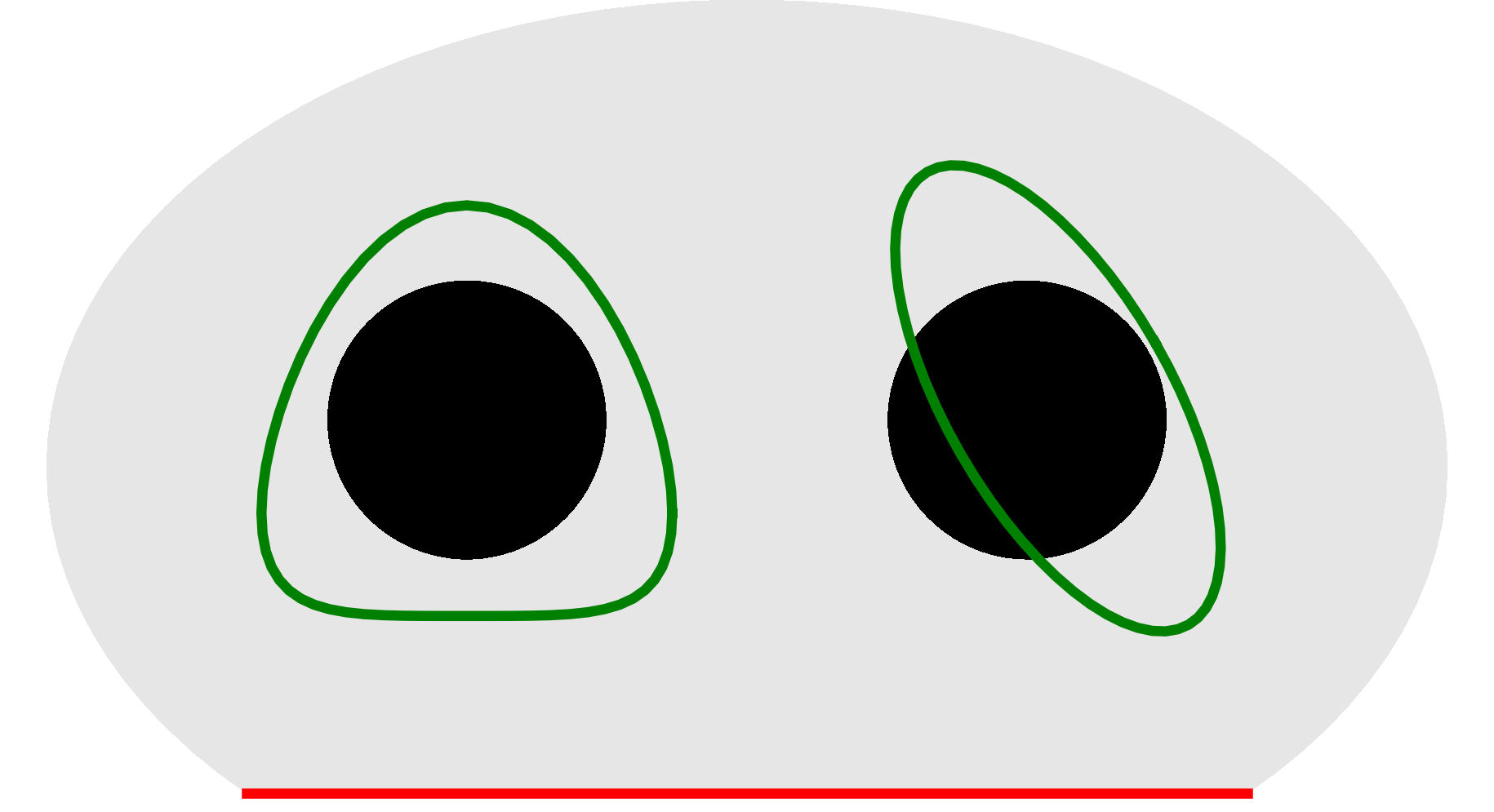}
\hfill
\includegraphics[width=0.48\linewidth]{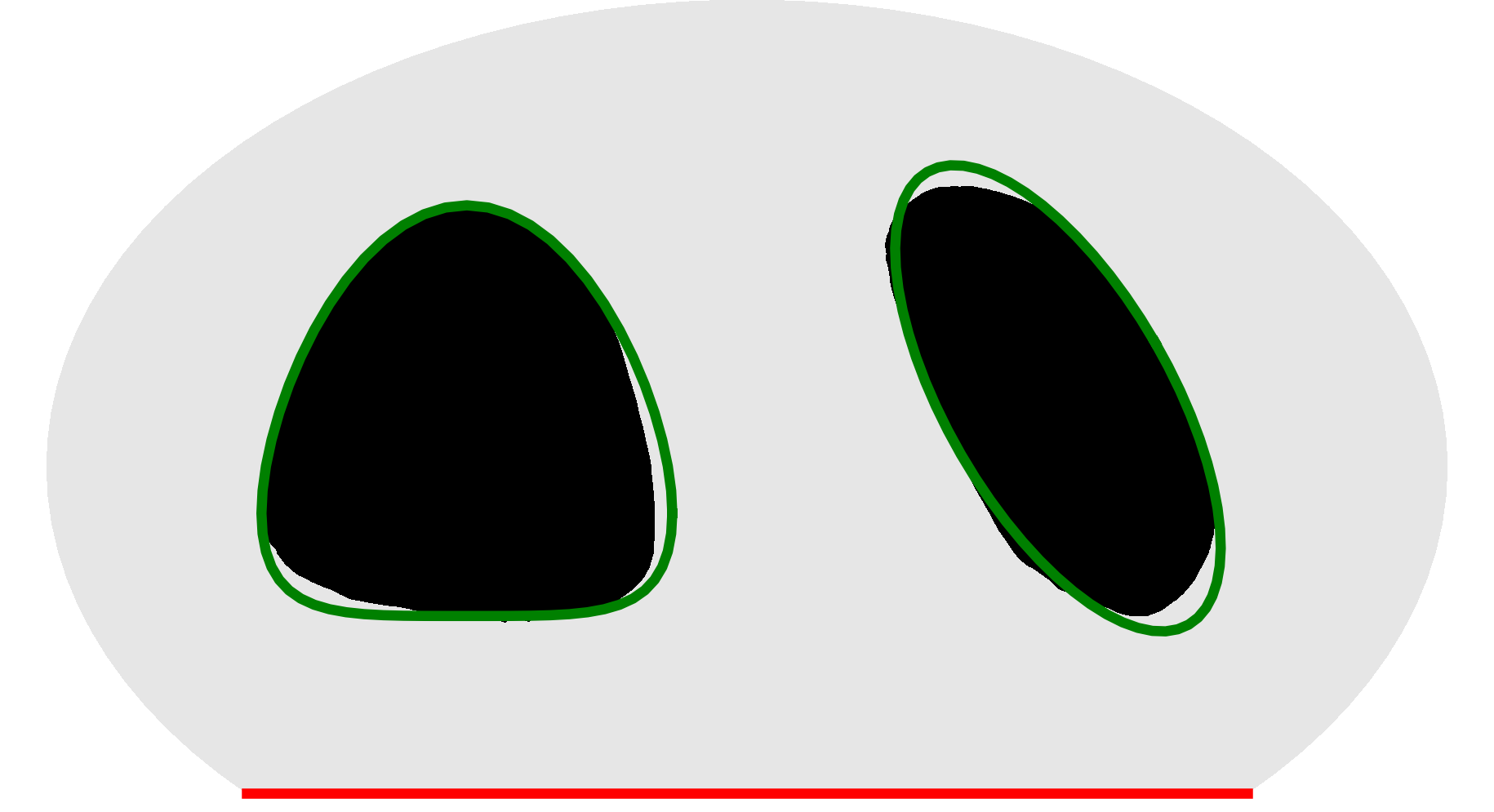}
\captionof{figure}{
	Inverse elasticity, \textit{Example 2}.
	Initial and recovered inclusions. 
	A mesh with $13\,391$ triangles was employed.
	The green line represents the ground truth.
}
\label{fig:inv2}
\end{figure}

Finally, we pass all these variables to an object of the \code{InverseElasticity} model:
\begin{lstlisting}
md = InverseElasticity(
	dim, domain, space, forces, ds_forces, ds1,
	dirbc_partial, dirbc_total, test_path
)
\end{lstlisting}
Before running the examples, we define a level set function to set
the initial inclusions in $\mathcal{D}$.
For the second example (see Fig. ~\ref{fig:inv2}), two balls of radius $0.15$ are considered as initial inclusions:
\begin{lstlisting}
centers = np.array([(-0.3, 0.4), (0.3, 0.4)])
radii = np.array([0.15, 0.15])
md.create_initial_level(centers, radii, factor=-1.0)
\end{lstlisting}
Notice that we set \code{factor=-1.0} to have the initial domain $\Omega^0$
as the union of the balls; see Subsection \ref{sub:other} for more details.

We now present the numerical results.
The zero-displacement boundaries are indicated in red, and
the true inclusion boundaries are shown in green.
In both examples, we employ the bilinear form
\begin{equation*}
	B(\theta, \xi) \coloneqq \int_{\mathcal{D}} {10}^{-1}\theta \cdot \xi + D\theta : D\xi,
	\quad
	\theta, \xi \in \mathbb{H} = H_0^1(\mathcal{D})^d,
\end{equation*}
and the following parameter values
\begin{lstlisting}
# In runDP, runTP, and runMP methods
niter=200,dfactor=0.1,lv_time=(1e-3,1.0),cost_tol=0.1
\end{lstlisting}

\begin{figure}[H]
\centering
		\includegraphics[width=0.48\linewidth]{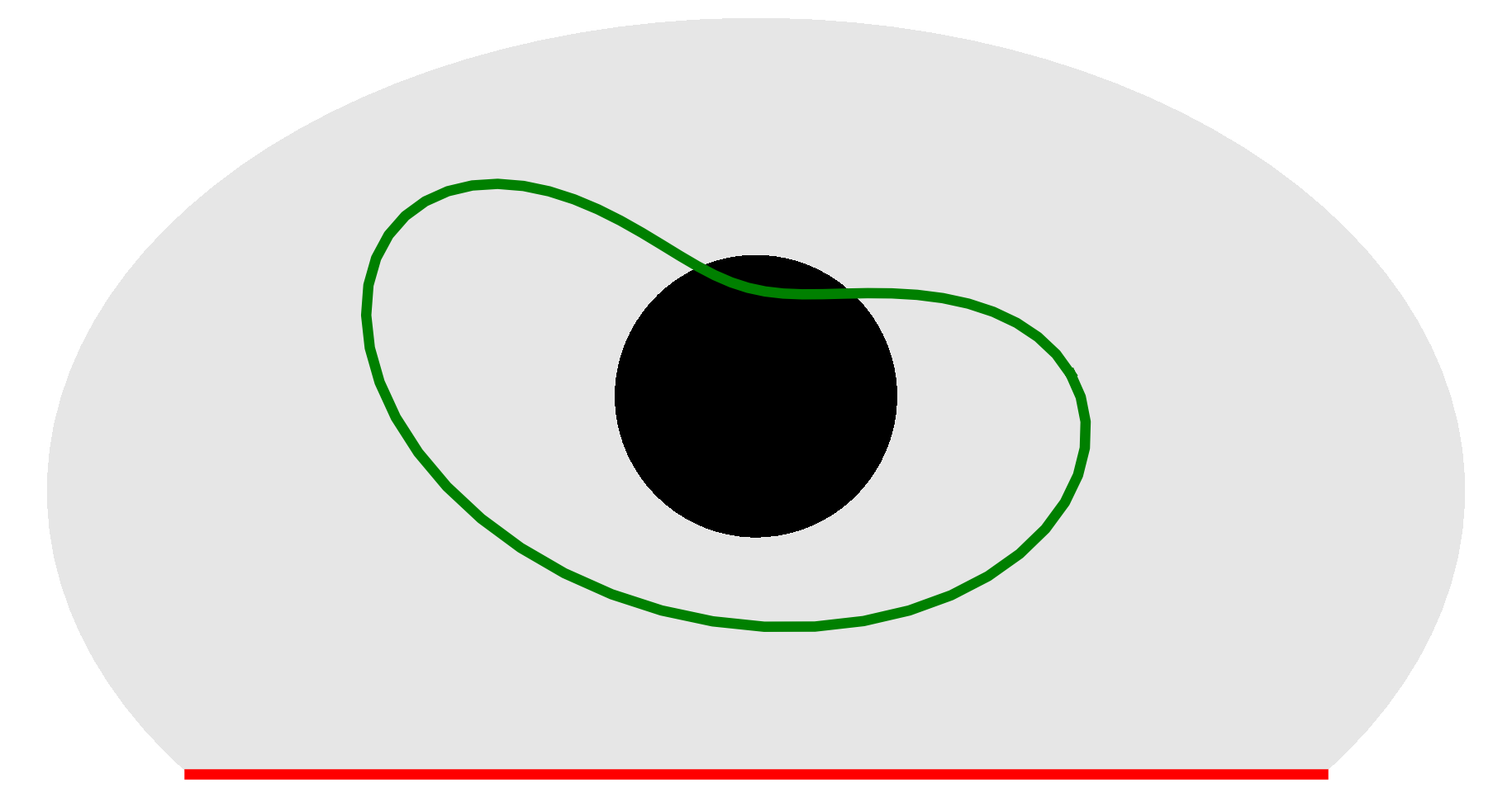}
\hfill
\includegraphics[width=0.48\linewidth]{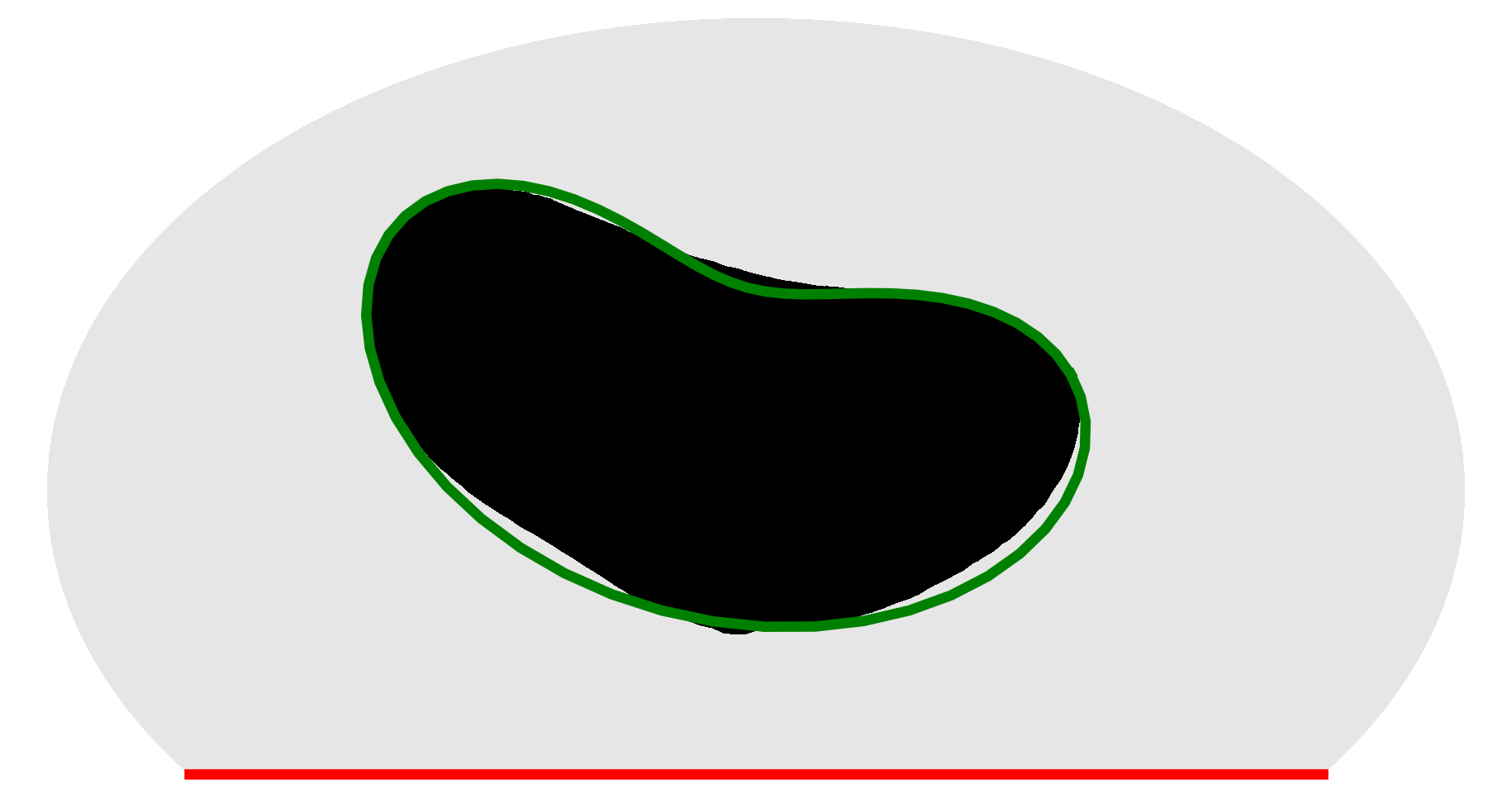}
\captionof{figure}{
	Inverse elasticity, \textit{Example 1}.
	Initial and recovered inclusion.
	A mesh with $12\,713$ triangles was employed.
	The green line represents the ground truth.
}
\label{fig:inv1}
	\end{figure}

\noindent\textit{Example 1.} (Recovery of a single inclusion)
The data consists of three pairs of force $f_k$ and displacement $g_k$.
The forces are applied separately on different parts of $\Gamma_1$,
all directed toward the center of the body.
The resulting boundary displacements are measured, only in $\Gamma_1$
since $\Gamma_0$ is fixed.
Thus, we have six state equations and their corresponding adjoints.
We conducted three separate experiments: the first using data parallelism with $6$ processes,
the second using task parallelism with $6$ processes, and the third using mixed parallelism with $12$ processes
(organized into $6$ groups of $2$ processes each). The commands and corresponding execution times were:
\begin{lstlisting}
mpirun -np 6 python test.py 06 # Data (33 sec)
mpirun -np 6 python test.py 07 # Task (29 sec)
mpirun -np 12 python test.py 08 # Mixed (27 sec)
\end{lstlisting}
The number of iterations for these runs was $124$, $114$, and $116$, respectively.
See the result in Figure~\ref{fig:inv1}.

\noindent\textit{Example 2.} (Recovery of two inclusions)
The data consists of eight pairs of force $f_k$ and displacement $g_k$.
The forces are applied separately on differents parts of $\Gamma_1$,
all directed toward the center of the body,
and covering the whole boundary $\Gamma_1$;
then resulting boundary displacements are measured.
Thus, we have sixteen state equations and their corresponding adjoints.
We conducted three separate experiments: the first using data parallelism with $16$ processes,
the second using task parallelism with $16$ processes,
and the third using mixed parallelism with $32$ processes
(organized into $16$ groups of $2$ processes each).
The commands and corresponding execution times were:
\begin{lstlisting}
mpirun -np 16 python test.py 34 # Data (146 sec)
mpirun -np 16 python test.py 35 # Task (660 sec)
mpirun -np 32 python test.py 36 # Mixed (43 sec)
\end{lstlisting}
The number of iterations for these runs was $174$, $173$, and $174$, respectively.
These tests were carried out on the server.
See the recovered inclusions in Figure \ref{fig:inv2}.

\subsection{Compliance minimization}\label{sub:compli}

Compliance minimization is a standard problem in structural mechanics, for which an abundant literature exists, see  \cite{MR1911658,MR2033390,Sigmund2014,vanDijk2013,MR1951408,MR2108636}.
Compliance minimization using the distributed shape derivative is also the main topic of the previous work \cite{MR3843884}.
We use it here to illustrate the performance of \fopt on benchmark problems.
Instead of working directly with the variable domain $\Omega$, it is  convenient to reformulate the problem on the fixed domain $\hold$. To this end, we follow the standard approach of approximating the original problem by filling the complement $\overline{\mathcal{D}} \setminus \Omega$ with an ``ersatz material''.

We minimize the compliance in a linear elasticity problem:
\begin{equation} \label{eq:compli_cost}
	\min_{\Omega \subset \mathcal{D}}
	J(\Omega) \coloneqq \int_{\mathcal{D}} A_{\Omega}\sigma(u):\epsilon(u)
\end{equation}
subject to
\begin{equation}\label{eq:compli_constr}
	\int_{\mathcal{D}} \chi_{\Omega} = V,
\end{equation}
where $u$ is the solution to the variational formulation:
Find $u \in  H^1_{\Gamma_0}(\mathcal{D})^d$ such that
\begin{equation}\label{eq:compli_weak}
\int_{\mathcal{D}} A_{\Omega}\sigma(u):\epsilon(v) = \int_{\Gamma_1} g \cdot v
\qquad
\forall \, v \in  H^1_{\Gamma_0}(\mathcal{D})^d,
\end{equation}
and $A_{\Omega}:\overline{\mathcal{D}} \rightarrow \mathds{R}$ is given by
\begin{equation}\label{eq:compli_AOmega}
	A_{\Omega} = \chi_{\Omega} + {10}^{-4}\chi_{\overline{\mathcal{D}}\setminus\Omega}.
\end{equation}
The corresponding strong formulation is the following transmission problem:
\begin{equation}\label{eq:compli_prob}
	\begin{aligned}
	   -\operatorname{div} A_{\Omega} \sigma(u) &= 0 && \text{in } \Omega \text{ and }\overline{\mathcal{D}}\setminus\Omega\\
											 u  &= 0 && \text{on } \Gamma_{0}\\
					     A_{\Omega} \sigma(u)n  &= g && \text{on } \Gamma_{1} \\
					     A_{\Omega} \sigma(u)n  &= 0 && \text{on } \Gamma \setminus (\Gamma_{0} \cup \Gamma_{1})	\\
				      (A_{\Omega} \sigma(u)n)^+ &= (A_{\Omega} \sigma(u)n)^- && \text{on }  \partial\Omega\\      
						                    u^+ &= u^- && \text{on } \partial\Omega   	.				 
	\end{aligned}
\end{equation}

The solution to the adjoint problem is $p=-2u$.
From \eqref{eq:compli_constr}, the constraint function is given by
\begin{equation*}
	C(\Omega) = \frac{1}{V} \int_{\mathcal{D}} \chi_{\Omega}.
\end{equation*} 
The derivative components of the compliance $J(\Omega)$ and the constraint function $C(\Omega)$ are
\begin{align*}
	S^J_0 &= \boldsymbol{0} ,\quad S^J_1 = A_{\Omega} (2{Du}^{\top}\sigma(u)- \sigma(u):\epsilon(u))I,\\
	S^C_0 &= \boldsymbol{0},\quad S^C_1 = \frac{1}{V} \chi_\Omega I.
\end{align*}

The \code{Compliance} and \code{CompliancePlus} classes
contain the equations for compliance minimization.
\code{CompliancePlus} extends \code{Compliance} by considering
multiple forces $\left\{g_k\right\}_k$, along with
the cost functional $J(\Omega)$ as a sum of compliances associated to each $g_k$, namely
\begin{equation} \label{eq:compli_plus}
	J(\Omega) \coloneqq \sum_{k} \int_{\mathcal{D}} A_{\Omega}\sigma(u_k):\epsilon(u_k),
\end{equation}
where $u_k$ is the solution to \eqref{eq:compli_weak} with $g=g_k$.

In the following examples, the region with zero displacement and
the force application areas are shown in red and blue, respectively.

\begin{figure}[H]
\centering
		\includegraphics[width=0.48\linewidth]{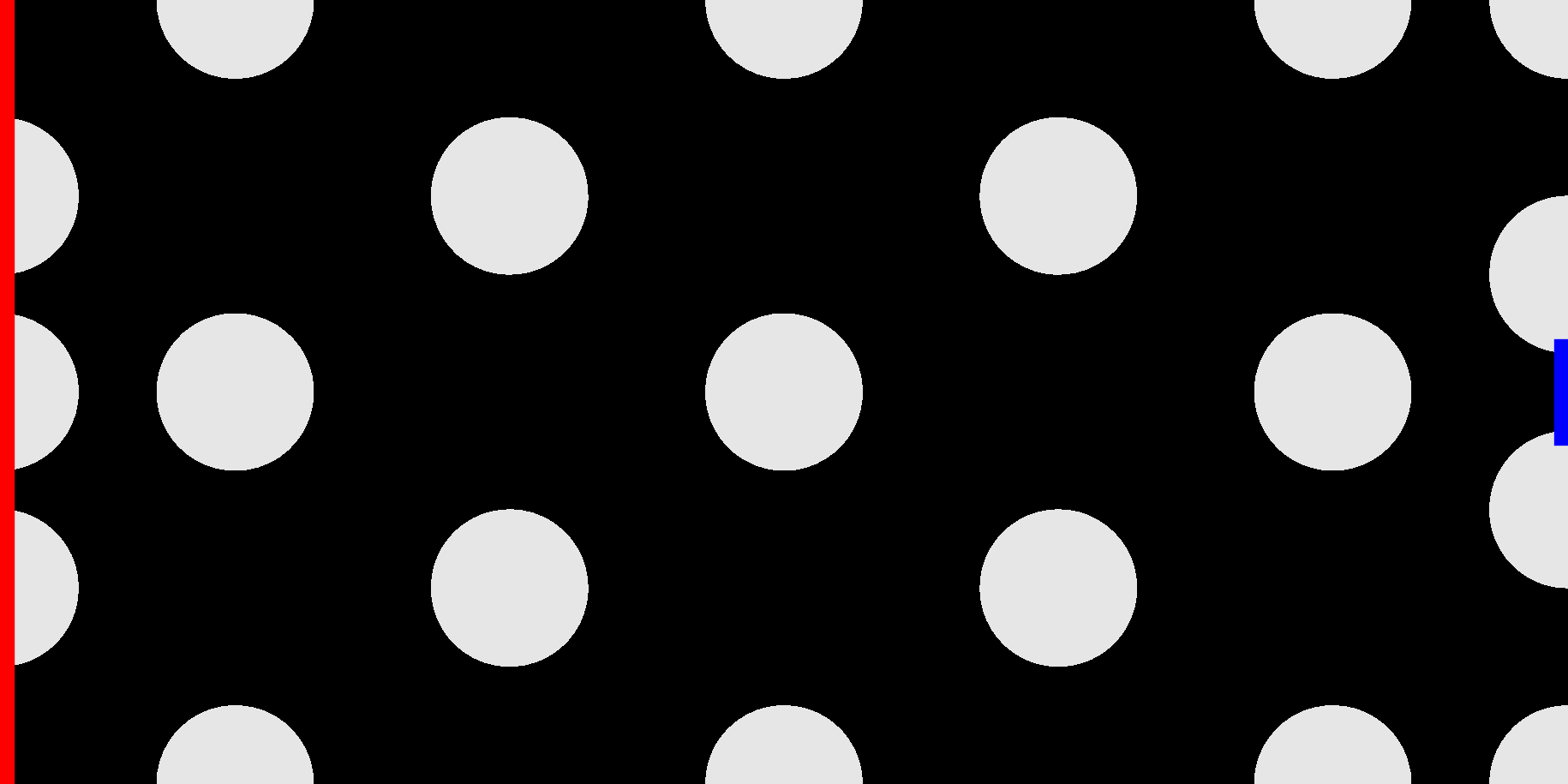}
\hfill
\includegraphics[width=0.48\linewidth]{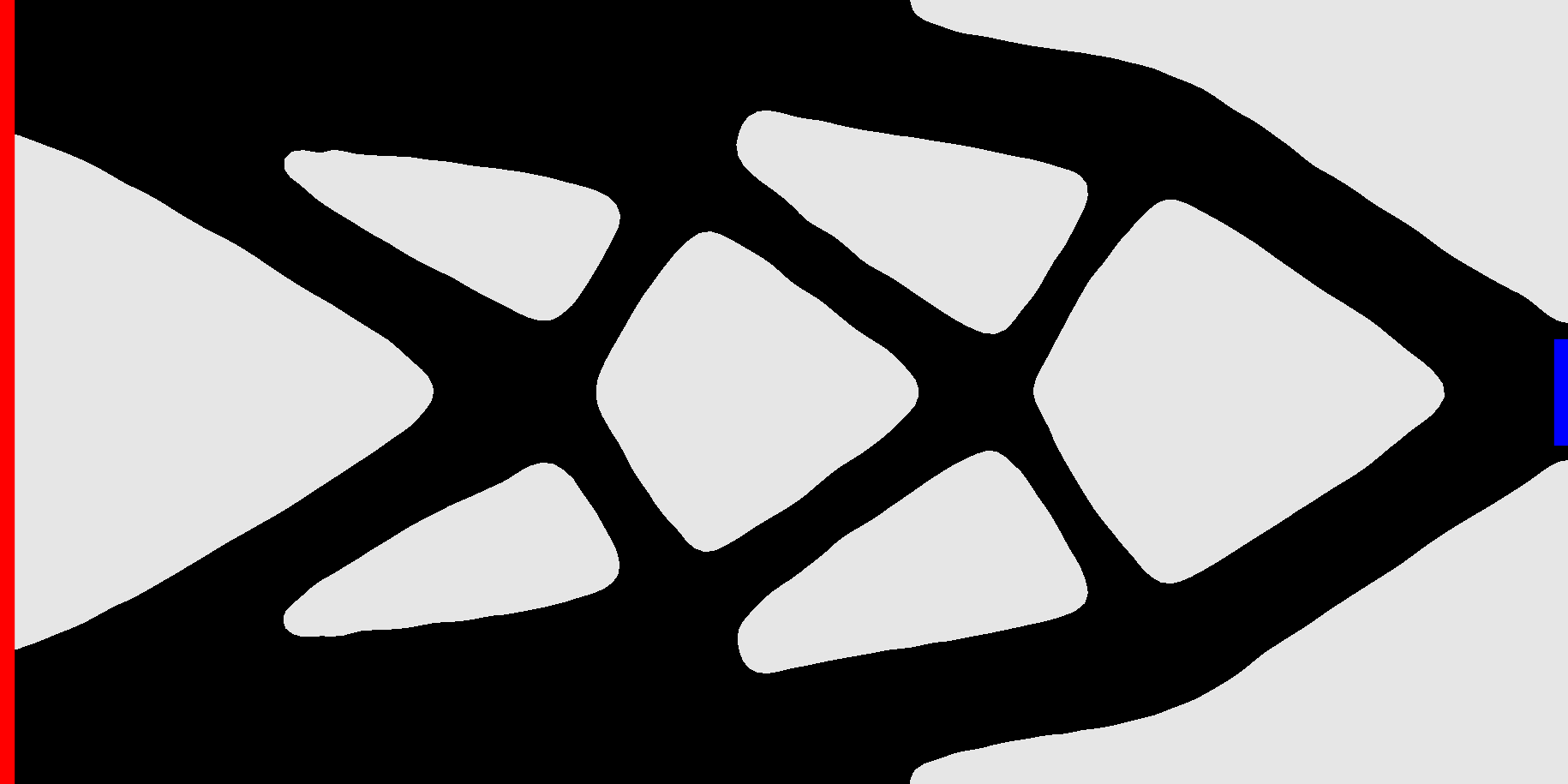}
\includegraphics[width=0.85\linewidth]{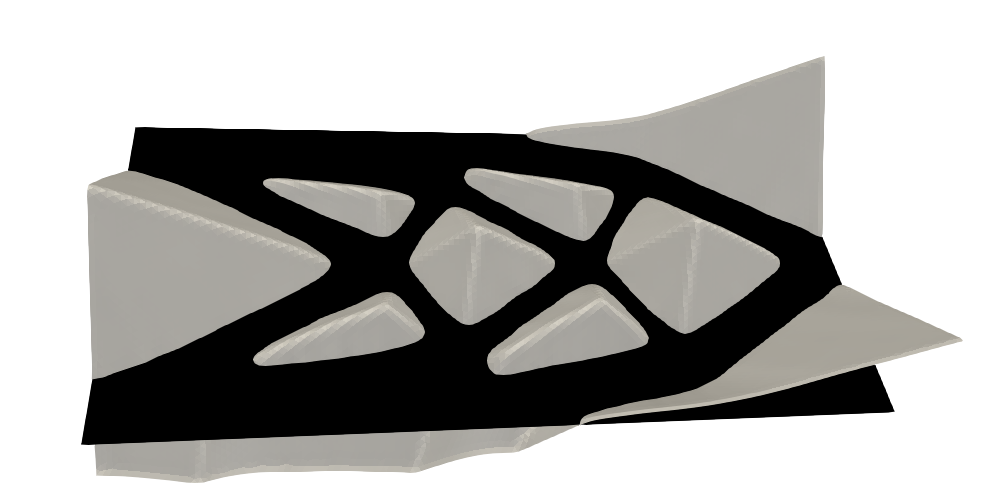}
\captionof{figure}{
Compliance minimization, \textit{Example 1}.
Initial guess (top-left) and optimized design (top-right) at iteration $i=50$.
The resulting level set function $\phi^{i}$ approximates the
distance function associated to $\Omega^{i}$ (bottom).
The domain $\mathcal{D}$ was discretized with $20\,946$ triangles.
$\Gamma_0$ appears in red and $\Gamma_1$ in blue.
}
\label{fig:compli1}
\end{figure}

\noindent\textit{Example 1.} (Symmetric cantilever)
We consider the rectangular domain $\mathcal{D} = (0, 2)\times(0, 1)$
with boundary subsets $\Gamma_0 = \{0\} \times (0, 1)$
and $\Gamma_1 = \{2\} \times (0.45,0.55)$.
The vertical force $g = (0, -2)^\top$ is applied on $\Gamma_1$ and
the material area is constrained to $V = 1$.
We run this example using data parallelism with four processes:
\begin{lstlisting}
mpirun -np 4 python test.py 01 # Data (6 sec)
\end{lstlisting}
See the results in Figure \ref{fig:compli1}. 
We have performed a reinitialization of the level set function
every four steps, with $20$ iterations and a final time of $0.1$.
This yields a well approximated distance function.
The other parameters passed to the \code{runDP} method
have the following values (with default values used for all unspecified parameters): 
\begin{lstlisting}
md.runDP(
	ctrn_tol=1e-3, dfactor=1e-1,
	reinit_step=4, reinit_pars=(20, 0.1), smooth=True
)
\end{lstlisting} 
Note that smoothness of the level set function was enforced by setting \code{smooth=True}.

The bilinear form used in this example is
\begin{align}
\notag B(\theta, \xi) &\coloneqq \int_{\mathcal{D}} {10}^{-1}\theta \cdot \xi + D\theta :D\xi+{10}^4(\theta\cdot\xi)\chi_{\Omega_0}\\
\label{eq:compli1_bili}		   &+ \int_{\partial\mathcal{D}} {10}^4 (\theta \cdot n)(\xi \cdot n),
\end{align}	
with $\theta, \xi \in \mathbb{H} = H^1(\mathcal{D})^d$. Thus, $B$
allows only tangential displacements along the boundary $\partial \mathcal{D}$, and
penalizes the material around the force application boundary through the term ${10}^4(\theta\cdot\xi)\chi_{\Omega_0}$, where
the subdomain
\[
\Omega_0 \coloneqq (1.95,2]\times (0.42, 0.58) \subset \overline{\mathcal{D}}
\]
contains the boundary $\Gamma_1$.

\noindent\textit{Example 2.} (Three-dimensional symmetric cantilever)
We consider the rectangular box domain
$\mathcal{D} = (0, 2)\times(0, 1)\times(0, 1)$
with boundary subsets
$\Gamma_0 = \{0\} \times (0,1) \times (0,1)$
and
$\Gamma_1 = \{1\} \times (0.4,0.6) \times (0.4,0.6)$.
The force $g = (0, 0, -4)^\top$ is applied on $\Gamma_1$, and
volume material is constrained to $V = 1$.
We run this example using data parallelism with $1$, $2$, and $3$ processes:
\begin{minipage}{\linewidth}
\begin{lstlisting}
mpirun -np 1 python test.py 02 # (4.357 hours)
mpirun -np 2 python test.py 02 # (2.332 hours)
mpirun -np 4 python test.py 02 # (1.486 hours)
\end{lstlisting}
\end{minipage}
See some iterations in Figure \ref{fig:compli2}.
We employ the bilinear form \eqref{eq:compli1_bili} with
$\Omega_0 = (1.90,2]\times (0.35, 0.65) \times (0.35, 0.65)$
in order to penalize the material around $\Gamma_1$.
The level set function used as initial guess was constructed
with infinity-norm balls:
\begin{lstlisting}
md.create_initial_level(centers, radii, ord=np.inf)
\end{lstlisting} 
This choice is consistent with the the regular cubic lattice of
the mesh generated by the \code{create\_box} function of \code{FEniCSX}.
The parameters passed to the \code{runDP} method
have the following values:
\begin{lstlisting}
md.runDP(
	ctrn_tol=1e-3, dfactor=1e-1,
	reinit_step=4, reinit_pars=(4,1e-2), smooth=True
)
\end{lstlisting} 
\begin{figure}[H]
\centering
		\includegraphics[width=0.48\linewidth]{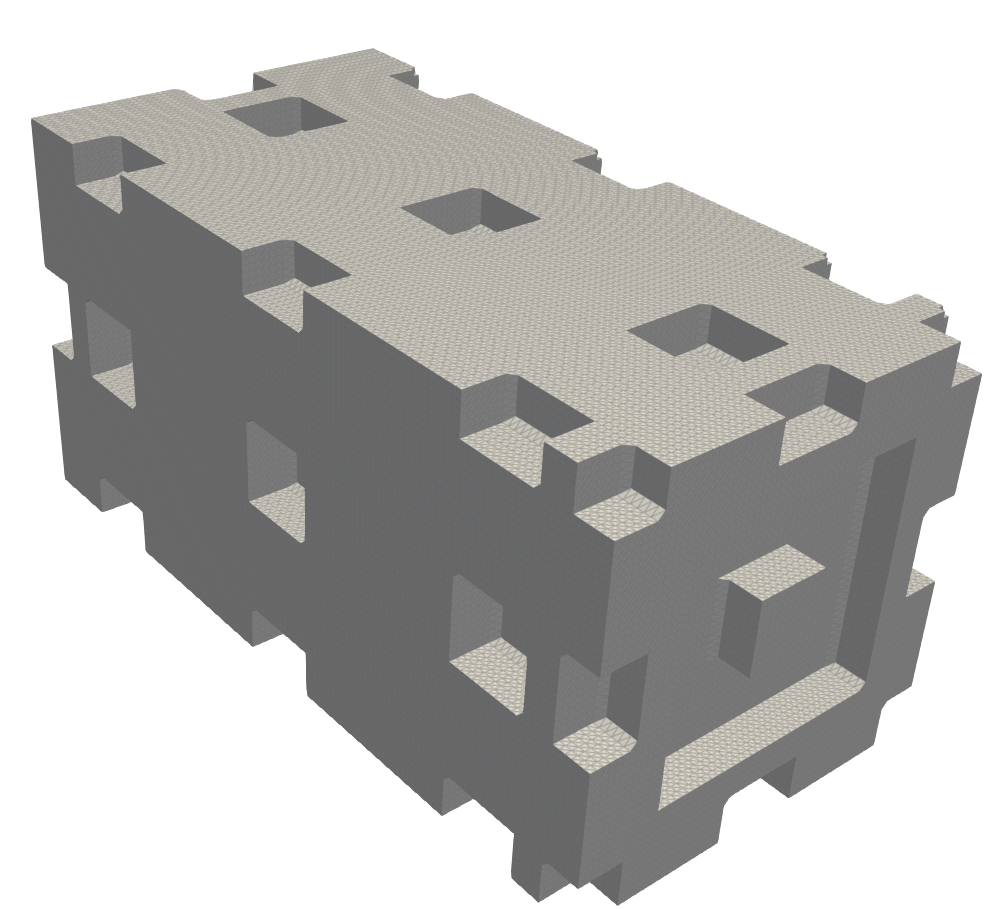}
\hfill
\includegraphics[width=0.48\linewidth]{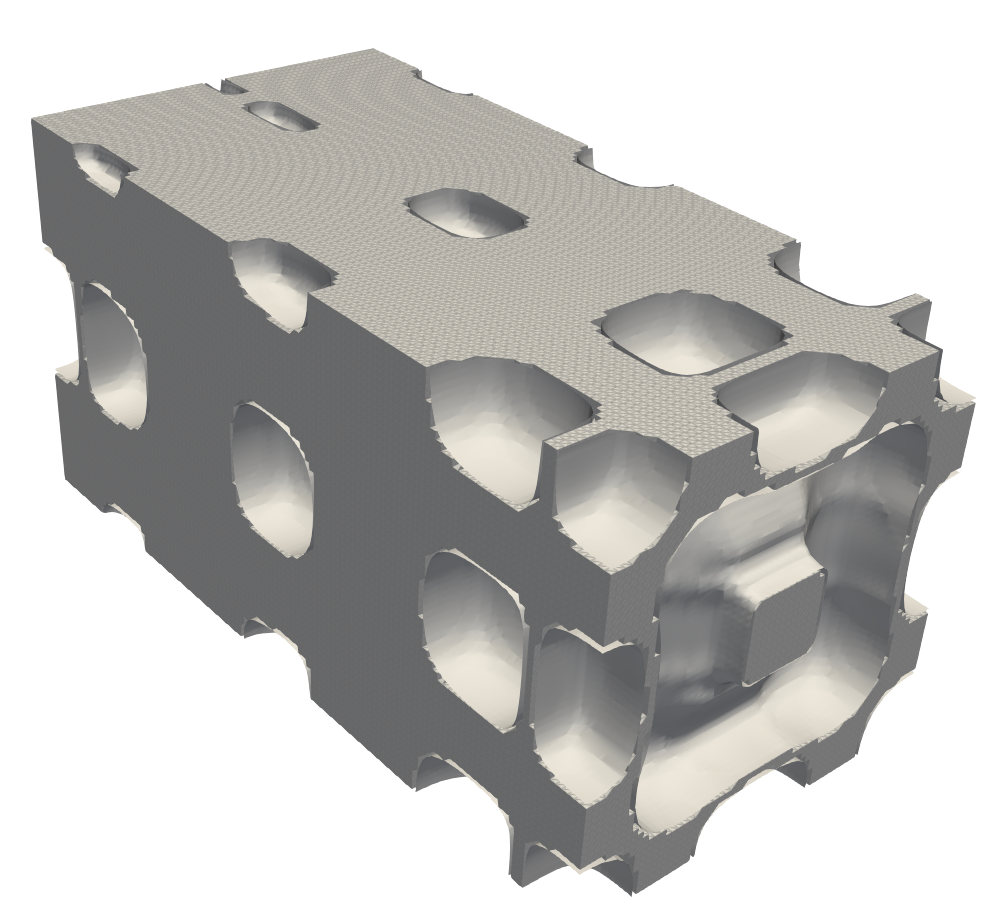}
\includegraphics[width=0.48\linewidth]{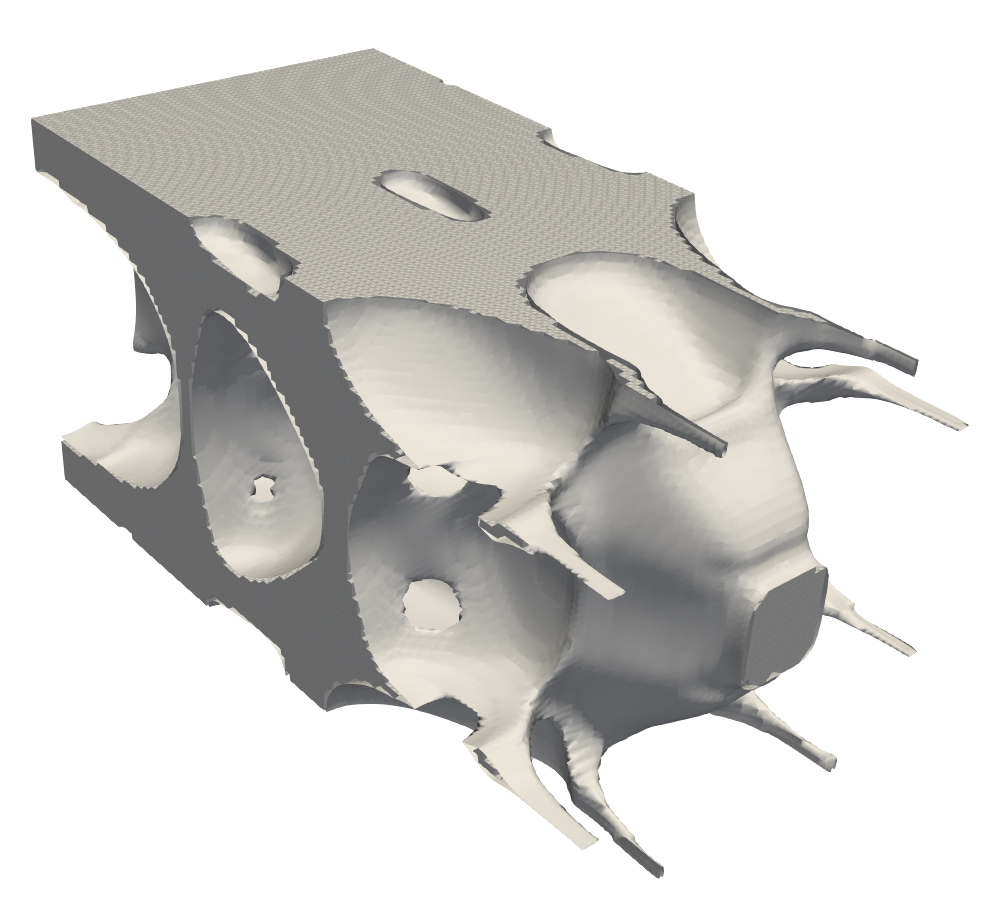}
\hfill
\includegraphics[width=0.48\linewidth]{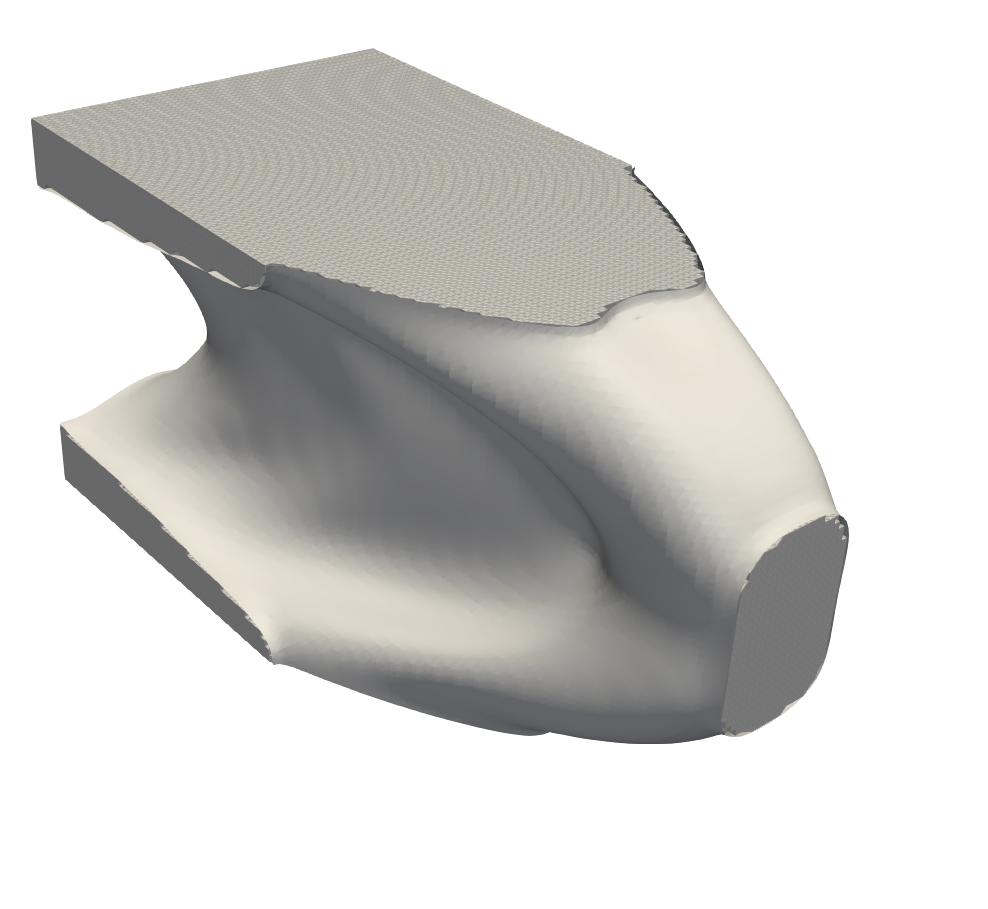}
\captionof{figure}{
Compliance minimization, \textit{Example 2}.
The three-dimensional cantilever at	
$i=0$ and $i=14$ (top);
$i=30$ and $i=81$ (bottom).
A mesh with $2\,592\,000$ tetrahedra was employed.
}
\label{fig:compli2}
\end{figure}

\noindent\textit{Example 3.} (Cantilever with two loads I)
Here $\mathcal{D}$ is the unit square and 
we solve the problem with multiple forces $g^{\mathrm{I}} = g^{\mathrm{II}} = (0, -2)^\top$,
applied on $\Gamma_1^{\mathrm{I}} = \{x=1\} \times (0,0.1)$ (right-upper side)
and $\Gamma_1^{\mathrm{II}} = \{x=1\} \times (0.9,1)$ (right-bottom side),
respectively.
The boundary of zero displacement is $\Gamma_0 = \{x=0\} \times (0, 1)$,
and the material area is constrained to $V = 0.5$.
Thus, we minimize the cost functional \eqref{eq:compli_plus},
which is implemented in the \code{cost} method of the \code{CompliancePlus} class.
To compare performances, we conducted three experiments: the first using data parallelism with $2$ processes, the second using
task parallelism with $2$ processes, and the third using
mixed parallelism with $4$ processes
(organized into $2$ groups of $2$ processes each).
The commands and corresponding execution times were:
\begin{lstlisting}
mpirun -np 2 python test.py 03 # Data (20 sec)
mpirun -np 2 python test.py 04 # Task (26 sec)	
mpirun -np 4 python test.py 05 # Mixed (15 sec)		
\end{lstlisting}
The optimized design and the number of iterations are the same in all cases,
see Figure \ref{fig:compli3}.
We have used the bilinear form \eqref{eq:compli1_bili} without any penalized subdomain.
The configuration passed to the \code{runDP}, \code{runTP}, and \code{runMP} methods was:
\begin{lstlisting}
ctrn_tol=1e-3, dfactor=1e-1,
reinit_step=4, reinit_pars=(16, 0.05), smooth=True
\end{lstlisting}
\begin{figure}[H]
	\centering
\includegraphics[width=0.48\linewidth]{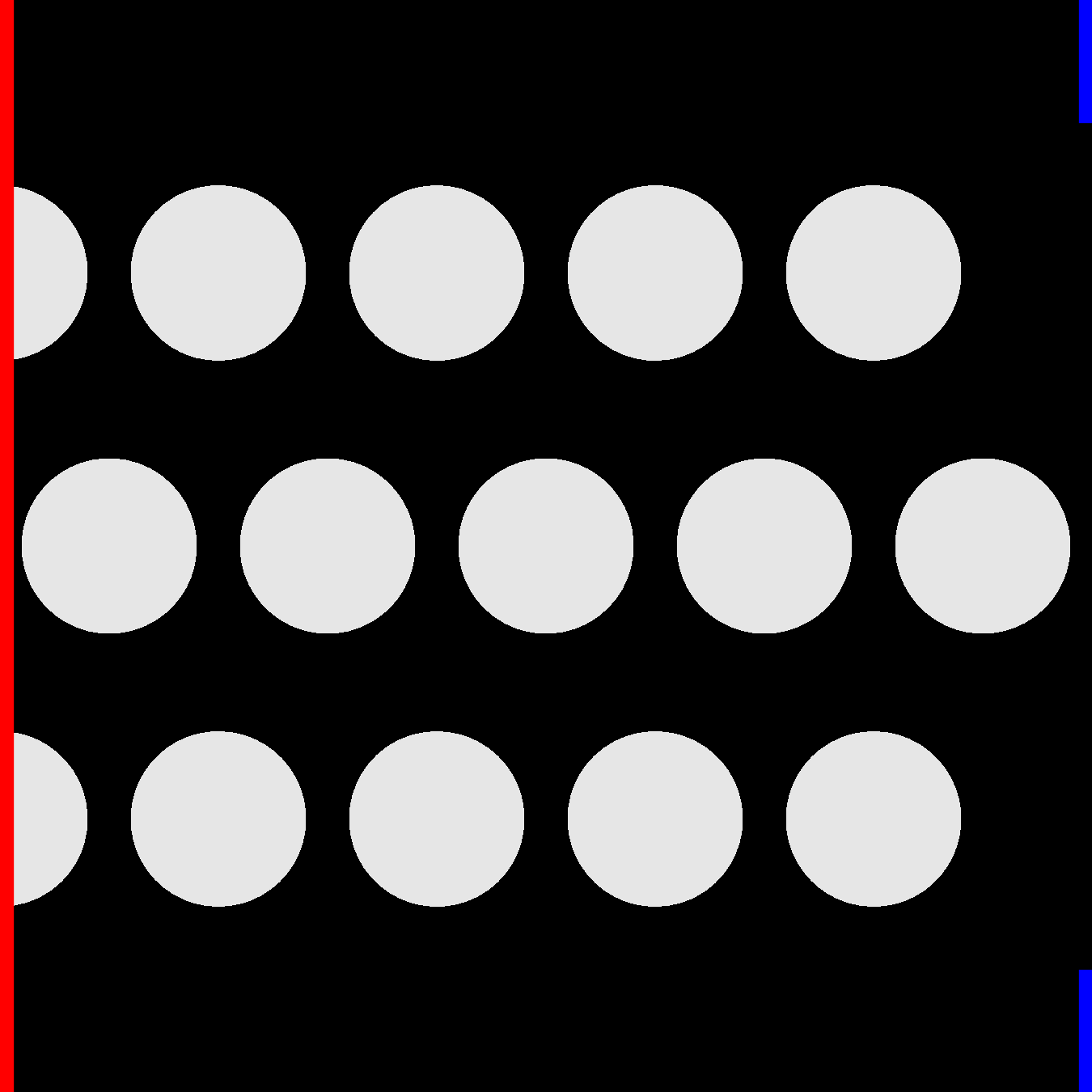}
\hfill
\includegraphics[width=0.48\linewidth]{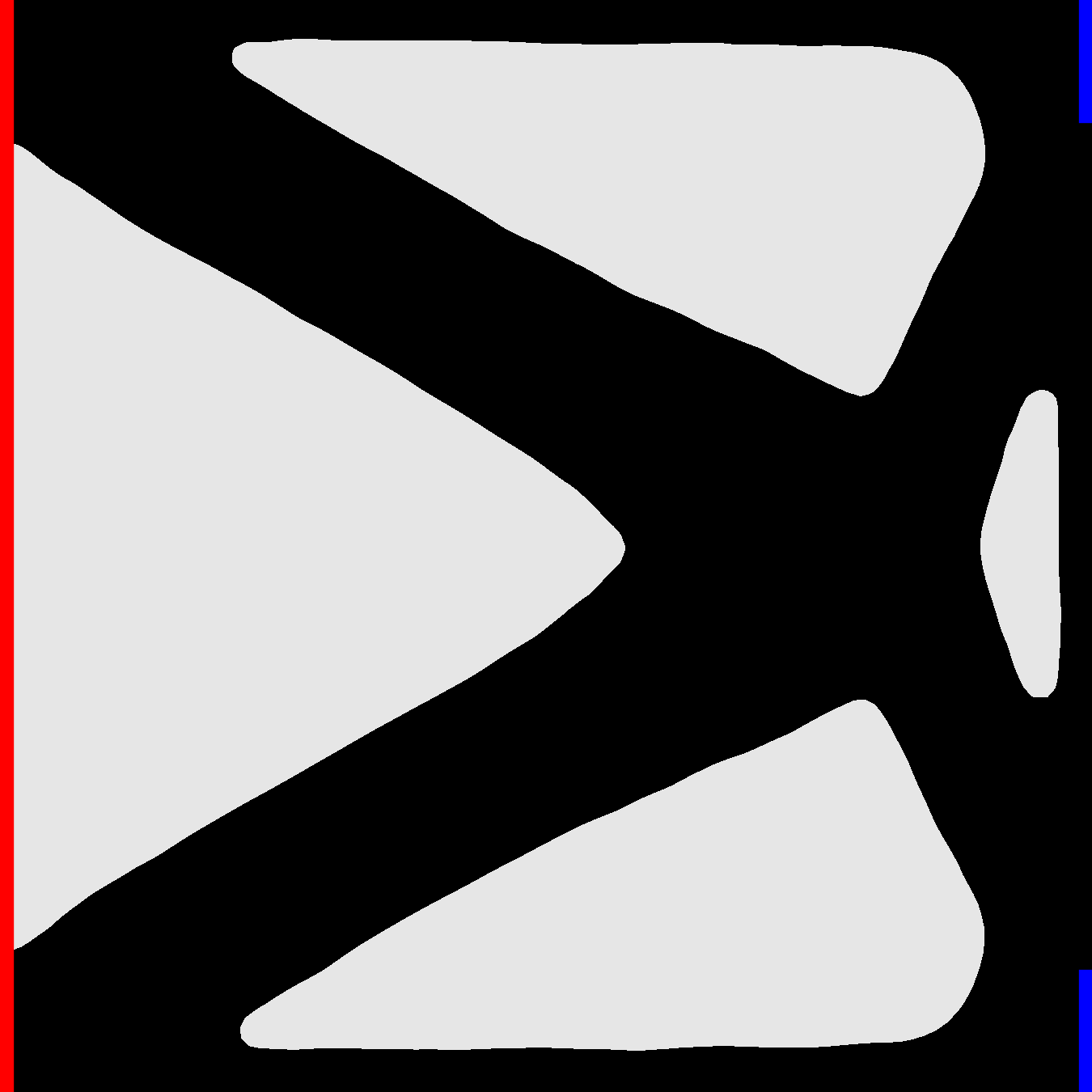}
\captionof{figure}{
	Compliance minimization, \textit{Example 3}.
	Initial and optimal domains.
	A total of $70$ iterations were perfomed for the three modes of parallelism.
	A mesh with $16\,425$ triangles was employed.
}
\label{fig:compli3}
\end{figure}

\begin{figure}[H]
	\centering
\includegraphics[width=0.48\linewidth]{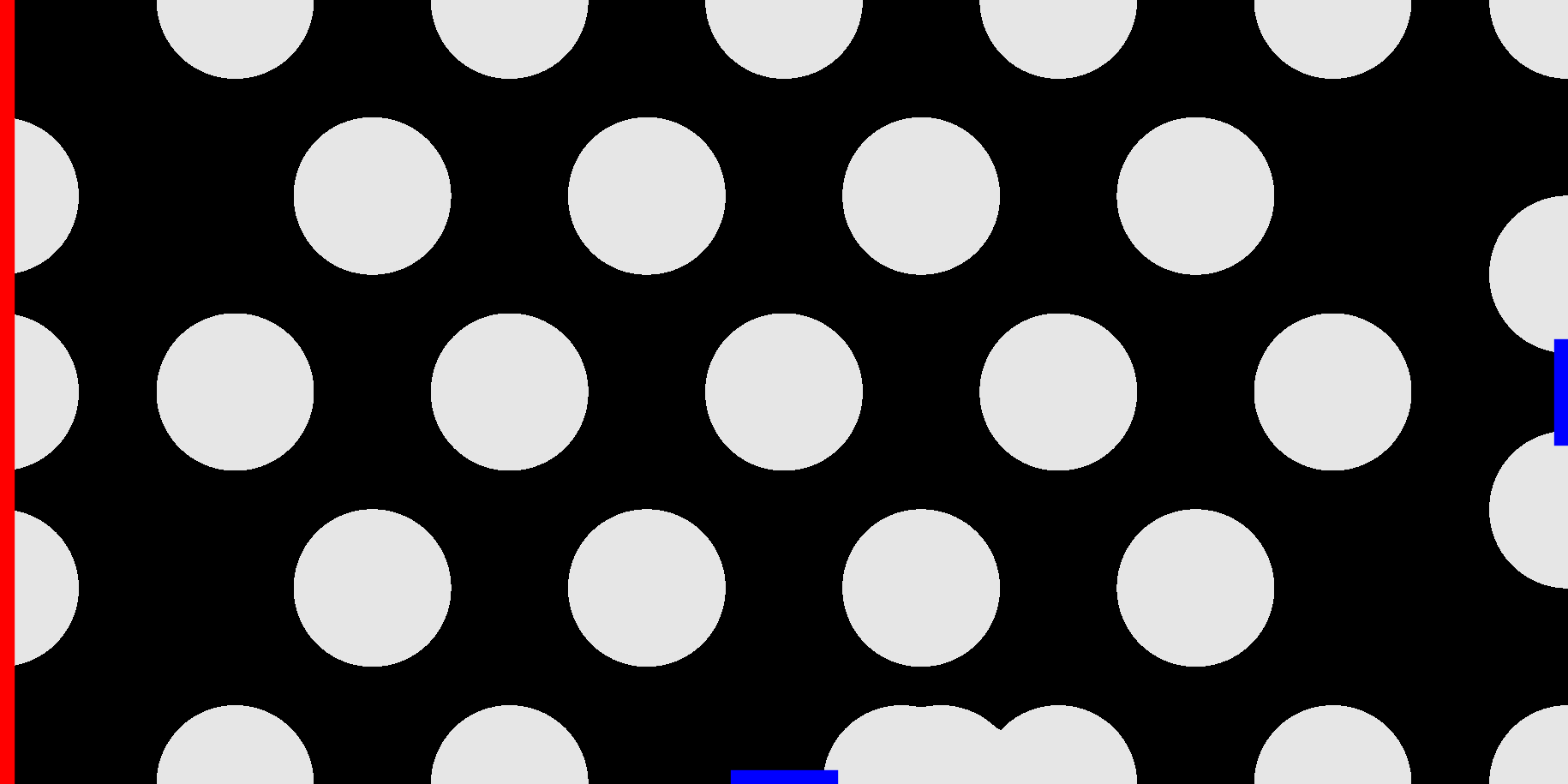}
\hfill
\includegraphics[width=0.48\linewidth]{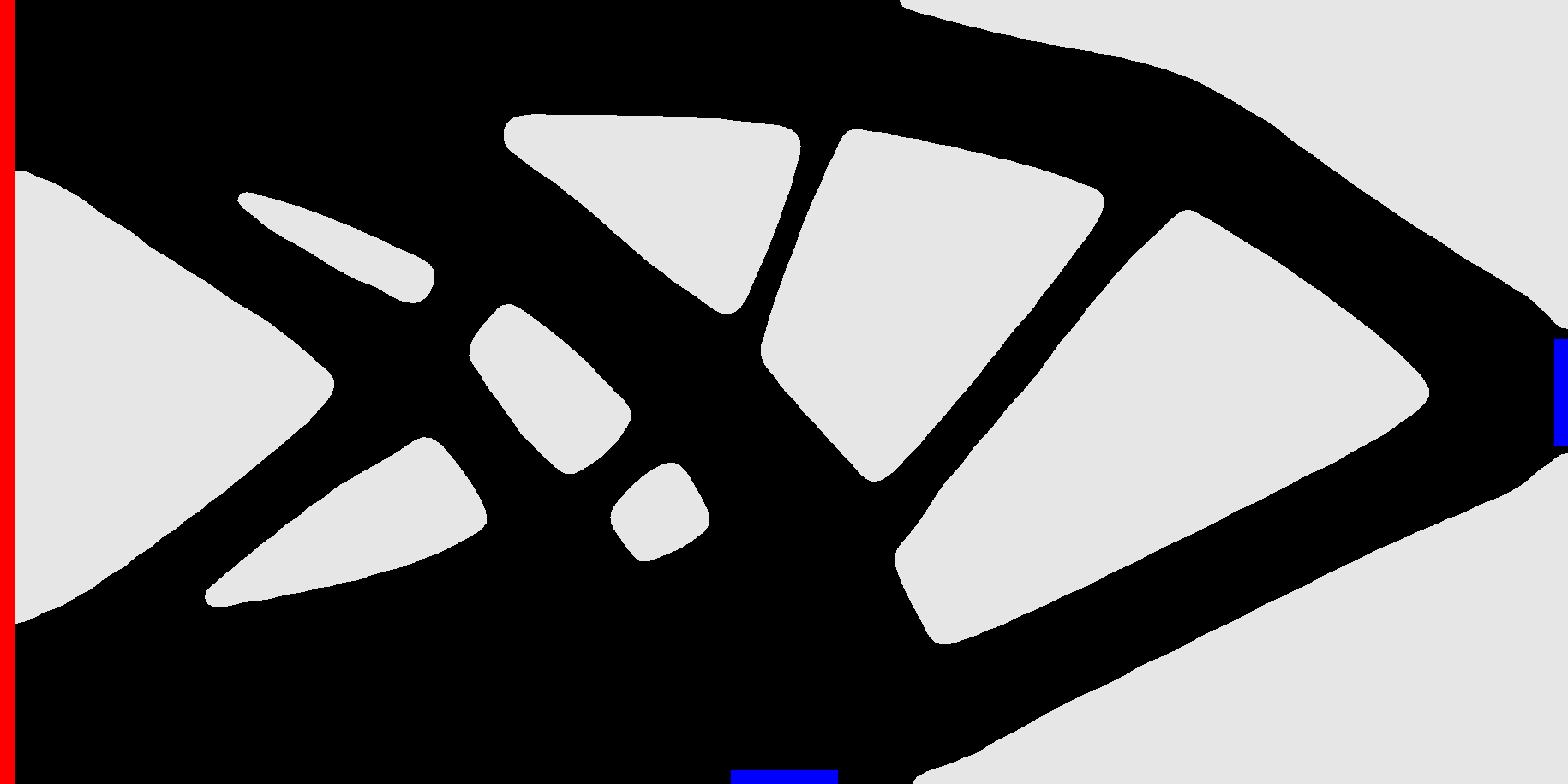}
\captionof{figure}{
	Compliance minimization, \textit{Example 4}.
	Initial and optimal domains.
	Both parallelism modalities require $79$ iterations.
}
\label{fig:compli4}
\end{figure}

\noindent\textit{Example 4.} (Cantilever problem with two loads II)
This example also applies multiple loads.
We consider the rectangular domain $\mathcal{D} = (0, 2)\times(0, 1)$.
The forces $g^{\mathrm{I}} = (0, -2)^\top$ and $g^{\mathrm{II}} = (0, 2)^\top$
are applied separately on $\Gamma_1^{\mathrm{I}} = (0.95,1.05) \times \{y=0\}$ (bottom-center side)
and $\Gamma_1^{\mathrm{II}} = \{x=2\} \times (0.45,0.55)$ (left-center side), respectively.
The boundary of zero displacement is $\Gamma_0 = \{x=0\} \times (0, 1)$ 
and the area constraint is $V = 1.1$.
The performance of data and task parallelisms are compared using the same parameter values,
on a mesh of $52\,155$ triangles:
\begin{lstlisting}
mpirun -np 2 python test.py 18 # Data (70 sec)
mpirun -np 2 python test.py 19 # Task (98 sec)
\end{lstlisting}
\begin{lstlisting}
ctrn_tol=1e-3, lgrn_tol=1e-3, dfactor=1e-1,
reinit_step=4, reinit_pars=(20, 0.01), smooth=True
\end{lstlisting}
The bilinear form $\eqref{eq:compli1_bili}$ is employed,
with the subdomain
\[
\Omega_0\coloneqq (0.94, 1.06)\times[0,0.05)
		\cup (1.95, 2]\times(0.42,0.58)
\]
to penalize the material around $\Gamma_1^{\mathrm{I}}$
and $\Gamma_1^{\mathrm{II}}$.
Recall that in task parallelism the number of processes must be
equal to the number of state/adjoint problems ($2$ in this example).
See the result in Figure~\ref{fig:compli4}.

\subsection{Heat conduction}\label{sub:heat}

Shape and topology optimization play a crucial role for efficient thermal management 
in modern applications such as battery thermal regulation \cite{Mo2021},
heat exchangers~\cite{GAO2008805},
additive-manufactured cooling channels~\cite{LamarcheGagnon2024}, and
temperature-sensitive components in aerospace systems.

Let $\Gamma_0\subset\Gamma=\partial\hold$ and  $V$ be a prescribed volume, smaller than that of $\mathcal{D}$.
Consider the minimization problem
\begin{equation} \label{eq:heat_cost}
	\min_{\Omega \subset \mathcal{D}}
	J(\Omega) \coloneqq \int_{\mathcal{D}} A_{\Omega}\left|\nabla u\right|^{2}
\end{equation}
subject to
\begin{equation}
	\int_{\mathcal{D}} \chi_{\Omega} = V,
\end{equation}
where $u$ is the solution to the following transmission problem
\begin{equation} \label{eq:heat_prob}
	\begin{aligned}
	-\mathrm{div} A_{\Omega} \nabla u &= f && \text{in } \Omega \text{ and } \mathcal{D}\setminus\overline{\Omega}\\
                                    u &= 0 && \text{on } \Gamma_{0}\\
                       \partial_{n} u &= 0 && \text{on } \Gamma\setminus\Gamma_{0}\\
		(A_{\Omega} \partial_{n} u)^+ &= (A_{\Omega} \partial_{n} u)^- && \text{on }  \partial\Omega\\
							      u^+ &= u^- && \text{on }  \partial\Omega        
	\end{aligned}	
\end{equation}
and $A_{\Omega}:\overline{\mathcal{D}} \rightarrow \mathds{R}$ is given by
\begin{equation}\label{eq:heat_AOmega}
A_{\Omega} = \chi_{\Omega} + {10}^{-3}\chi_{\overline{\mathcal{D}}\setminus\Omega}.
\end{equation}
According to \eqref{eq:problem}, the constraint function for this problem is written as
\begin{equation}
	C(\Omega) = \frac{1}{V} \int_{\mathcal{D}} \chi_{\Omega}.
\end{equation}
The weak formulation of \eqref{eq:heat_prob} reads:
Find $u \in  H^1_{\Gamma_0}(\mathcal{D})$ such that
\begin{equation} \label{eq:heat_weak}
	\int_{\mathcal{D}} A_{\Omega} \nabla u \cdot \nabla v =
	\int_{\mathcal{D}} f v
	\quad
	\forall \, v \in  H^1_{\Gamma_0}(\mathcal{D}).
\end{equation}
In this case, the adjoint problem is the same as \eqref{eq:heat_weak}
(except by a minus sign on the right-hand side), and thus $p=-u$.
The shape derivative components are given by
\begin{align*}
	S^J_0 &= 2u\nabla f, \\
	S^J_1 &= (2uf-A_{\Omega}{|\nabla u|}^{2})I +
				2A_{\Omega}\nabla u\otimes\nabla u
\end{align*}
for the cost functional, and
\begin{equation*}
	S^C_0 = \boldsymbol{0},\quad S^C_1 = \frac{1}{V} \chi_\Omega I
\end{equation*}
for the constraint function.
In all examples, we use the square domain $\mathcal{D} = (0, 1)\times(0, 1)$,
and red color to highlight $\Gamma_0$, where the temperature is fixed to zero.

\begin{figure}[H]
	\centering
\includegraphics[width=0.48\linewidth]{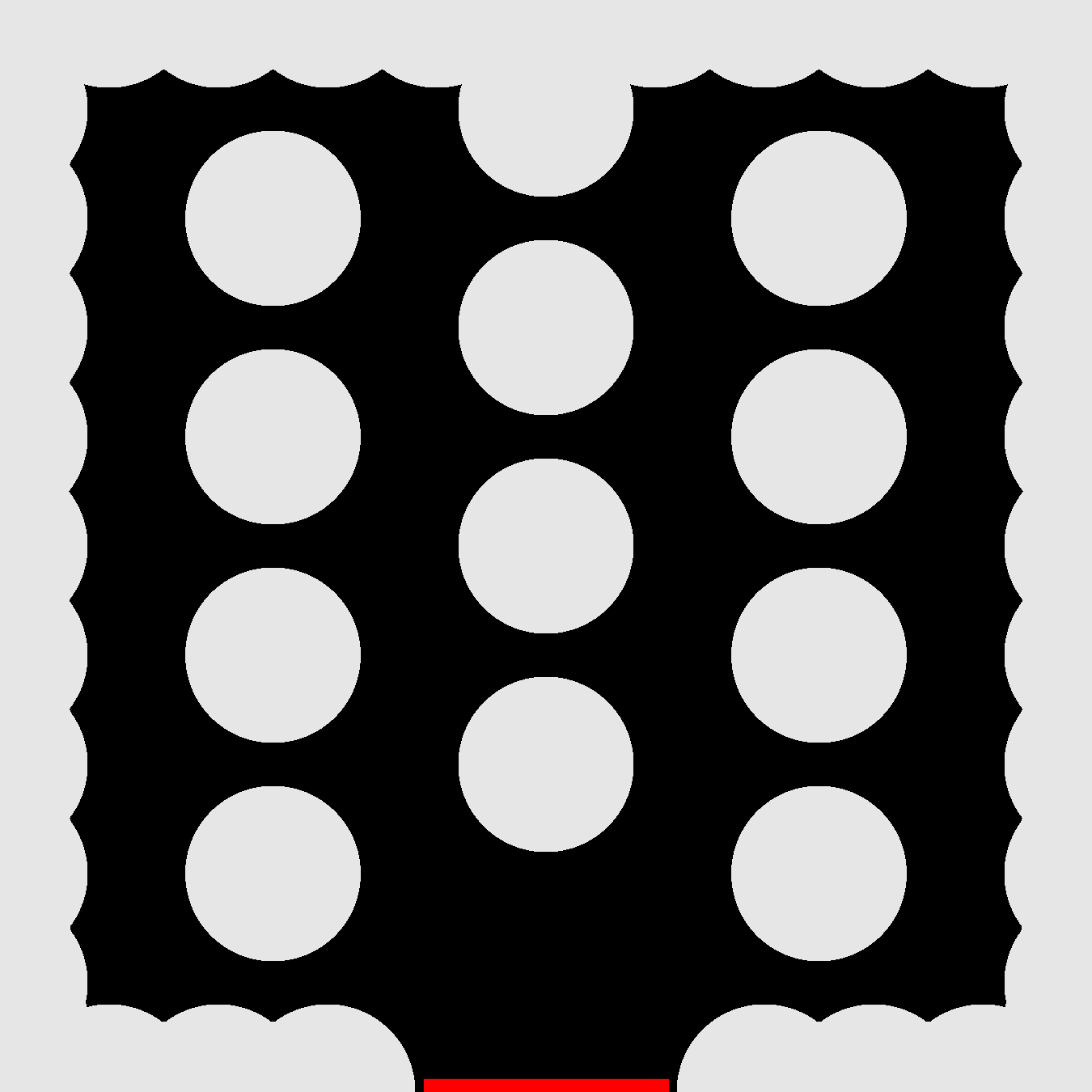}
\hfill
\includegraphics[width=0.48\linewidth]{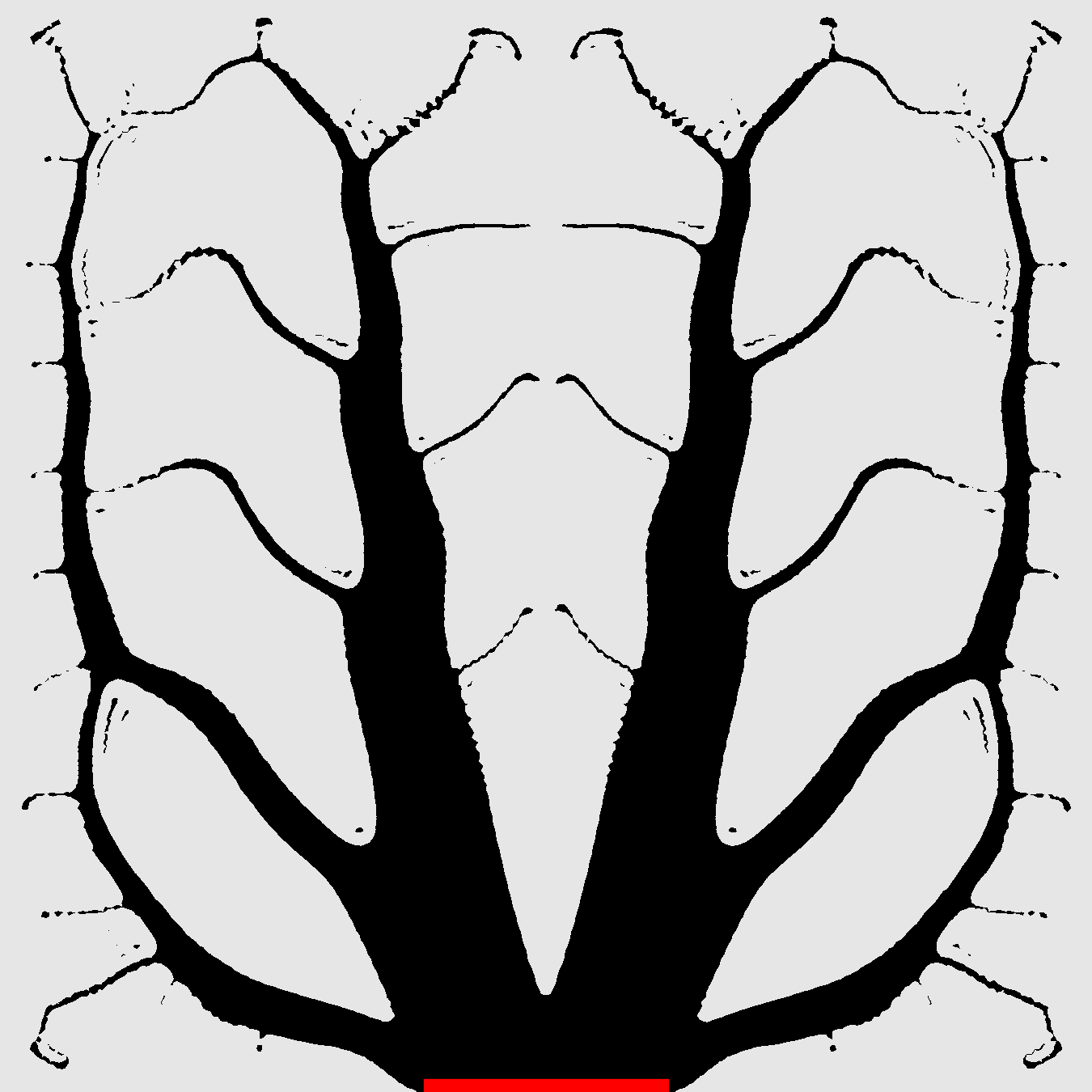}
\captionof{figure}{
	Heat conduction, \textit{Example 1}.
	Initial guess and optimized design at iteration $i=204$,
	computed on a mesh with {92\,582} triangles.
	The black-colored region $\Om$ has the highest conductivity.
}
\label{fig:heat1}
\end{figure}

\noindent\textit{Example 1.}
We apply a uniform heat $f \equiv 1$,
with volume constraint $V = 0.25$ and $\Gamma_0 = (0.4, 0.6)\times\{y=0\}$.
The test is carried out using data parallelism with two processes:
\begin{lstlisting}
mpirun -np 2 python test.py 09 # Data (58 sec)
\end{lstlisting}
See the results in Figure \ref{fig:heat1}.
In this example, we have employed the bilinear form
\begin{equation*}
	B(\theta, \xi) = \int_{\mathcal{D}} D\theta : D\xi
	+ {10}^4(\theta \cdot \xi) \chi_{\Omega_0}, 
\end{equation*}
with $\theta, \xi \in \mathbb{H} = H^1_0(\mathcal{D})^d$
and $\Omega_0 = (0.3, 0.7)\times[0, 0.05)$, 
in order to penalize the material around $\Gamma_0$.

\begin{figure}[H]
	\centering
\includegraphics[width=0.48\linewidth]{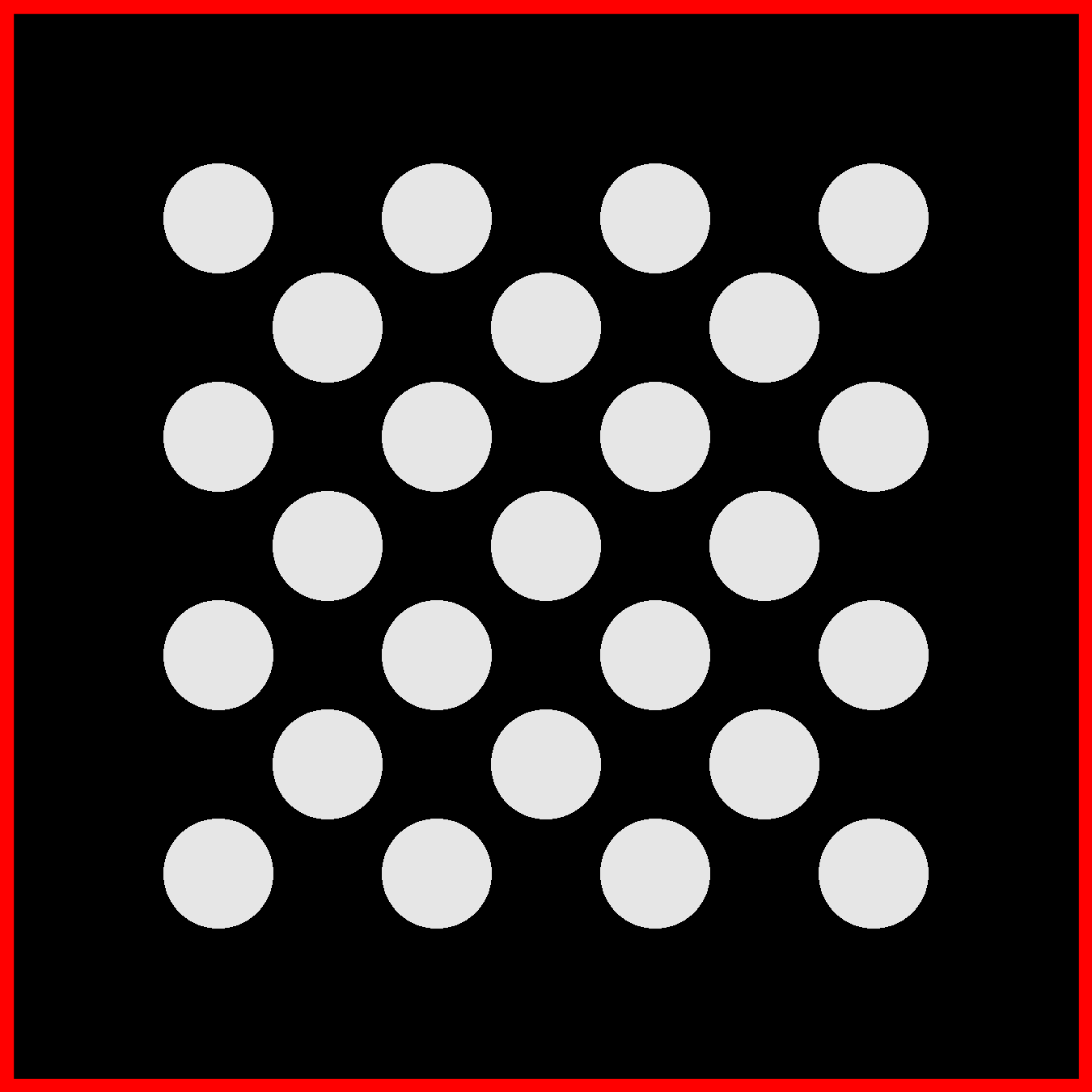}
\hfill
\includegraphics[width=0.48\linewidth]{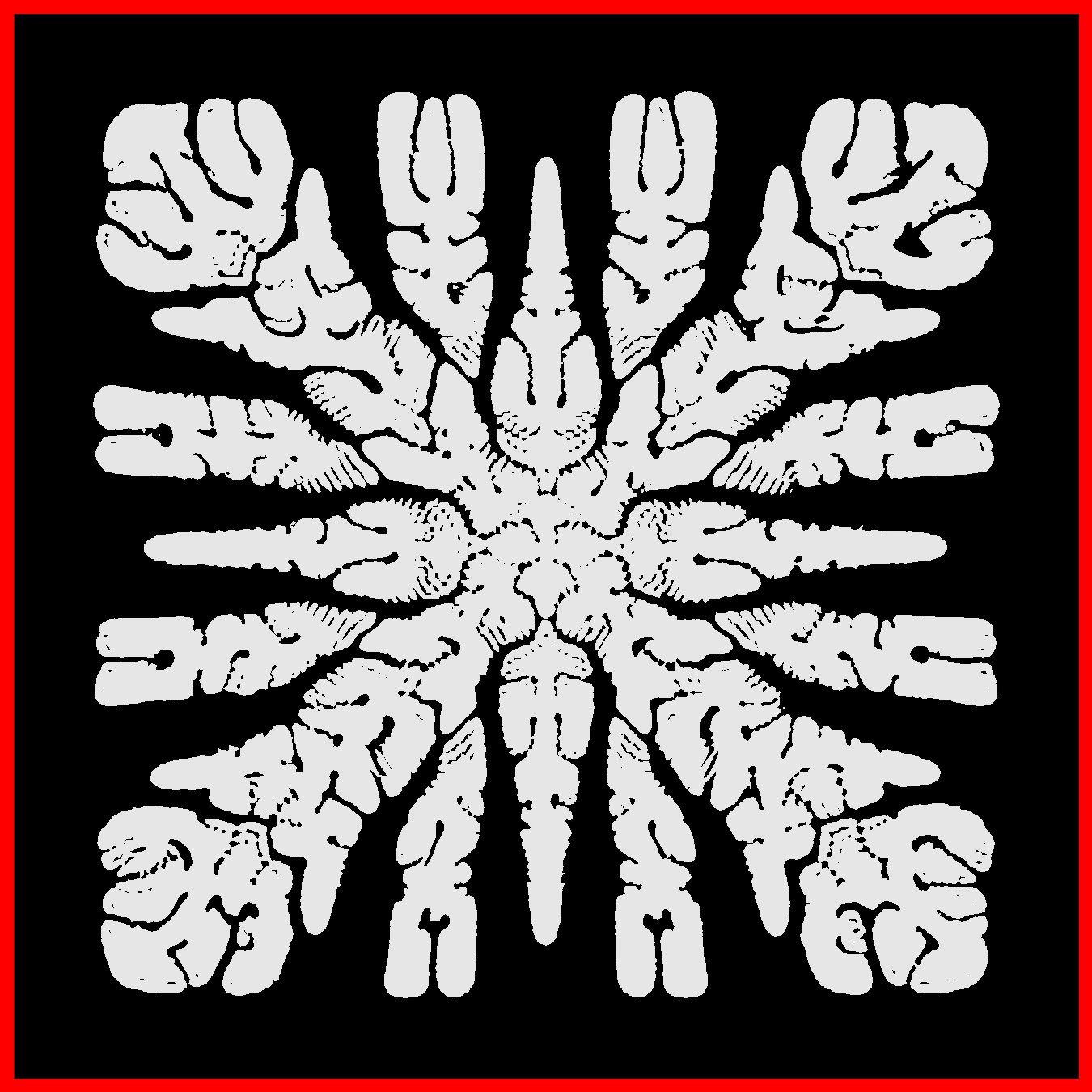}
\captionof{figure}{
	Heat conduction, \textit{Example 2.}
	Initial guess and optimized design at iteration $i=172$,
	computed on a mesh with {144\,712} triangles.
	The black-colored region $\Om$ has the highest conductivity.
}
\label{fig:heat2}
\end{figure}

\noindent\textit{Example 2.}
We apply a uniform heat $f\equiv 1$, with volume constraint $V = 0.6$
and $\Gamma_0 = \Gamma$ (i.e., the whole boundary is kept at zero temperature).
We run this example using data parallelism with four processes:
\begin{lstlisting}
mpirun -np 4 python test.py 10 # Data (52 sec)
\end{lstlisting}
See the results in Figure \ref{fig:heat2}.
The bilinear form $B$ is the same as in \textit{Example 1},
but with $\Omega_0$ defined as
\[
\Omega_0 \coloneqq \{x<0.1\} \cup \{x>0.9\} \cup \{y<0.1\} \cup \{y>0.9\}.
\]
Thus, the material is penalized close to the boundary.
Notice the typical fractal structure of the optimized geometries
in Figures~\ref{fig:heat1} and~\ref{fig:heat2};
compare the results with those in~\cite{GAO2008805}.

\begin{figure}[H]\centering
\includegraphics[width=0.48\linewidth]{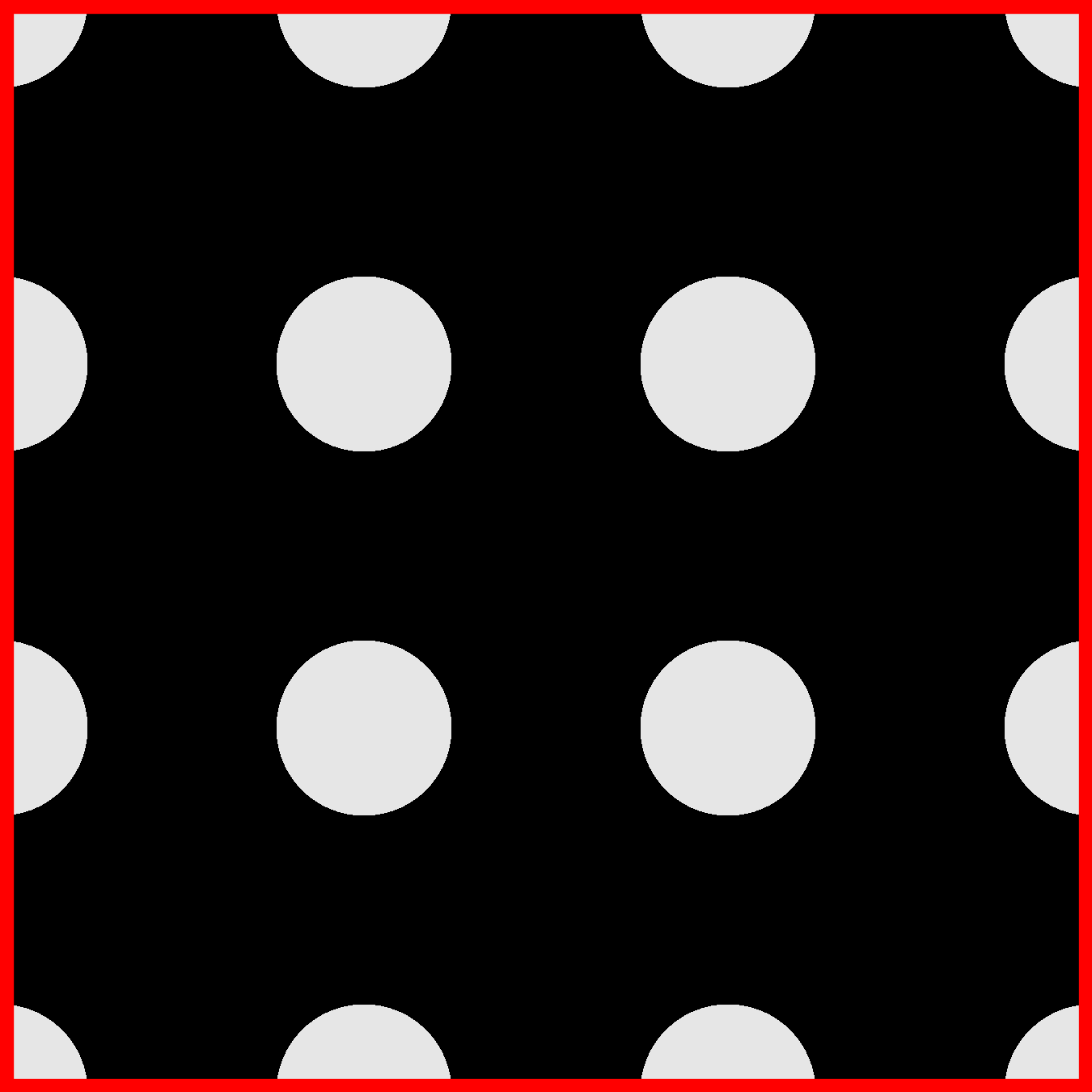}
\hfill
\includegraphics[width=0.48\linewidth]{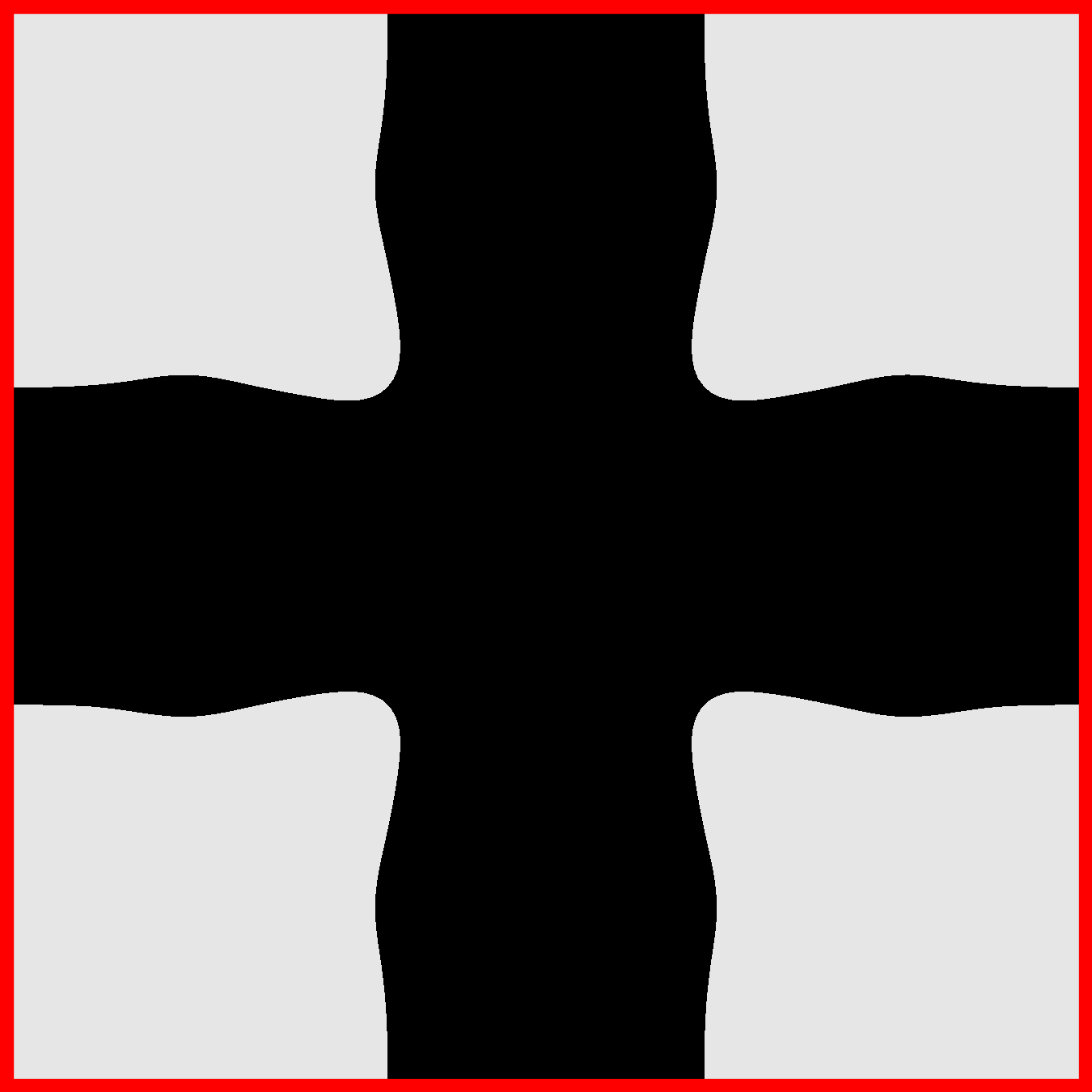}
\captionof{figure}{
	Heat conduction, \textit{Example 3}.
	Initial guess and optimized design at iteration $i=181$,
	computed on a mesh with {23\,652} triangles.
}
\label{fig:heat3}
\end{figure}

\noindent\textit{Example 3.} 
We apply a localized heat source given by the function
$f(x,y) \coloneqq \omega(\vert(x,y)-(0.5, 0.5)\vert_2)$,
where $\omega$ is the radial function defined as
\begin{equation} \label{eq:rad}
	\omega(r) \coloneqq
	\begin{cases}
	25(1 + \cos(10 \pi r )) & \text{if } r < 0.1, \\
	0  & \text{otherwise}.
	\end{cases}
\end{equation}
Thus, $f$ is radially symmetric around the point 
$(0.5,0.5)$ and infinitely differentiable on the entire plane.
The support of $f$ is the disk centered at $(0.5,0.5)$ with radius $0.1$.
The volume constraint is $V = 0.5$ and $\Gamma_0=\Gamma$.
We run this example using data parallelism with six processes:
\begin{lstlisting}
mpirun -np 6 python test.py 11 # Data (11 sec)
\end{lstlisting}
See the results in Figure \ref{fig:heat3}.
The term \eqref{eq:bd_plt} is added to the bilinear form
$B(\theta, \xi) = \int_{\mathcal{D}} D\theta : D\xi$ and
the velocity problem \eqref{eq:velocity} is solved in $\mathbb{H}=H^1(\mathcal{D})^d$,
allowing tangential displacements along the boundary $\partial\mathcal{D}$.

\begin{figure}[H]\centering
\includegraphics[width=0.48\linewidth]{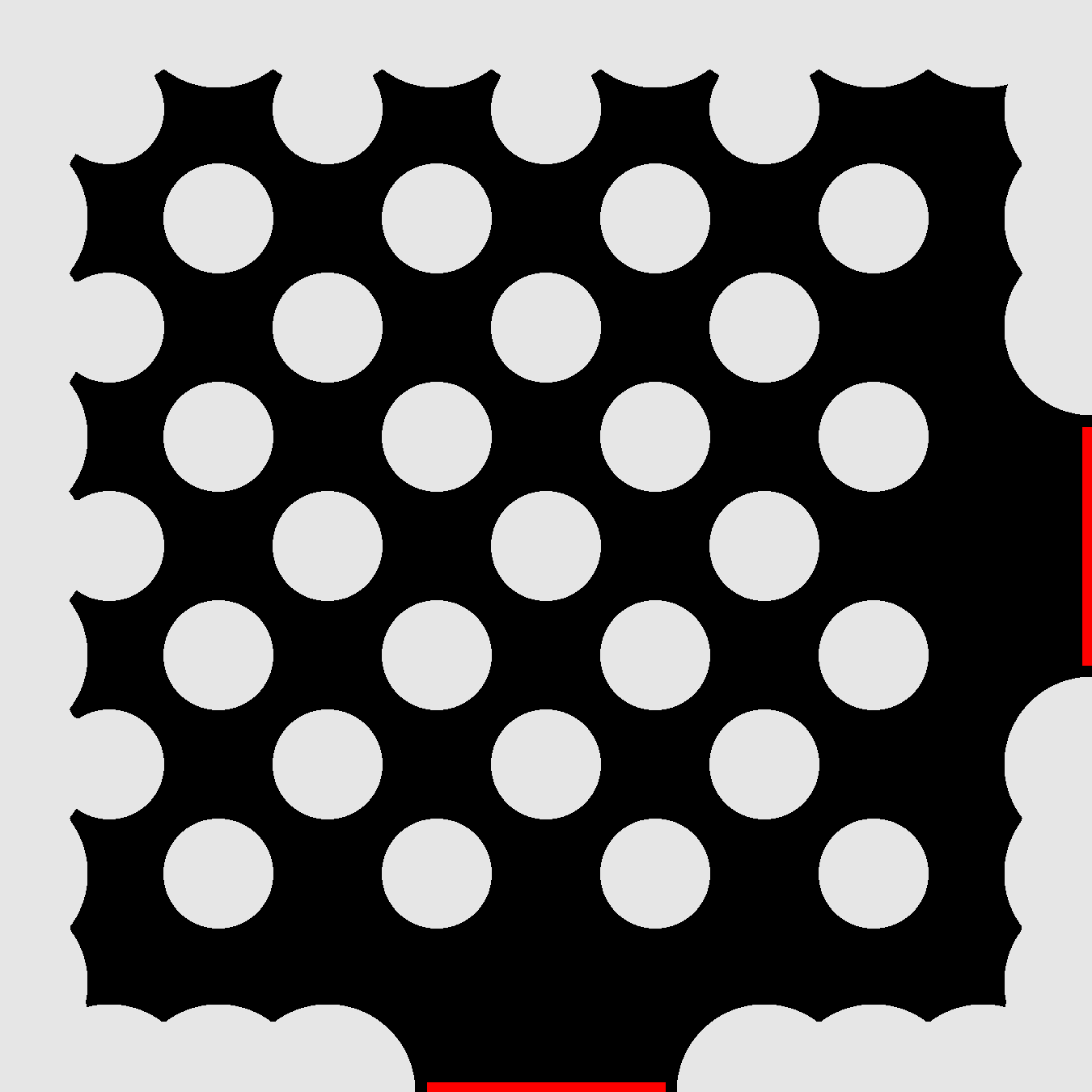}
\hfill
\includegraphics[width=0.48\linewidth]{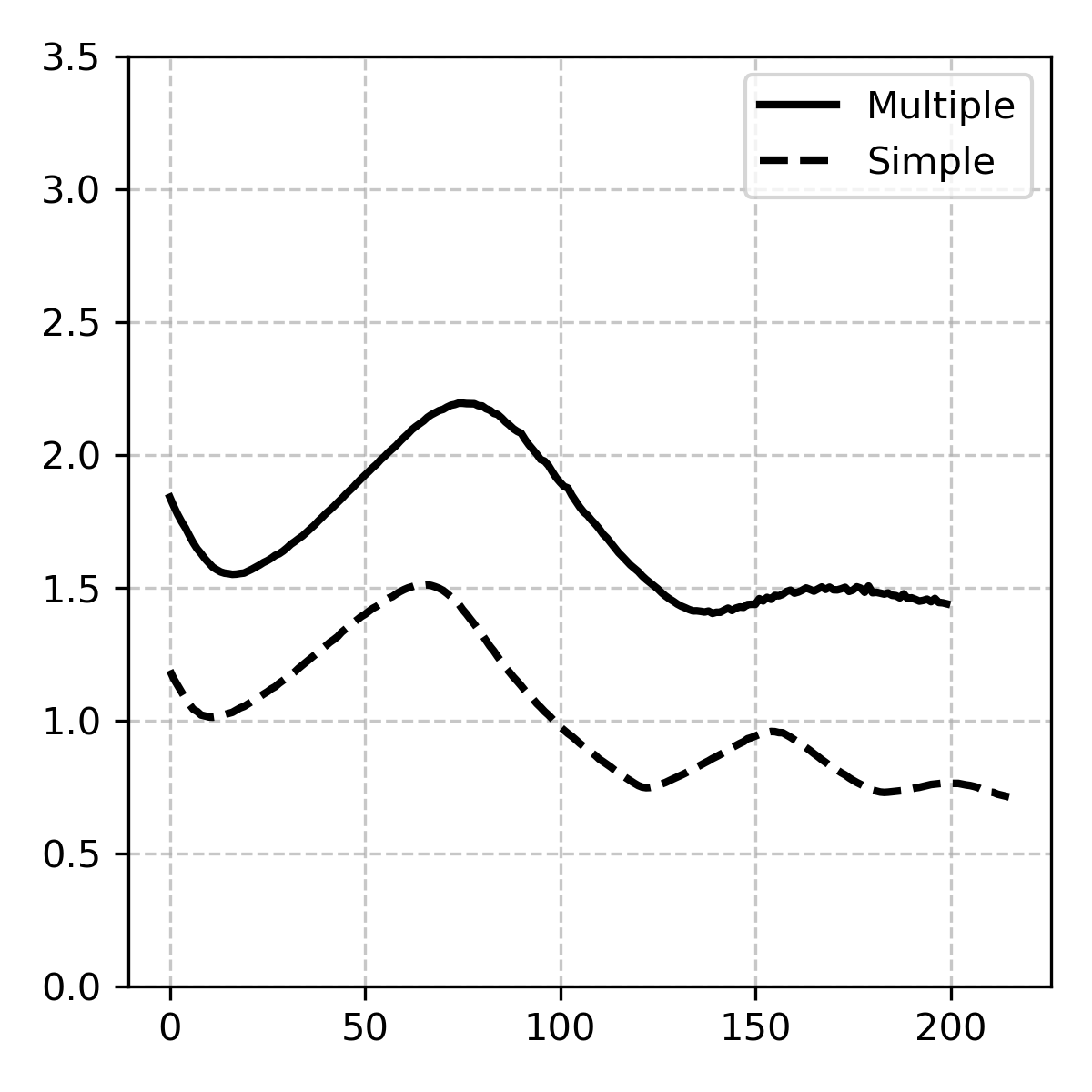}
\includegraphics[width=0.48\linewidth]{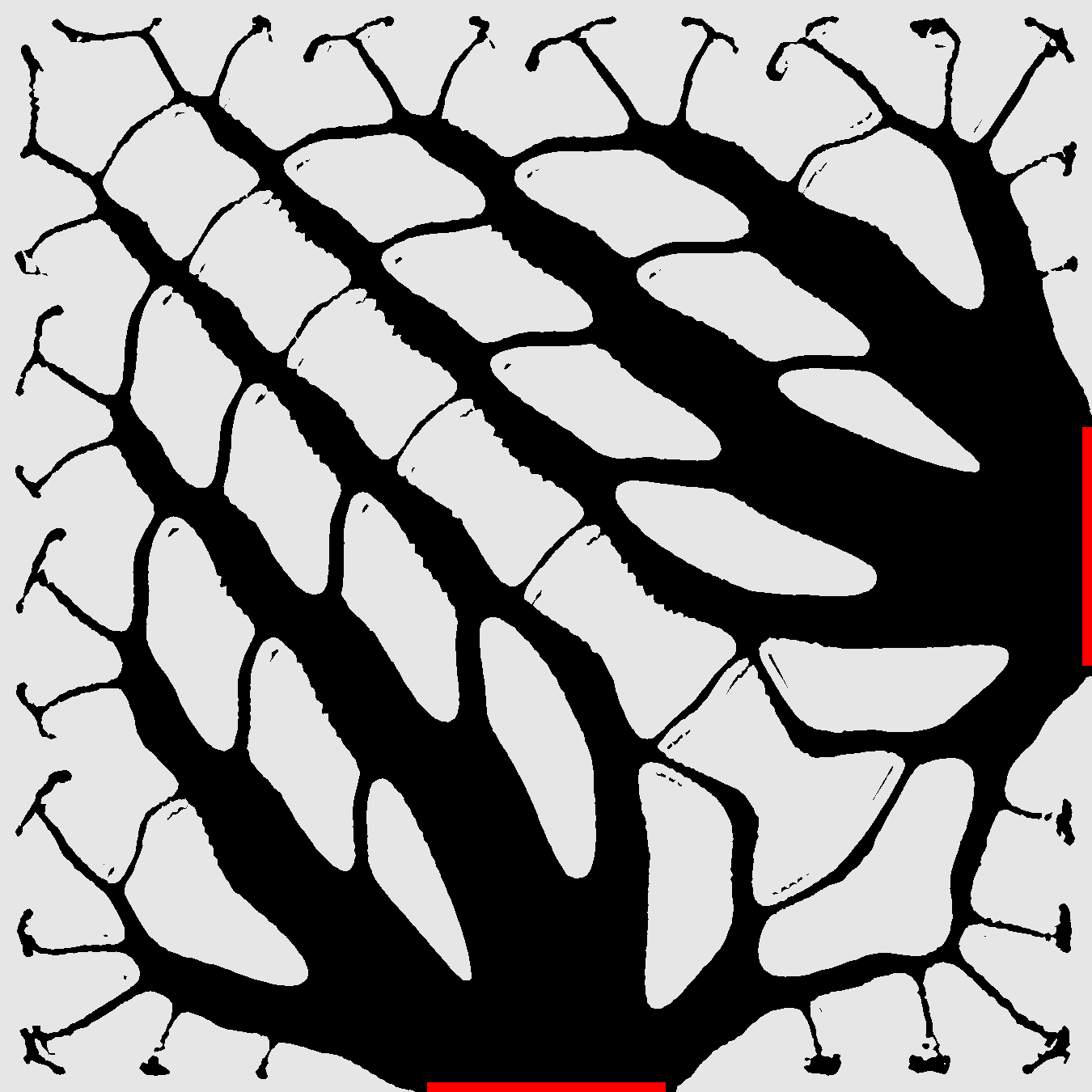}
\hfill
\includegraphics[width=0.48\linewidth]{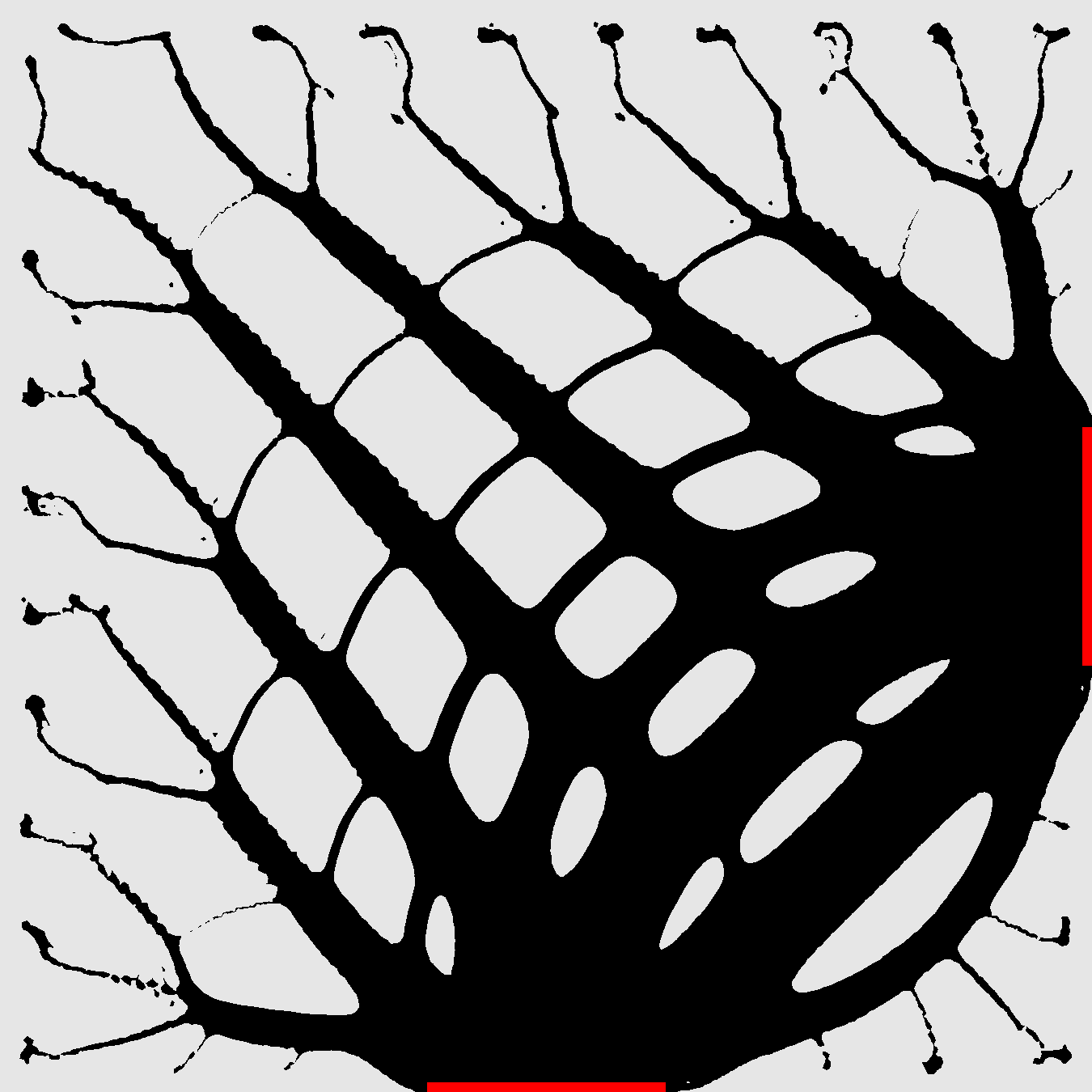}
\captionof{figure}{
	Heat conduction, \textit{Example 4}.
	Initial guess (top-left) and cost values (top-right),
	along with the optimized designs of the single case (bottom-left)
	and multiple case (bottom-right), at iteration $i=215$ and $i=199$, respectively.
	The final cost value of the multiple sink test
	approximately duplicates the one of the single sink test.
}
\label{fig:heat4}
\end{figure}

\noindent\textit{Example 4.}
In this example, we perform two tests.
In the first test, we solve the problem with one heat sink
$\Gamma_0 = \Gamma_0^\mathrm{I} \cup \Gamma_0^\mathrm{II}$,
where
\[
\Gamma_0^\mathrm{I} = (0.4,0.6)\times\{y=0\}\;,\;\Gamma_0^\mathrm{II} = \{x=1\}\times(0.4,0.6).
\]
In the second test, we solve the problem with multiple heat sinks
$\Gamma_0^\mathrm{I}$ and $\Gamma_0^\mathrm{II}$. In this case,
the cost functional is given by the sum of the corresponding thermal compliances:
\[
	J(\Omega) \coloneqq \frac{1}{2}\int_{\mathcal{D}} A_{\Omega}\left|\nabla u^\mathrm{I}\right|^{2}+
	\frac{1}{2}\int_{\mathcal{D}} A_{\Omega}\left|\nabla u^\mathrm{II}\right|^{2},
\]
where $u^\mathrm{I}$ solves \eqref{eq:heat_prob} with $\Gamma_0 = \Gamma_0^\mathrm{I}$
and $u^\mathrm{II}$ solves \eqref{eq:heat_prob} with $\Gamma_0 = \Gamma_0^\mathrm{II}$.
In both cases, $f\equiv 1$ and $V = 0.4$.
The first test is carried out using data parallelism with two processes:
\begin{lstlisting}
mpirun -np 2 python test.py 12 # Data (68 sec)
\end{lstlisting}
Since in the second test there are two PDEs to be solved,
we employ task parallelism with two processes:
\begin{lstlisting}
mpirun -np 2 python test.py 13 # Task (106 sec)
\end{lstlisting}
The results are shown in Figure \ref{fig:heat4}.

\begin{figure}[H]\centering
\includegraphics[width=0.48\linewidth]{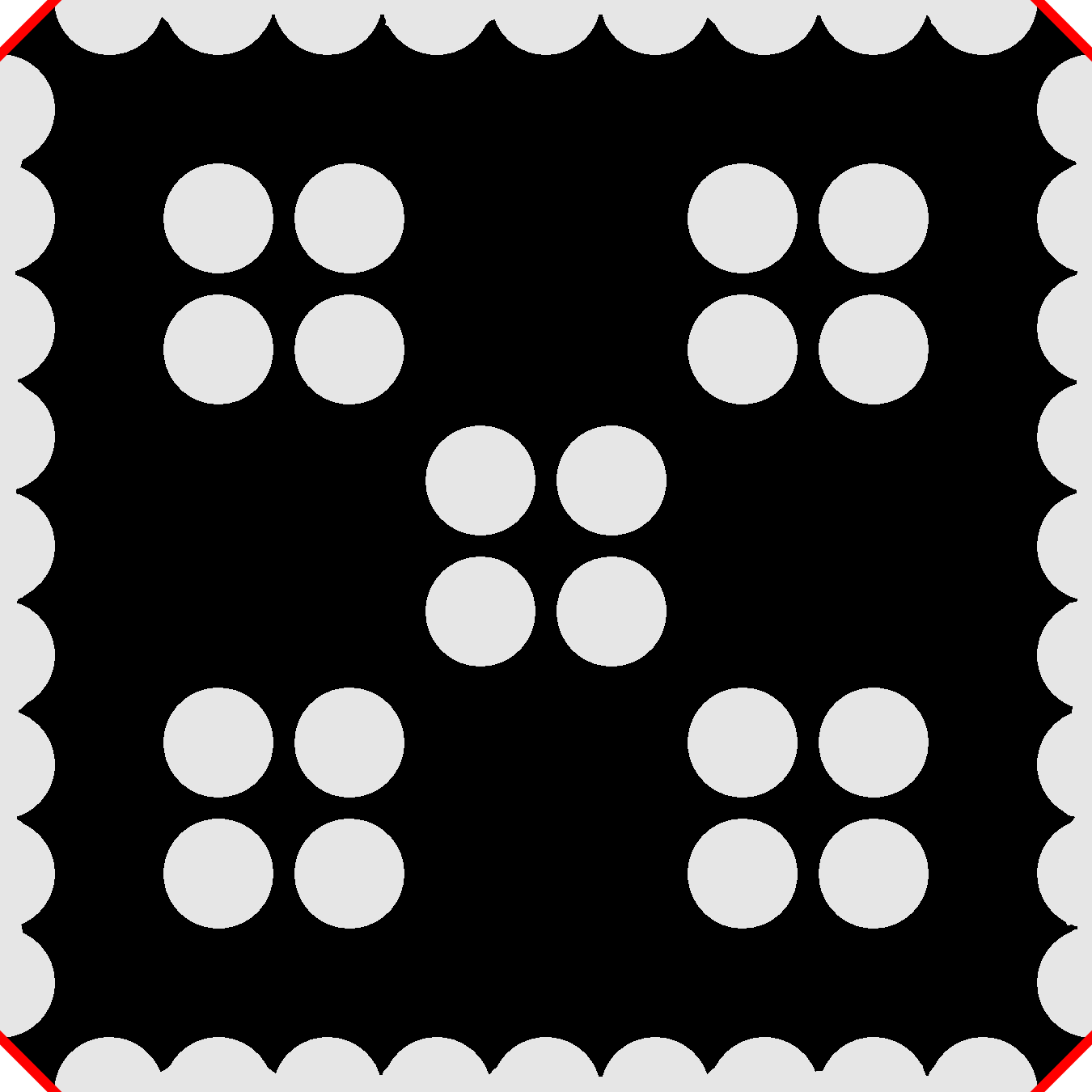}
\hfill
\includegraphics[width=0.48\linewidth]{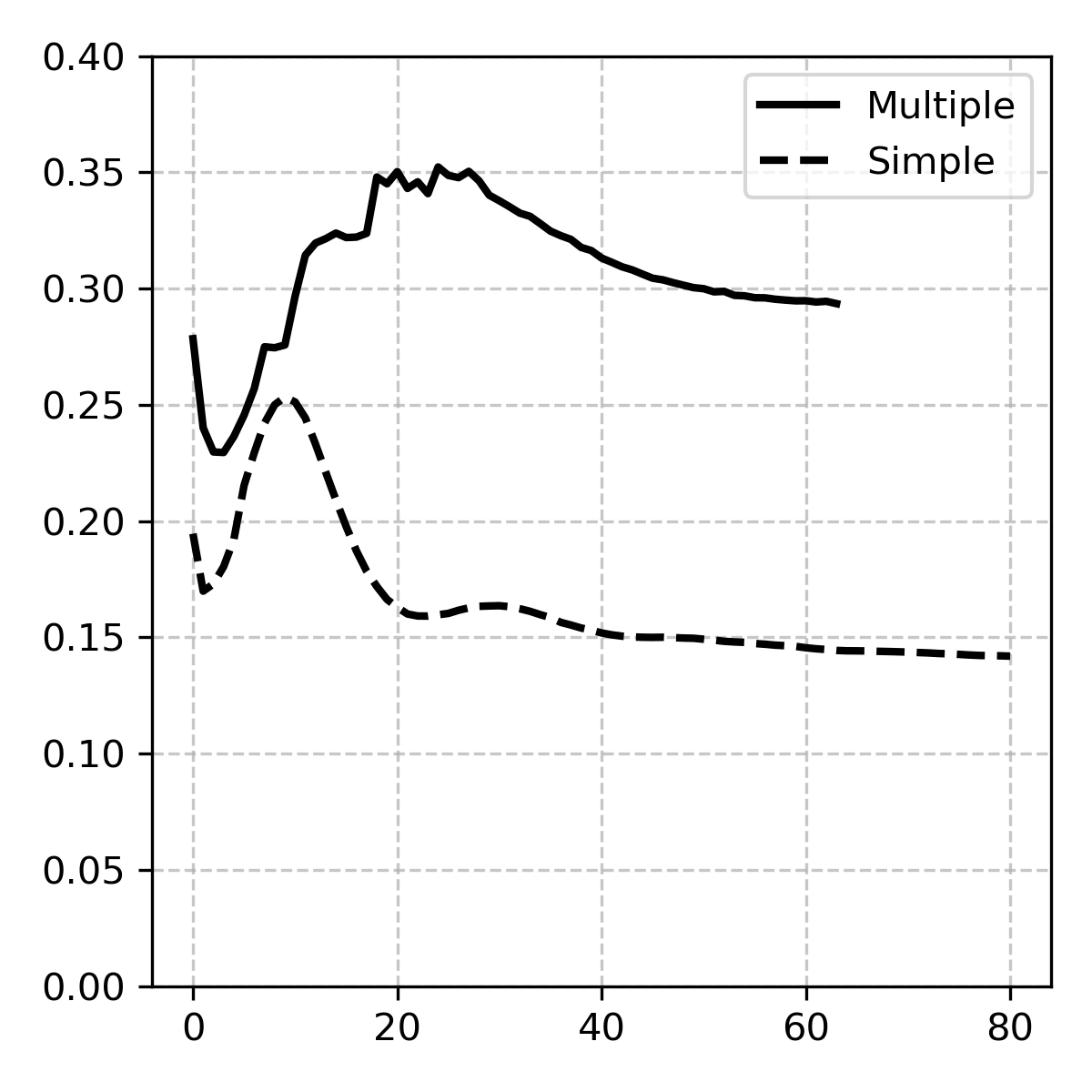}
\includegraphics[width=0.48\linewidth]{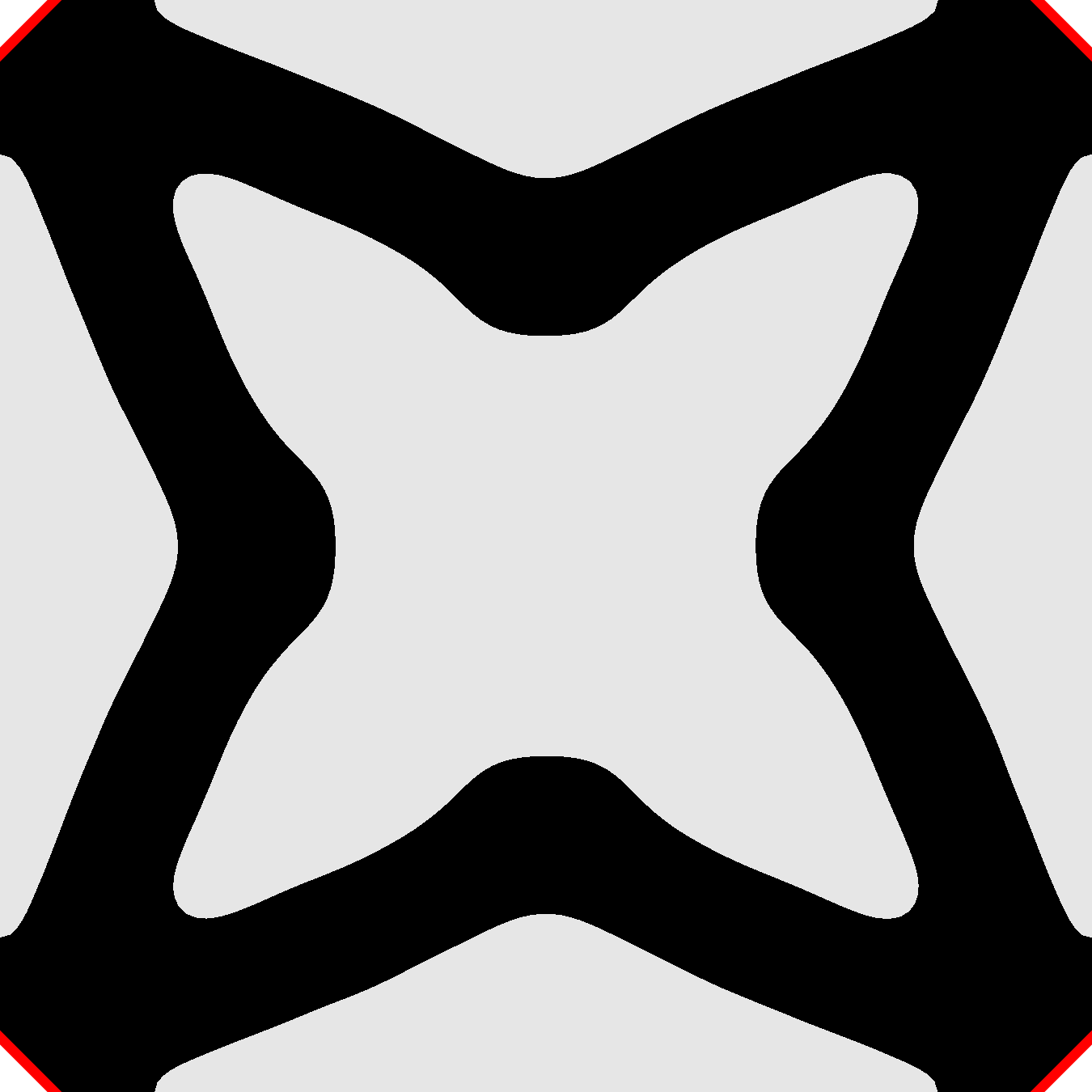}
\hfill
\includegraphics[width=0.48\linewidth]{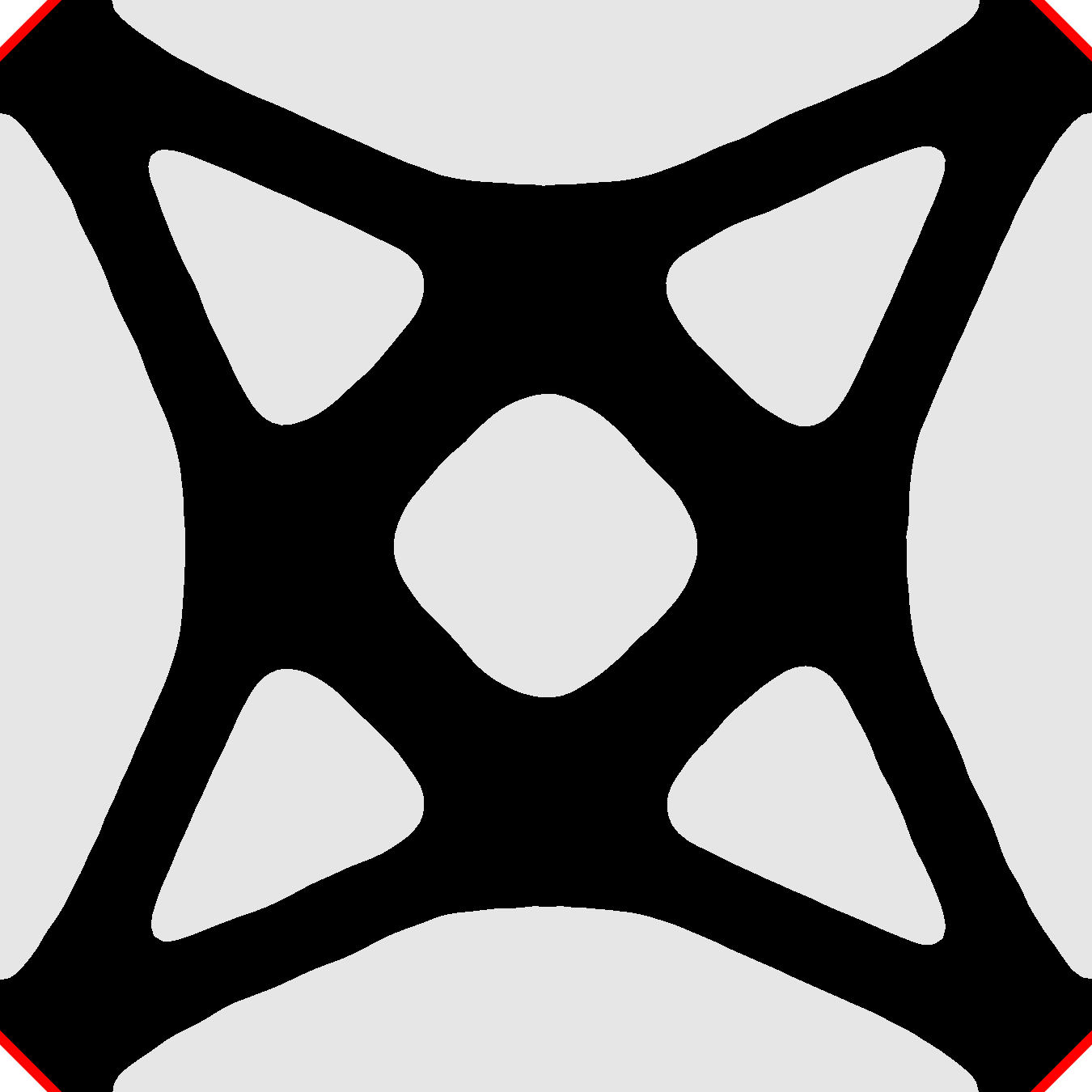}
\captionof{figure}{
	Heat conduction, \textit{Example 5}.
	Initial guess (top-left) and cost values (top-right),
	along with the optimized designs of the single case (bottom-left)
	and multiple case (bottom-right), at iteration $i=80$ and $i=63$, respectively.
	The final cost value of the multiple source test
	approximately duplicates the one of the single source test.
}
\label{fig:heat5}
\end{figure}

\noindent\textit{Example 5.}
As in the previous example, here we consider the single and multiple cases,
but for the heat source.
In the single case test, we apply one heat source 
$f = f^\mathrm{I} + f^\mathrm{II} + f^\mathrm{III} + f^\mathrm{IV}$,
where
\[
	f^\mathrm{I}(x,y) \coloneqq \omega(\vert(x,y)-(0.5, 0.25)\vert_2),
\]
\[
	f^\mathrm{II}(x,y) \coloneqq \omega(\vert(x,y)-(0.75, 0.5)\vert_2),
\]
\[
	f^\mathrm{III}(x,y) \coloneqq \omega(\vert(x,y)-(0.5, 0.75)\vert_2),
\]
\[
	f^\mathrm{IV}(x,y) \coloneqq \omega(\vert(x,y)-(0.25, 0.5)\vert_2).
\]
The function $\omega$ was defined in \eqref{eq:rad}.
In the multiple case, we consider four heat sources 
$f=f^\mathrm{I}$, $f=f^\mathrm{II}$, $f=f^\mathrm{III}$, and $f=f^\mathrm{IV}$,
along with the cost functional
\[
	J(\Omega) \coloneqq \sum_{i=\mathrm{I}}^{\mathrm{IV}}\frac{1}{4}
	\int_{\mathcal{D}} A_{\Omega}\left|\nabla u^i\right|^{2},
\]
where $u^i$ solves \eqref{eq:heat_prob} with $f = f^i$.
In both cases, the volume constraint is $V = 0.45$, and 
$\Gamma_0$ consists of four small subsets
of the boundary in the corners of the domain, each one with length $0.05\sqrt{2}$.
The first test is carried out using data parallelism with two processes:
\begin{lstlisting}
mpirun -np 2 python test.py 21 # Data (15 sec)
\end{lstlisting}
Since in the second test there are four PDEs to be solved,
we employ task parallelism with four processes:
\begin{lstlisting}
mpirun -np 4 python test.py 22 # Task (20 sec)
\end{lstlisting}
The numerical results are shown in Figure~\ref{fig:heat5}.
They agree with the results reported in \cite{MR2278190}, which we are able to reproduce within our framework.

\subsection{A nonlinear PDE arising in population dynamics}\label{sub:pop}

We present an example in which the PDE constraint is nonlinear.
Consider the problem of spatially distributing resources in order
to maximize the population size of a species that consumes those resources.
We assume that the population growth is governed by the logistic equation with a diffusive term.
Since the temporal derivative is absent, the objective is to maximize the population size in the long-time regime,
when the species distribution has reached a stationary state.
In this setting, the PDE constraint takes the form
of a nonlinear diffusion equation \cite{MR1908418}.

Let $V$ be a prescribed volume, smaller than the volume of~$\mathcal{D}$.
The problem reads
\begin{equation} \label{eq:log_cost}
	\min_{\Omega \subset \mathcal{D}}
	J(\Omega) \coloneqq -\int_{\mathcal{D}} u
\end{equation}
subject to
\begin{equation}
	\int_{\mathcal{D}} \chi_{\Omega} = V,
\end{equation}
where $u$ is the solution to
\begin{equation}\label{eq:log_prob}
	\begin{aligned}
	-\Delta u &= r u \left(1-\frac{u}{K_\Omega}\right) && \text{in } \Omega \text{ and } \mathcal{D}\setminus\overline{\Omega}\\
    \partial_{n} u &= 0 && \text{on } \Gamma\\     
	(\partial_n u)^+ &= (\partial_n u)^- && \text{on } \partial\Omega \\ 
				u^+ &= u^- && \text{on } \partial\Omega
	\end{aligned}
\end{equation}
and $K_{\Omega}:\overline{\mathcal{D}} \rightarrow \mathds{R}$ is given by
\begin{equation}\label{eq:log_KOmega}
	K_{\Omega} = \chi_{\Omega} + {10}^{-2}\chi_{\overline{\mathcal{D}}\setminus\Omega}.
\end{equation}
According to \eqref{eq:problem}, the constraint function for this problem is written as
\begin{equation} \label{eq:log_const}
	C(\Omega) = \frac{1}{V} \int_{\mathcal{D}} \chi_{\Omega}.
\end{equation}
Here, $u$ represents the equilibrium population density,
$r$ is a positive constant that represents the growth rate of the population in the logistic model,
and the positive function $K_{\Omega}$ represents the spatially varying carrying capacity, i.e.,
the maximum equilibrium population density that the available resources can sustain at each point.
Equation \eqref{eq:log_KOmega} indicates the presence of two types of resources,
with one permitting significantly higher growth than the other.
Moreover, we constrain the distribution of the main resource
to be supported in a region of volume $V$.

The weak formulation of \eqref{eq:log_prob} is the following:
Find $u\in H^1(\mathcal{D})$ such that
\begin{equation}\label{eq:log_weak}
	\int_{\mathcal{D}}\nabla u \cdot \nabla v =
    \int_{\mathcal{D}}ru\left(1-\frac{u}{K_{\Omega}}\right)v
    \quad
    \forall \, v \in  H^1(\mathcal{D}).
\end{equation}
From the Lagrangian functional, which is constructed as in Section \ref{sec:model_prob},
we obtain the adjoint equation:
Find $p\in H^1(\mathcal{D})$ such that
\begin{equation*}
	\int_{\mathcal{D}}\nabla p \cdot \nabla q +
    \int_{\mathcal{D}}r\left(\frac{2u}{K_{\Omega}}-1\right)pq = 
    \int_{\mathcal{D}}q
    \quad
    \forall \, q \in  H^1(\mathcal{D}).
\end{equation*}
Observe that although the state problem is nonlinear, the adjoint problem is linear.
The derivative components of the cost functional \eqref{eq:log_cost} and the constraint \eqref{eq:log_const} are
\begin{align*}
	S^J_0 &= \boldsymbol{0},\\
	S^J_1 &= \left(\nabla u\cdot\nabla p-u-ru\left(1-\frac{u}{K_{\Omega}}\right)p\right)I \\
		  &-(\nabla u\otimes\nabla p+\nabla p\otimes\nabla u),	
\end{align*}
and
\begin{equation*}
	S_0^C = \boldsymbol{0}, \quad S_1^C = \frac{1}{V} \chi_{\Omega}I.
\end{equation*}

\begin{figure}[H]
		\centering
	\includegraphics[width=0.48\linewidth]{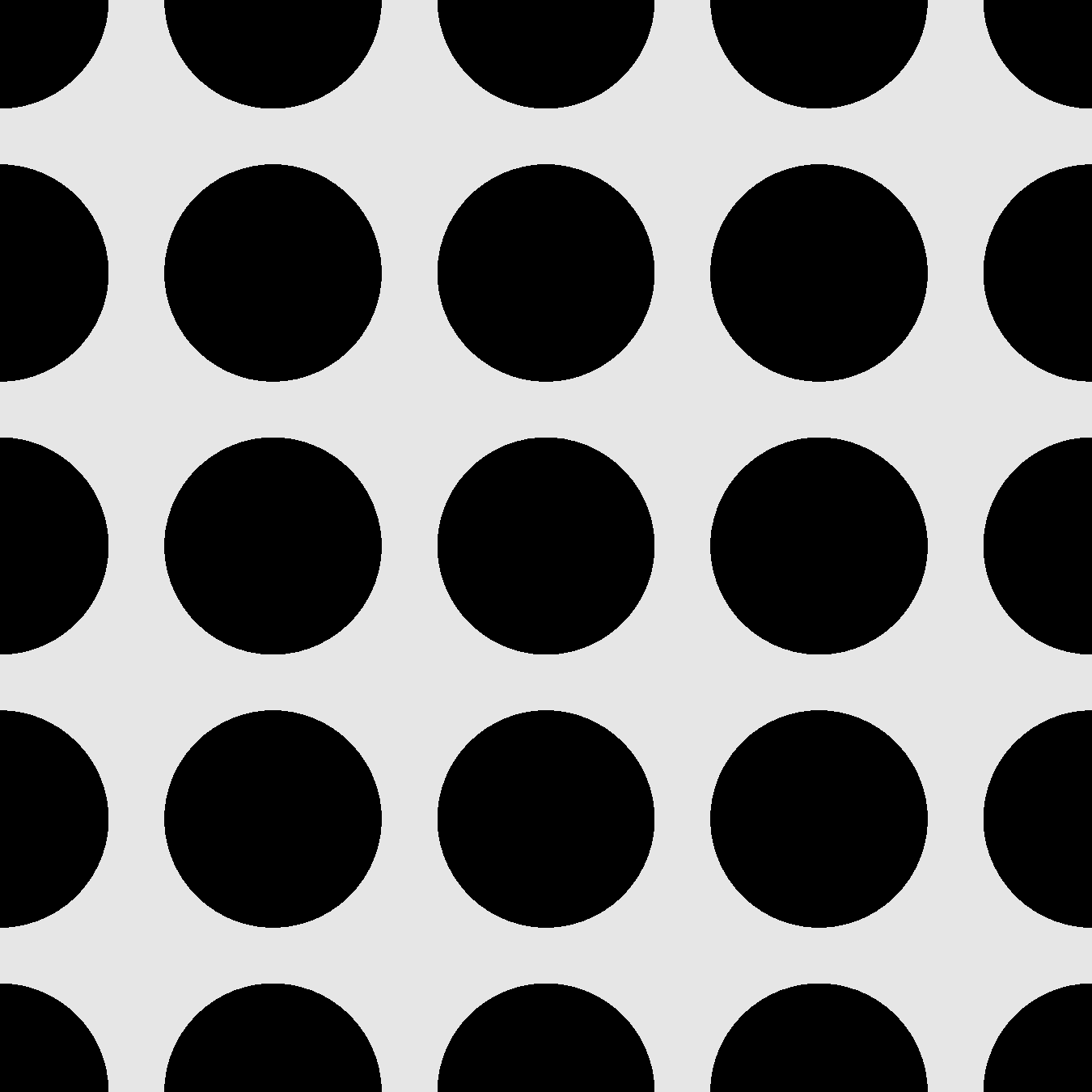}
	\hfill
	\includegraphics[width=0.48\linewidth]{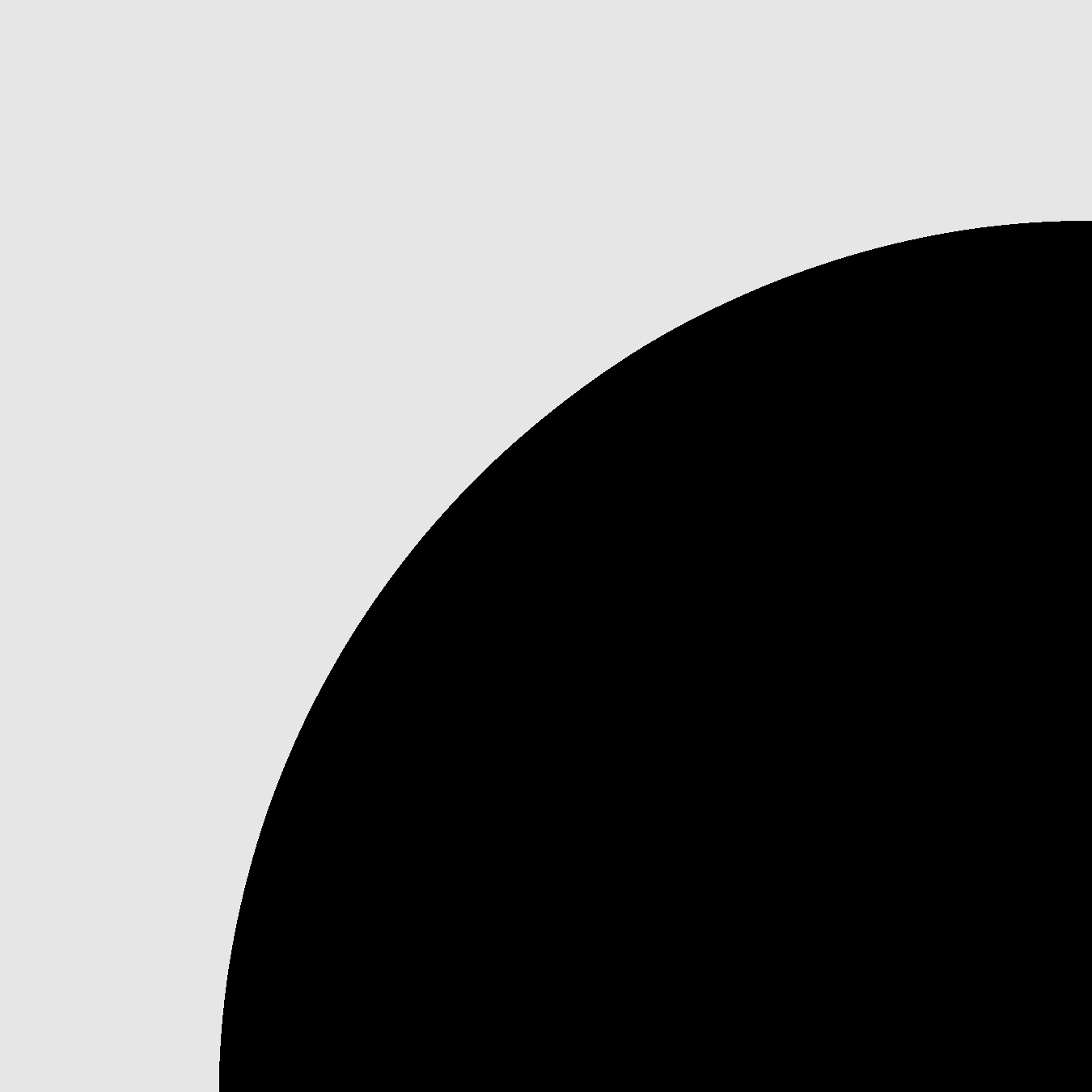}
	\includegraphics[width=0.48\linewidth]{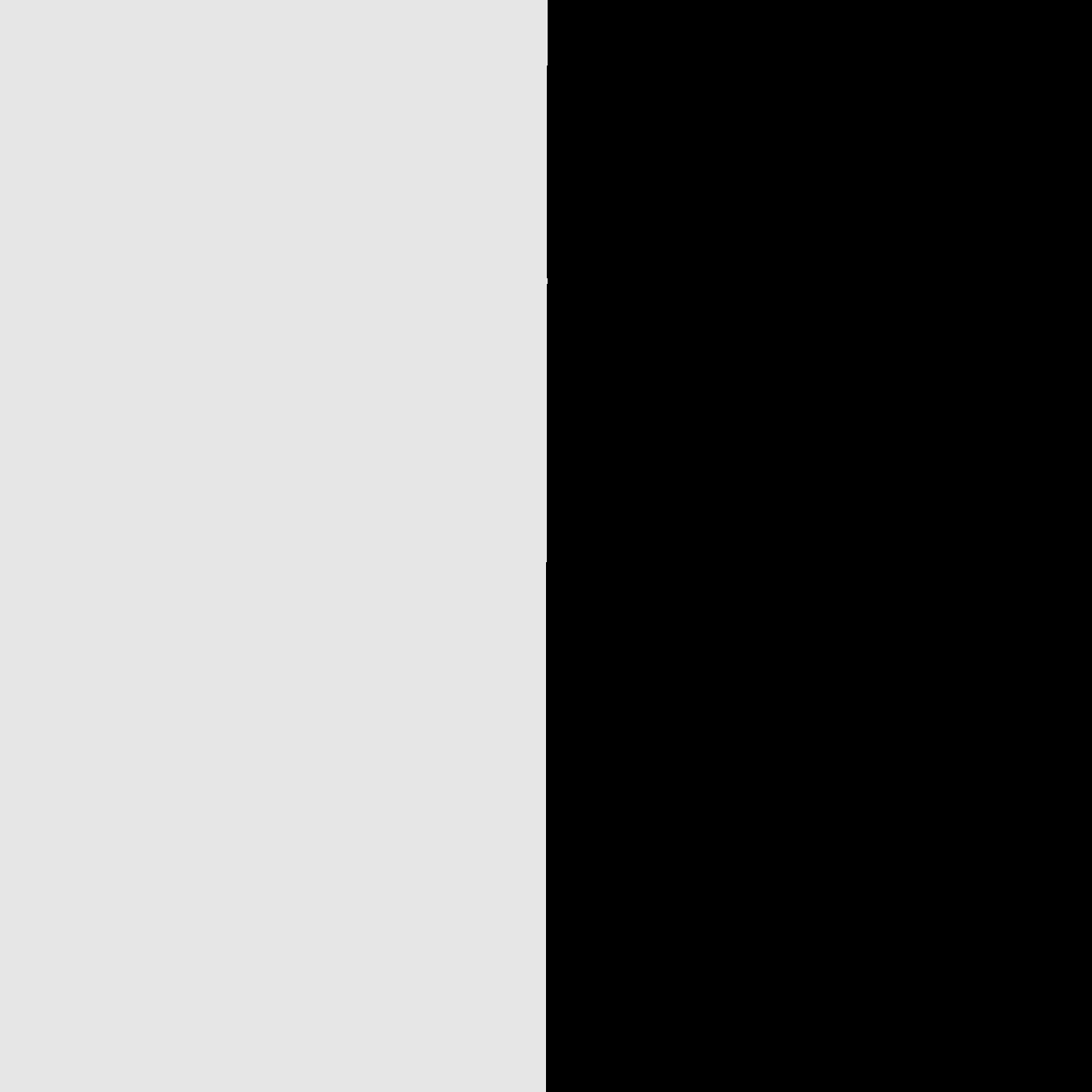}
	\hfill
	\includegraphics[width=0.48\linewidth]{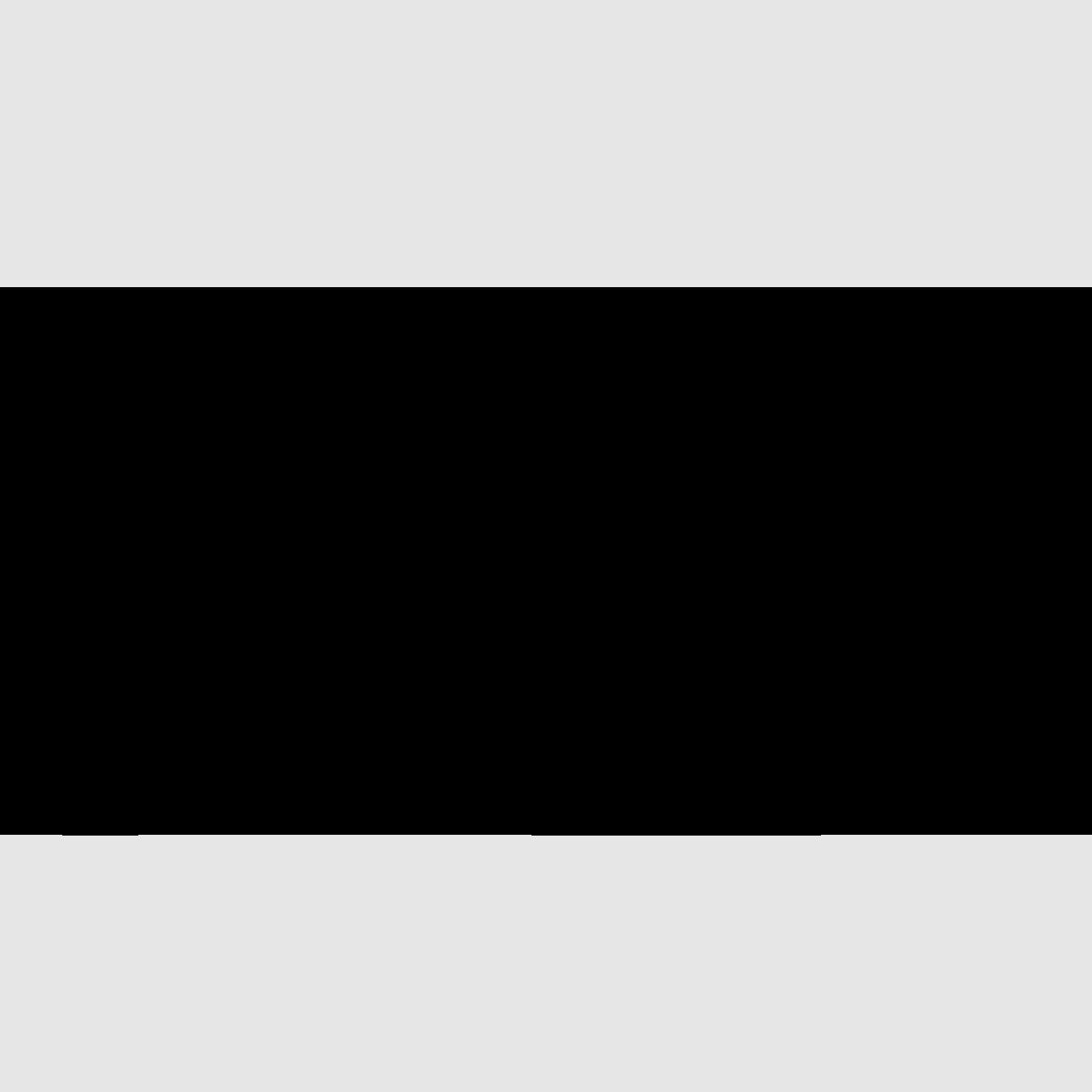}
	\captionof{figure}{
		Examples of resource allocation optimization in population dynamics.
		Initial guess (top-left) and optimized designs for
		$r=10$ (top-right), $r=40$ (bottom-left), and $r=100$ (bottom-right).
		The stopping criterion was satisfied at $147$, $107$, and $167$ iterations, respectively.
	}
	\label{fig:log}
	\end{figure}

\noindent\textit{Example.}
We solve problem \eqref{eq:log_cost}--\eqref{eq:log_prob} on the unit square $\mathcal{D} = (0,1)\times(0,1)$
with a resource constraint $V=0.5$. The bilinear form used in this problem is
\begin{equation*}
	B(\theta, \xi) \coloneqq
	\int_{\mathcal{D}} D\theta : D\xi +
	{10}^4 \int_{\partial \mathcal{D}}  (\theta \cdot n) (\xi \cdot n),	
\end{equation*}
with $\theta, \xi \in \mathbb{H}=H^1(\mathcal{D})$.
Thus, we allow the resource to move freely along the boundary $\partial \mathcal{D}$.
We perform tests with growth rates $r = 10$, $40$, and $100$, all using data parallelism across two processes:
\begin{lstlisting}
mpirun -np 2 python test.py 14 # r=10 (39 sec)
mpirun -np 2 python test.py 15 # r=40 (26 sec)
mpirun -np 2 python test.py 16 # r=100 (42 sec)	 
\end{lstlisting}
See the optimized resource distributions in Figure~\ref{fig:log}.

Newton's method was employed to solve the state equation.
For this purpose, the \code{pde} method returns the equation \eqref{eq:log_weak} in the form $F(u)=0$,
an empty list \code{[]} representing the absence of Dirichlet boundary contions,
the Jacobian $DF(u)(\delta u)$, the \code{UFL} variable that represents the state $u$,
and the positive function $u_0(x,y) = 1 + 0.2 \sin(6\pi x)\sin(6\pi y)$ as initial guess:
\begin{lstlisting}
def pde(self, phi):
	
	u = Coefficient(self.space)
	v = TestFunction(self.space)
	du = TrialFunction(self.space)

	F = dot(grad(u), grad(v))*self.dx
	F -= self.r*(1 - u/self.K(phi))*u*v*self.dx

	DF = dot(grad(du), grad(v))*self.dx
	DF -= self.r*(1 - 2*u/self.K(phi))*du*v*self.dx

	return [(F, [], DF, u, self.ini_func)]
\end{lstlisting}
The function \code{self.ini\_func} is a callable that will be interpolated on the domain mesh.
We provide it when a model of the \code{Logistic} class is created:
\begin{lstlisting}
u0 = lambda x: (
	1+0.2*np.sin(6*np.pi*x[0])*np.sin(6*np.pi*x[1])
)
md = Logistic(dim,domain,space,vol,r,u0,test_path)
\end{lstlisting}
Our choice of $u_0$ is motivated by the fact that 
$u$ represents a population density; hence, it must be positive and, additionally, cannot exceed the maximum value~$K$.

\subsection{Performance investigation}

\begin{figure}[H]\centering
\begin{minipage}{0.48\linewidth}
	\centering
	\includegraphics[width=\linewidth]{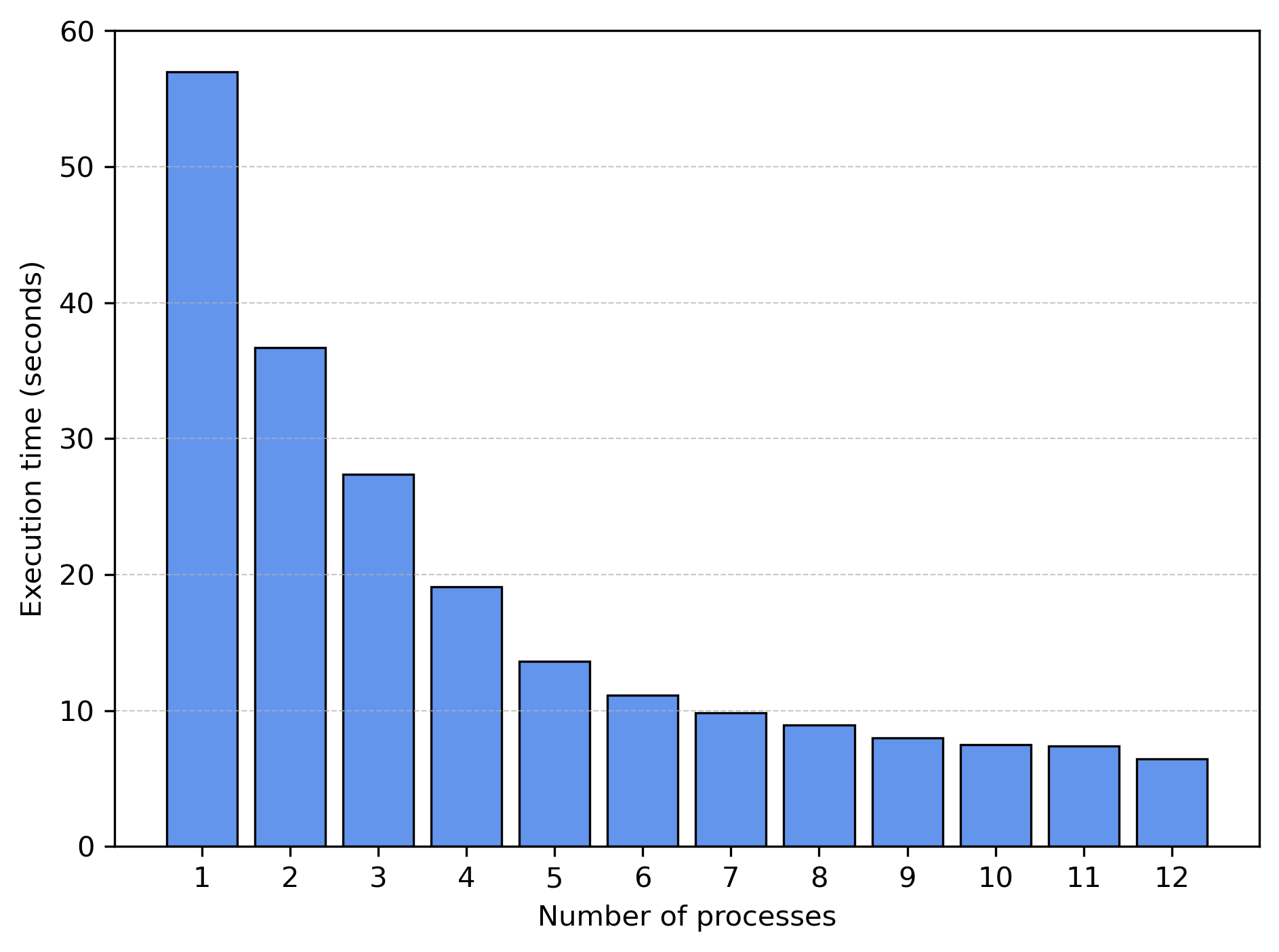}
\end{minipage}\hfill
\begin{minipage}{0.48\linewidth}
	\centering
	\includegraphics[width=0.9\linewidth]{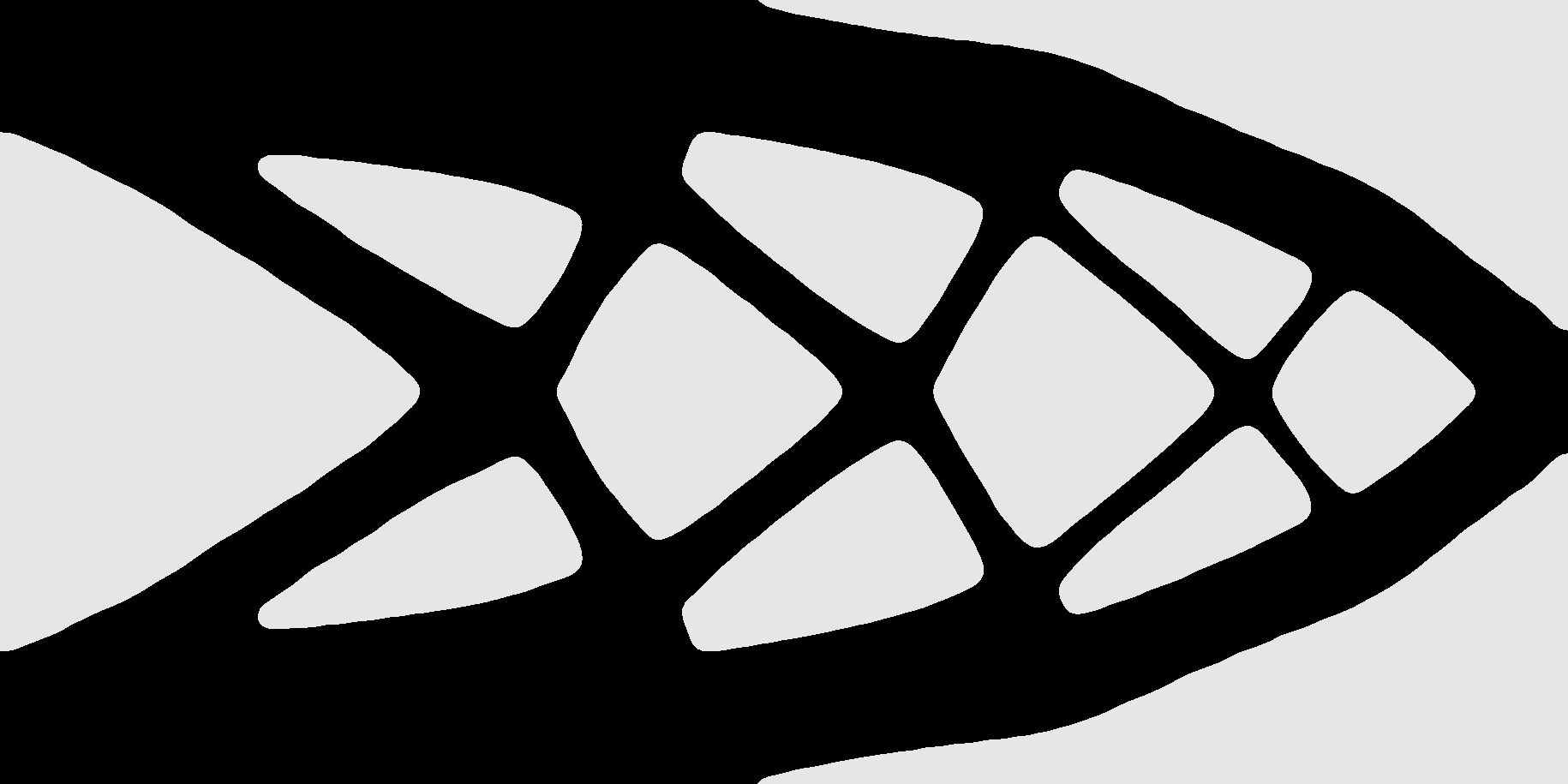}
\end{minipage}
\captionof{figure}{
	Performance results on a mesh with $52\,085$ triangles for the cantilever with one load. 
	Number of processes versus execution times in seconds (left)
	and optimized design at iteration $i=45$ (right).
}
\label{fig:perform1}
\end{figure}

\begin{figure}[H]\centering
\begin{minipage}{0.48\linewidth}
    \centering
    \includegraphics[width=\linewidth]{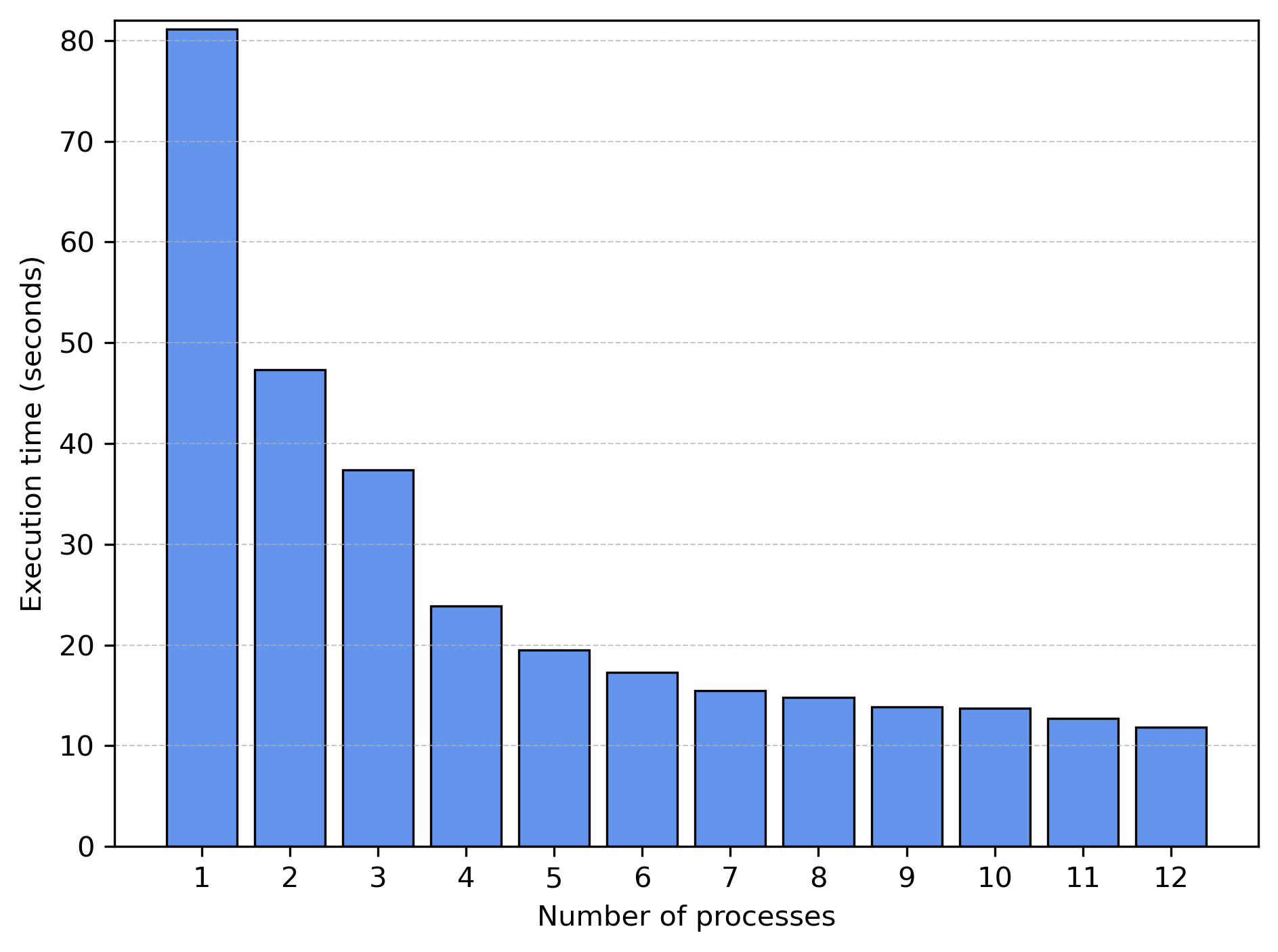}
\end{minipage}\hfill
\begin{minipage}{0.48\linewidth}
    \centering
    \includegraphics[width=0.70\linewidth]{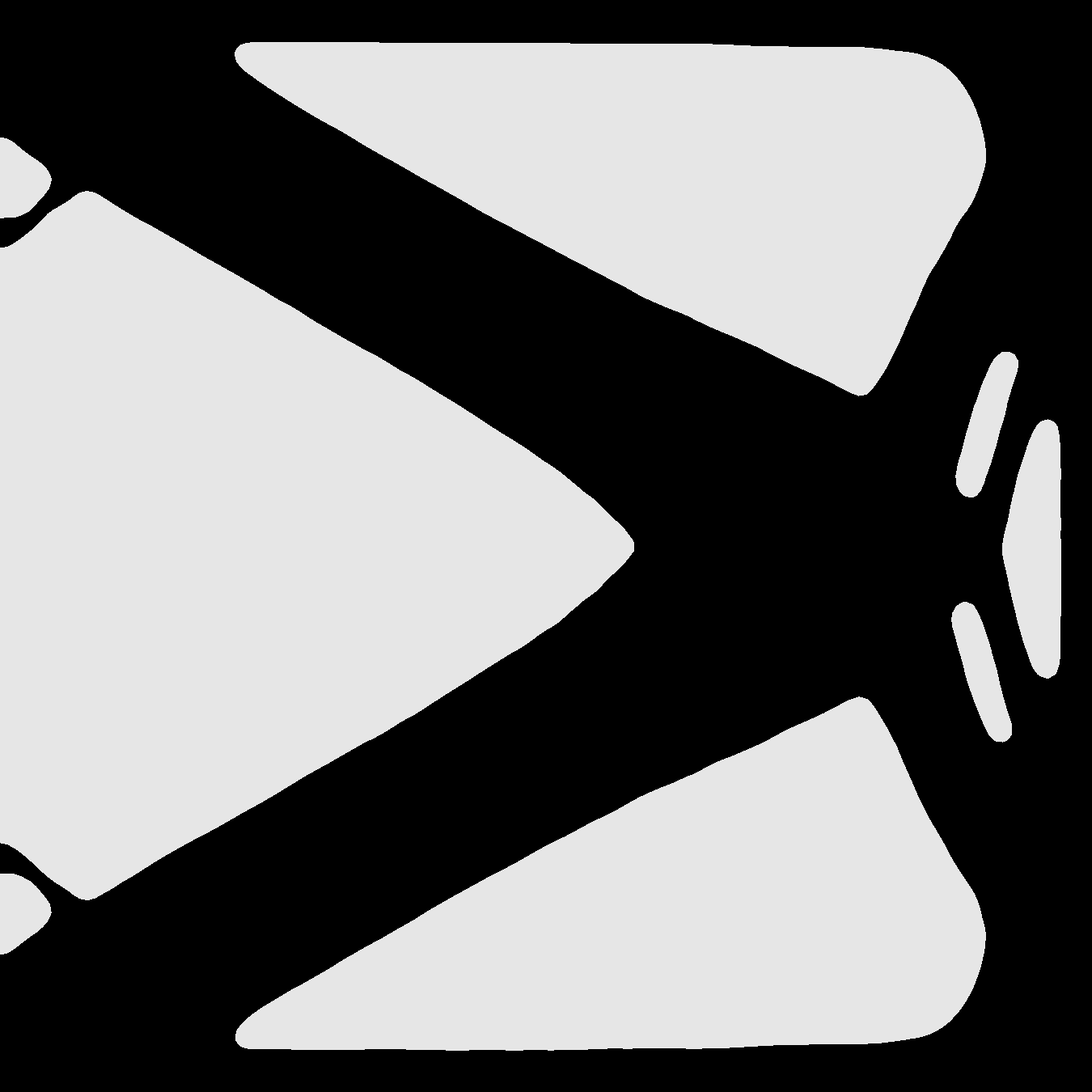}
\end{minipage}
\captionof{figure}{
	Performance results on a mesh with $36\,375$ triangles for the cantilever with two loads.
	Number of processes versus execution times in seconds (left)
	and optimized design at iteration $i=61$ (right).
}
\label{fig:perform2}
\end{figure}

We now present two performance tests aimed at assessing parallel scalability.
We revisit the symmetric cantilevers from \textit{Example~1} and \textit{Example~3} in Subsection~\ref{sub:compli},
with modified initial conditions to produce different optimization outcomes.
Both tests were run on the server using meshes with more than a 200\% increase in the number of finite elements compared with the settings of \textit{Example~1} and \textit{Example~3} shown in  Figure~\ref{fig:compli1} and  Figure~\ref{fig:compli3}.
The finer discretization allows for a greater number of smaller holes in the initial shapes, in contrast to the coarser examples discussed earlier.
The results are shown in Figures \ref{fig:perform1} and \ref{fig:perform2}.
The execution time decreased monotonically with the number of processes, up to 12 processes.
An almost linear speed-up was observed between 1 and 5 processes, with a pronounced reduction in execution time, followed by a slower, more gradual decrease at higher process counts.
Optimized designs with thinner structural parts were obtained.

\vspace{-0.8cm}
\section{Conclusion}
\vspace{-0.2cm}
This work presents a toolbox for PDE-constrained shape optimization based on the distributed shape derivative and a level-set method.
The implementation significantly generalizes the approach of \cite{MR3843884}: while inspired by the same underlying concepts, the present implementation is considerably more general in both scope and applicability.

The \fopt toolbox introduces several notable features.
The Unified Form Language (UFL) provides an intuitive way to represent weak formulations of PDEs in a notation close to the mathematical one. 
Our module encapsulates the numerical methods in classes
(level set, velocity problem, reinitialization),
which can be modified independently,
leaving the user with the task of writing the problem equations
(weak formulations, cost functional, constraints, derivative components)
in UFL format.
This modular structure separates the formulation of the equations from the definition of the mesh, finite elements, and boundary conditions.
\fopt  supports both two- and three-dimensional problems with flexible geometries, enabled by the finite element capabilities of \fx. Another key development is the integration of parallel computing with three distinct modes, which significantly enhances efficiency and enables large-scale simulations. In addition, the toolbox includes a built-in Proximal-Perturbed Lagrangian method for handling shape constraints, particularly useful for volume and perimeter constraints, which are ubiquitous in shape optimization.

Several extensions of this work will be pursued in future research. The extension to higher-order tensor representations of distributed shape derivatives is expected to be straightforward and will be valuable for applications, such as problems involving the bi-Laplacian~\cite{LL2024}. Incorporating tensor representations that include boundary terms is another important direction, though more challenging from an implementation perspective, as it requires advanced discretization strategies such as unfitted finite element methods, including CutFEM \cite{Burman1} and XFEM \cite{belytschko2001arbitrary}. 
Other relevant extensions include the treatment of time-dependent problems, the use of topological derivatives \cite{MR4916810,Novotny2013} and the integration of automatic differentiation tools \cite{Allaire2024,BLAUTH2023101577,MR4233697,Ham2019,Paganini2021} to facilitate the computation of adjoints and shape derivatives.

\end{multicols} 

\bibliographystyle{plain}
\bibliography{formopt}

@article {MR2278190,
    AUTHOR = {Zhuang, ChunGang and Xiong, ZhenHua and Ding, Han},
     TITLE = {A level set method for topology optimization of heat
              conduction problem under multiple load cases},
   JOURNAL = {Comput. Methods Appl. Mech. Engrg.},
  FJOURNAL = {Computer Methods in Applied Mechanics and Engineering},
    VOLUME = {196},
      YEAR = {2007},
    NUMBER = {4-6},
     PAGES = {1074--1084},
      ISSN = {0045-7825,1879-2138},
   MRCLASS = {80M50 (49N90 65M99)},
  MRNUMBER = {2278190},
       DOI = {10.1016/j.cma.2006.08.005},
       URL = {https://doi.org/10.1016/j.cma.2006.08.005},
}

@article{Allaire2024,
  title = {Autofreefem: automatic code generation with FreeFEM and LaTex output for shape and topology optimization of non-linear multi-physics problems},
  volume = {67},
  ISSN = {1615-1488},
  url = {http://dx.doi.org/10.1007/s00158-024-03917-5},
  DOI = {10.1007/s00158-024-03917-5},
  number = {12},
  journal = {Structural and Multidisciplinary Optimization},
  publisher = {Springer Science and Business Media LLC},
  author = {Allaire,  Grégoire and Gfrerer,  Michael H.},
  year = {2024},
  month = dec 
}

@misc{fenicsx,
  title     = {{DOLFINx}: the next generation {FEniCS} problem solving environment},
  author    = {Baratta, Igor A. and Dean, Joseph P. and Dokken, J{\o}rgen S. and Habera, Michal and Hale, Jack S. and Richardson, Chris N. and Rognes, Marie E. and Scroggs, Matthew W. and Sime, Nathan and Wells, Garth N.},
  doi       = {10.5281/zenodo.10447666},
  year      = {2023},
  howpublished = {preprint}
}

@article{Aage2014,
  title = {Topology optimization using PETSc: An easy-to-use,  fully parallel,  open source topology optimization framework},
  volume = {51},
  ISSN = {1615-1488},
  url = {http://dx.doi.org/10.1007/s00158-014-1157-0},
  DOI = {10.1007/s00158-014-1157-0},
  number = {3},
  journal = {Structural and Multidisciplinary Optimization},
  publisher = {Springer Science and Business Media LLC},
  author = {Aage,  Niels and Andreassen,  Erik and Lazarov,  Boyan Stefanov},
  year = {2014},
  month = aug,
  pages = {565–572}
}

@article{Ham2019,
  title = {Automated shape differentiation in the Unified Form Language},
  volume = {60},
  ISSN = {1615-1488},
  url = {http://dx.doi.org/10.1007/s00158-019-02281-z},
  DOI = {10.1007/s00158-019-02281-z},
  number = {5},
  journal = {Structural and Multidisciplinary Optimization},
  publisher = {Springer Science and Business Media LLC},
  author = {Ham,  David A. and Mitchell,  Lawrence and Paganini,  Alberto and Wechsung,  Florian},
  year = {2019},
  month = aug,
  pages = {1813–1820}
}

@article{KohnVogelius1984,
  author = {Robert V. Kohn and Michael Vogelius},
  title = {Determining conductivity by boundary measurements},
  journal = {Communications on Pure and Applied Mathematics},
  volume = {37},
  number = {3},
  pages = {289--298},
  year = {1984},
  doi = {10.1002/cpa.3160370302}
}

@article{Morassi2009,
  title = {Uniqueness and stability in determining a rigid inclusion in an elastic body},
  volume = {200},
  ISSN = {1947-6221},
  url = {http://dx.doi.org/10.1090/memo/0938},
  DOI = {10.1090/memo/0938},
  number = {938},
  journal = {Memoirs of the American Mathematical Society},
  publisher = {American Mathematical Society (AMS)},
  author = {Morassi,  Antonino and Rosset,  Edi},
  year = {2009},
  pages = {0–0}
}

@article {MR4261115,
    AUTHOR = {Albuquerque, Yuri F. and Laurain, Antoine and Yousept, Irwin},
     TITLE = {Level set-based shape optimization approach for
              sharp-interface reconstructions in time-domain full waveform
              inversion},
   JOURNAL = {SIAM J. Appl. Math.},
  FJOURNAL = {SIAM Journal on Applied Mathematics},
    VOLUME = {81},
      YEAR = {2021},
    NUMBER = {3},
     PAGES = {939--964},
      ISSN = {0036-1399},
   MRCLASS = {35Q93 (35R05 35R30 86A15)},
  MRNUMBER = {4261115},
       DOI = {10.1137/20M1378090},
       URL = {https://doi.org/10.1137/20M1378090},
}

@article {MR2225309,
    AUTHOR = {de Gournay, Fr{\'e}d{\'e}ric},
     TITLE = {Velocity extension for the level-set method and multiple
              eigenvalues in shape optimization},
   JOURNAL = {SIAM J. Control Optim.},
  FJOURNAL = {SIAM Journal on Control and Optimization},
    VOLUME = {45},
      YEAR = {2006},
    NUMBER = {1},
     PAGES = {343--367},
      ISSN = {0363-0129},
     CODEN = {SJCODC},
   MRCLASS = {74P10 (49Q12 74H99 74S30 90C90)},
  MRNUMBER = {2225309},
MRREVIEWER = {Andrzej M. My{\'s}li{\'n}ski},
       DOI = {10.1137/050624108},
       NOOPurl = {http://dx.doi.org/10.1137/050624108},
}

@book{Novotny2013,
  title = {Topological Derivatives in Shape Optimization},
  ISBN = {9783642352454},
  ISSN = {1860-6253},
  url = {http://dx.doi.org/10.1007/978-3-642-35245-4},
  DOI = {10.1007/978-3-642-35245-4},
  journal = {Interaction of Mechanics and Mathematics},
  publisher = {Springer Berlin Heidelberg},
  author = {Novotny,  Antonio André and Sokołowski,  Jan},
  year = {2013}
}

@article{vanDijk2013,
  title = {Level-set methods for structural topology optimization: a review},
  volume = {48},
  ISSN = {1615-1488},
  url = {http://dx.doi.org/10.1007/s00158-013-0912-y},
  DOI = {10.1007/s00158-013-0912-y},
  number = {3},
  journal = {Structural and Multidisciplinary Optimization},
  publisher = {Springer Science and Business Media LLC},
  author = {van Dijk,  N. P. and Maute,  K. and Langelaar,  M. and van Keulen,  F.},
  year = {2013},
  month = mar,
  pages = {437–472}
}

@article{belytschko2001arbitrary,
  title={Arbitrary discontinuities in finite elements},
  author={Belytschko, Ted and Mo{\"e}s, Nicolas and Usui, Shigeru and Parimi, Satyendra},
  journal={International Journal for Numerical Methods in Engineering},
  volume={50},
  number={4},
  pages={993--1013},
  year={2001},
  publisher={Wiley Online Library},
  doi={10.1002/nme.104}
}

@article{Burman1,
author = {Burman, Erik and Claus, Susanne and Hansbo, Peter and Larson, Mats G. and Massing, André},
title = {CutFEM: Discretizing geometry and partial differential equations},
journal = {International Journal for Numerical Methods in Engineering},
volume = {104},
number = {7},
pages = {472-501},
doi = {https://doi.org/10.1002/nme.4823},
url = {https://onlinelibrary.wiley.com/doi/abs/10.1002/nme.4823},
eprint = {https://onlinelibrary.wiley.com/doi/pdf/10.1002/nme.4823},
year = {2015}
}

@article{Mo2021,
  title = {Topology optimization of cooling plates for battery thermal management},
  volume = {178},
  ISSN = {0017-9310},
  url = {http://dx.doi.org/10.1016/j.ijheatmasstransfer.2021.121612},
  DOI = {10.1016/j.ijheatmasstransfer.2021.121612},
  journal = {International Journal of Heat and Mass Transfer},
  publisher = {Elsevier BV},
  author = {Mo,  Xiaobao and Zhi,  Hui and Xiao,  Yizhi and Hua,  Haiyu and He,  Liang},
  year = {2021},
  month = oct,
  pages = {121612}
}

@article{LamarcheGagnon2024,
  title = {Additively manufactured conformal cooling channels through topology optimization},
  volume = {67},
  ISSN = {1615-1488},
  url = {http://dx.doi.org/10.1007/s00158-024-03846-3},
  DOI = {10.1007/s00158-024-03846-3},
  number = {8},
  journal = {Structural and Multidisciplinary Optimization},
  publisher = {Springer Science and Business Media LLC},
  author = {Lamarche-Gagnon,  Marc-Étienne and Molavi-Zarandi,  Marjan and Raymond,  Vincent and Ilinca,  Florin},
  year = {2024},
  month = jul 
}

@book{10.5555/2789330,
author = {Ayachit, Utkarsh},
title = {The ParaView Guide: A Parallel Visualization Application},
year = {2015},
isbn = {1930934300},
publisher = {Kitware, Inc.},
address = {Clifton Park, NY, USA}
}

@article{GAO2008805,
title = {Topology optimization of heat conduction problem involving design-dependent heat load effect},
journal = {Finite Elements in Analysis and Design},
volume = {44},
number = {14},
pages = {805-813},
year = {2008},
issn = {0168-874X},
doi = {https://doi.org/10.1016/j.finel.2008.06.001},
url = {https://www.sciencedirect.com/science/article/pii/S0168874X08000905},
author = {T. Gao and W.H. Zhang and J.H. Zhu and Y.J. Xu and D.H. Bassir}
}

@article{BLAUTH2021100646,
title = {cashocs: A Computational, Adjoint-Based Shape Optimization and Optimal Control Software},
journal = {SoftwareX},
volume = {13},
pages = {100646},
year = {2021},
issn = {2352-7110},
doi = {https://doi.org/10.1016/j.softx.2020.100646},
url = {https://www.sciencedirect.com/science/article/pii/S2352711020303599},
author = {Sebastian Blauth}
}

@article{BLAUTH2023101577,
title = {Version 2.0 - cashocs: A Computational, Adjoint-Based Shape Optimization and Optimal Control Software},
journal = {SoftwareX},
volume = {24},
pages = {101577},
year = {2023},
issn = {2352-7110},
doi = {https://doi.org/10.1016/j.softx.2023.101577},
url = {https://www.sciencedirect.com/science/article/pii/S235271102300273X},
author = {Sebastian Blauth}
}

@article{MEFTAHI20151554,
title = {Sensitivity analysis for some inverse problems in linear elasticity via minimax differentiability},
journal = {Applied Mathematical Modelling},
volume = {39},
number = {5},
pages = {1554-1576},
year = {2015},
issn = {0307-904X},
doi = {https://doi.org/10.1016/j.apm.2014.09.026},
url = {https://www.sciencedirect.com/science/article/pii/S0307904X14004612},
author = {H. Meftahi and J.-P. Zolésio}
}

@article{Paganini2021,
  title = {Fireshape: a shape optimization toolbox for Firedrake},
  volume = {63},
  ISSN = {1615-1488},
  url = {http://dx.doi.org/10.1007/s00158-020-02813-y},
  DOI = {10.1007/s00158-020-02813-y},
  number = {5},
  journal = {Structural and Multidisciplinary Optimization},
  publisher = {Springer Science and Business Media LLC},
  author = {Paganini,  Alberto and Wechsung,  Florian},
  year = {2021},
  month = feb,
  pages = {2553–2569}
}

@article {MR4233697,
    AUTHOR = {Gangl, Peter and Sturm, Kevin and Neunteufel, Michael and
              Sch\"oberl, Joachim},
     TITLE = {Fully and semi-automated shape differentiation in {\tt
              {NGS}olve}},
   JOURNAL = {Struct. Multidiscip. Optim.},
  FJOURNAL = {Structural and Multidisciplinary Optimization},
    VOLUME = {63},
      YEAR = {2021},
    NUMBER = {3},
     PAGES = {1579--1607},
      ISSN = {1615-147X,1615-1488},
   MRCLASS = {49Q10 (49Q12 65N30 65Y15)},
  MRNUMBER = {4233697},
       DOI = {10.1007/s00158-020-02742-w},
       URL = {https://doi.org/10.1007/s00158-020-02742-w},
}

@book {MR2104179,
    AUTHOR = {Toselli, Andrea and Widlund, Olof},
     TITLE = {Domain decomposition methods---algorithms and theory},
    SERIES = {Springer Series in Computational Mathematics},
    VOLUME = {34},
 PUBLISHER = {Springer-Verlag, Berlin},
      YEAR = {2005},
     PAGES = {xvi+450},
      ISBN = {3-540-20696-5},
   MRCLASS = {65-02 (65N55 74S05 76M10)},
  MRNUMBER = {2104179},
MRREVIEWER = {R\'emi\ Vaillancourt},
       DOI = {10.1007/b137868},
       URL = {https://doi.org/10.1007/b137868},
}

@article {MR1772729,
    AUTHOR = {Ang, Dang Dinh and Trong, Dang Duc and Yamamoto, Masahiro},
     TITLE = {Identification of cavities inside two-dimensional
              heterogeneous isotropic elastic bodies},
   JOURNAL = {J. Elasticity},
  FJOURNAL = {Journal of Elasticity. The Physical and Mathematical Science
              of Solids},
    VOLUME = {56},
      YEAR = {1999},
    NUMBER = {3},
     PAGES = {199--212},
      ISSN = {0374-3535,1573-2681},
   MRCLASS = {74G75 (35R30)},
  MRNUMBER = {1772729},
MRREVIEWER = {Giovanni\ Franco\ Crosta},
       DOI = {10.1023/A:1007661505879},
       URL = {https://doi.org/10.1023/A:1007661505879},
}

@article{LAURAIN2020328,
author = "Laurain, A.",
title = "Distributed and boundary expressions of first and second order shape derivatives in nonsmooth domains",
journal = "Journal de Math{\'e}matiques Pures et Appliqu{\'e}es",
volume = "134",
pages = "328--368",
year = "2020",
doi = "10.1016/j.matpur.2019.09.002",
}

@article{LL2024,
author = {Laurain, Antoine  and Lopes, Pedro T. P. },
title = {On second-order tensor representation of derivatives in shape optimization},
journal = {Philosophical Transactions of the Royal Society A: Mathematical, Physical and Engineering Sciences},
volume = {382},
number = {2277},
pages = {20230300},
year = {2024},
doi = {10.1098/rsta.2023.0300},
URL = {https://royalsocietypublishing.org/doi/abs/10.1098/rsta.2023.0300},
eprint = {https://royalsocietypublishing.org/doi/pdf/10.1098/rsta.2023.0300}
}

@article {MR4531392,
    AUTHOR = {Laurain, A. and Lopes, P. T. P. and Nakasato, J. C.},
     TITLE = {An abstract {L}agrangian framework for computing shape
              derivatives},
   JOURNAL = {ESAIM. Control, Optimisation and Calculus of Variations},
    VOLUME = {29},
      YEAR = {2023},
     PAGES = {article 5},
      ISSN = {1292-8119,1262-3377},
   MRCLASS = {49Q10 (35Q93 35R37 49Q12)},
  MRNUMBER = {4531392},
       DOI = {10.1051/cocv/2022078},
       URL = {https://doi.org/10.1051/cocv/2022078},
}

@book {MR2512810,
    AUTHOR = {Henrot, A. and Pierre, M.},
     TITLE = {Variation et optimisation de formes},
    SERIES = {Math\'ematiques \& Applications (Berlin) [Mathematics \&
              Applications]},
    VOLUME = {48},
      NOTE = {Une analyse g{\'e}om{\'e}trique. [A geometric analysis]},
 PUBLISHER = {Springer},
   ADDRESS = {Berlin},
      YEAR = {2005},
     PAGES = {xii+334},
      ISBN = {978-3-540-26211-4; 3-540-26211-3},
   MRCLASS = {49-02 (35A15 35J20 35R35 49K20 49Q05 49Q10 74P10)},
  MRNUMBER = {2512810 (2009m:49003)},
}

@article {MR2108636,
    AUTHOR = {Wang, Michael Yu and Zhou, Shiwei},
     TITLE = {Phase field: a variational method for structural topology
              optimization},
   JOURNAL = {CMES Comput. Model. Eng. Sci.},
  FJOURNAL = {CMES. Computer Modeling in Engineering \& Sciences},
    VOLUME = {6},
      YEAR = {2004},
    NUMBER = {6},
     PAGES = {547--566},
      ISSN = {1526-1492},
   MRCLASS = {74P15 (49S05 74G15 74N99 74Q05 90C90)},
  MRNUMBER = {2108636},
}

@article {MR1951408,
    AUTHOR = {Wang, Michael Yu and Wang, Xiaoming and Guo, Dongming},
     TITLE = {A level set method for structural topology optimization},
   JOURNAL = {Comput. Methods Appl. Mech. Engrg.},
  FJOURNAL = {Computer Methods in Applied Mechanics and Engineering},
    VOLUME = {192},
      YEAR = {2003},
    NUMBER = {1-2},
     PAGES = {227--246},
      ISSN = {0045-7825},
     CODEN = {CMMECC},
   MRCLASS = {74S30 (49Q10 65M99 74P05)},
  MRNUMBER = {1951408},
MRREVIEWER = {Michal Ko{\v{c}}vara},
       DOI = {10.1016/S0045-7825(02)00559-5},
       URL = {http://dx.doi.org/10.1016/S0045-7825(02)00559-5},
}

@article {MR2235384,
    AUTHOR = {Amstutz, Samuel and Andr{\"a}, Heiko},
     TITLE = {A new algorithm for topology optimization using a level-set
              method},
   JOURNAL = {J. Comput. Phys.},
  FJOURNAL = {Journal of Computational Physics},
    VOLUME = {216},
      YEAR = {2006},
    NUMBER = {2},
     PAGES = {573--588},
      ISSN = {0021-9991},
     CODEN = {JCTPAH},
   MRCLASS = {65K10 (49Q10 74P05 74P10)},
  MRNUMBER = {2235384},
       DOI = {10.1016/j.jcp.2005.12.015},
       URL = {http://dx.doi.org/10.1016/j.jcp.2005.12.015},
}

@article {MR965860,
    AUTHOR = {Osher, Stanley and Sethian, James A.},
     TITLE = {Fronts propagating with curvature-dependent speed: algorithms
              based on {H}amilton-{J}acobi formulations},
   JOURNAL = {J. Comput. Phys.},
  FJOURNAL = {Journal of Computational Physics},
    VOLUME = {79},
      YEAR = {1988},
    NUMBER = {1},
     PAGES = {12--49},
      ISSN = {0021-9991},
     CODEN = {JCTPAH},
   MRCLASS = {80A25 (80-08)},
  MRNUMBER = {MR965860 (89h:80012)},
}

@article {MR3348199,
    AUTHOR = {Hiptmair, R. and Paganini, A. and Sargheini, S.},
     TITLE = {Comparison of approximate shape gradients},
   JOURNAL = {BIT},
  FJOURNAL = {BIT. Numerical Mathematics},
    VOLUME = {55},
      YEAR = {2015},
    NUMBER = {2},
     PAGES = {459--485},
      ISSN = {0006-3835},
   MRCLASS = {65N30 (49Q12)},
  MRNUMBER = {3348199},
       DOI = {10.1007/s10543-014-0515-z},
       URL = {http://dx.doi.org/10.1007/s10543-014-0515-z},
}

@incollection {MR2642680,
    AUTHOR = {Berggren, Martin},
     TITLE = {A unified discrete-continuous sensitivity analysis method for
              shape optimization},
 BOOKTITLE = {Applied and numerical partial differential equations},
    SERIES = {Comput. Methods Appl. Sci.},
    VOLUME = {15},
     PAGES = {25--39},
 PUBLISHER = {Springer, New York},
      YEAR = {2010},
   MRCLASS = {49Q10 (35A15 65N30 76M25)},
  MRNUMBER = {2642680 (2011g:49080)},
       DOI = {10.1007/978-90-481-3239-3_4},
       URL = {http://dx.doi.org/10.1007/978-90-481-3239-3_4},
}

@article {MR3535238,
    AUTHOR = {Laurain, A. and Sturm, K.},
     TITLE = {Distributed shape derivative {\it via} averaged adjoint method
              and applications},
   JOURNAL = {ESAIM. Mathematical Modelling and Numerical Analysis},
    VOLUME = {50},
      YEAR = {2016},
    NUMBER = {4},
     PAGES = {1241--1267},
      ISSN = {0764-583X},
   MRCLASS = {49Q10 (35R30)},
  MRNUMBER = {3535238},
MRREVIEWER = {Faustino Maestre},
       DOI = {10.1051/m2an/2015075},
       URL = {https://doi.org/10.1051/m2an/2015075},
}

@article {MR2459656,
    AUTHOR = {Fulma{\'n}ski, Piotr and Laurain, Antoine and Scheid,
              Jean-Fran{\c{c}}ois and Soko{\l}owski, Jan},
     TITLE = {Level set method with topological derivatives in shape
              optimization},
   JOURNAL = {Int. J. Comput. Math.},
  FJOURNAL = {International Journal of Computer Mathematics},
    VOLUME = {85},
      YEAR = {2008},
    NUMBER = {10},
     PAGES = {1491--1514},
      ISSN = {0020-7160},
   MRCLASS = {49Q10 (35J05 35J20 65K10 74P15)},
  MRNUMBER = {2459656},
       DOI = {10.1080/00207160802033350},
       NOOPurl = {http://dx.doi.org/10.1080/00207160802033350},
}

@Article{Challis2009,
author="Challis, Vivien J.",
title="A discrete level-set topology optimization code written in Matlab",
journal="Structural and Multidisciplinary Optimization",
year="2009",
volume="41",
number="3",
pages="453--464",
}

@book {SokZol92,
    AUTHOR = {Soko{\l}owski, Jan and Zol{\'e}sio, Jean-Paul},
     TITLE = {Introduction to shape optimization},
    SERIES = {Springer Series in Computational Mathematics},
    VOLUME = {16},
      NOTE = {Shape sensitivity analysis},
 PUBLISHER = {Springer-Verlag, Berlin},
      YEAR = {1992},
     PAGES = {ii+250},
      ISBN = {3-540-54177-2},
   MRCLASS = {49-02 (49Q15 73K40 73V25)},
  MRNUMBER = {1215733 (94d:49002)},
MRREVIEWER = {Fredi Tr{\"o}ltzsch},
       DOI = {10.1007/978-3-642-58106-9},
       NOOPurl = {http://dx.doi.org/10.1007/978-3-642-58106-9},
}

@book {MR2731611,
    AUTHOR = {Delfour, M. C. and Zol{\'e}sio, J.-P.},
     TITLE = {Shapes and geometries},
    SERIES = {Advances in Design and Control},
    VOLUME = {22},
   EDITION = {Second},
      NOTE = {Metrics, analysis, differential calculus, and optimization},
 PUBLISHER = {Society for Industrial and Applied Mathematics (SIAM)},
   ADDRESS = {Philadelphia, PA},
      YEAR = {2011},
     PAGES = {xxiv+622},
      ISBN = {978-0-898719-36-9},
   MRCLASS = {49-02 (26-02 54-02)},
  MRNUMBER = {2731611 (2012c:49002)},
MRREVIEWER = {Jan Soko{\l}owski},
       DOI = {10.1137/1.9780898719826},
       NOOPurl = {http://dx.doi.org/10.1137/1.9780898719826},
}

@article {MR2033390,
    AUTHOR = {Allaire, Gr{\'e}goire and Jouve, Fran{\c{c}}ois and Toader,
              Anca-Maria},
     TITLE = {Structural optimization using sensitivity analysis and a
              level-set method},
   JOURNAL = {J. Comput. Phys.},
  FJOURNAL = {Journal of Computational Physics},
    VOLUME = {194},
      YEAR = {2004},
    NUMBER = {1},
     PAGES = {363--393},
      ISSN = {0021-9991},
     CODEN = {JCTPAH},
   MRCLASS = {74P15 (49Q12 65M06 74G15)},
  MRNUMBER = {2033390},
MRREVIEWER = {Andrzej M. My{\'s}li{\'n}ski},
       DOI = {10.1016/j.jcp.2003.09.032},
       NOOPurl = {http://dx.doi.org/10.1016/j.jcp.2003.09.032},
}

@article {MR1911658,
    AUTHOR = {Allaire, Gr{\'e}goire and Jouve, Fran{\c{c}}ois and Toader,
              Anca-Maria},
     TITLE = {A level-set method for shape optimization},
   JOURNAL = {C. R. Math. Acad. Sci. Paris},
  FJOURNAL = {Comptes Rendus Math\'ematique. Acad\'emie des Sciences. Paris},
    VOLUME = {334},
      YEAR = {2002},
    NUMBER = {12},
     PAGES = {1125--1130},
      ISSN = {1631-073X},
   MRCLASS = {49Q10 (74G15 74P10)},
  MRNUMBER = {1911658},
       DOI = {10.1016/S1631-073X(02)02412-3},
       NOOPurl = {http://dx.doi.org/10.1016/S1631-073X(02)02412-3},
}

@article {MR2252743,
    AUTHOR = {Allaire, G. and Pantz, O.},
     TITLE = {Structural optimization with {F}ree{F}em++},
   JOURNAL = {Struct. Multidiscip. Optim.},
  FJOURNAL = {Structural and Multidisciplinary Optimization},
    VOLUME = {32},
      YEAR = {2006},
    NUMBER = {3},
     PAGES = {173--181},
      ISSN = {1615-147X},
   MRCLASS = {74P05 (74P15 74Q05 74S05)},
  MRNUMBER = {2252743},
       DOI = {10.1007/s00158-006-0017-y},
       NOOPurl = {http://dx.doi.org/10.1007/s00158-006-0017-y},
}

@Article{Gain2013,
author="Gain, Arun L.
and Paulino, Glaucio H.",
title="A critical comparative assessment of differential equation-driven methods for structural topology optimization",
journal="Structural and Multidisciplinary Optimization",
year="2013",
volume="48",
number="4",
pages="685--710",
issn="1615-1488",
doi="10.1007/s00158-013-0935-4",
NOOPurl="http://dx.doi.org/10.1007/s00158-013-0935-4",
}

@Article{Andreassen2010,
author="Andreassen, Erik
and Clausen, Anders
and Schevenels, Mattias
and Lazarov, Boyan S.
and Sigmund, Ole",
title="Efficient topology optimization in MATLAB using 88 lines of code",
journal="Structural and Multidisciplinary Optimization",
year="2010",
volume="43",
number="1",
pages="1--16",
}

@article{Sigmund2014,
author="Sigmund, O.",
title="A 99 line topology optimization code written in Matlab",
journal="Structural and Multidisciplinary Optimization",
year="2014",
volume="21",
number="2",
pages="120--127",
}

@book {MR2008524,
    AUTHOR = {Bends{\o}e, M. P. and Sigmund, O.},
     TITLE = {Topology optimization},
      NOTE = {Theory, methods and applications},
 PUBLISHER = {Springer-Verlag, Berlin},
      YEAR = {2003},
     PAGES = {xiv+370},
      ISBN = {3-540-42992-1},
   MRCLASS = {74P15 (49Q10 74P05 74Q05)},
  MRNUMBER = {2008524},
}

@article {MR4686093,
    AUTHOR = {Feppon, Florian},
     TITLE = {Density-based topology optimization with the null space
              optimizer: a tutorial and a comparison},
   JOURNAL = {Struct. Multidiscip. Optim.},
  FJOURNAL = {Structural and Multidisciplinary Optimization},
    VOLUME = {67},
      YEAR = {2024},
    NUMBER = {1},
     PAGES = {Paper No. 4, 34},
      ISSN = {1615-147X,1615-1488},
   MRCLASS = {74P15},
  MRNUMBER = {4686093},
       DOI = {10.1007/s00158-023-03710-w},
       URL = {https://doi.org/10.1007/s00158-023-03710-w},
}

@article {MR4916810,
    AUTHOR = {Filho, J. M. M. Luz and Gomes, A. T. A. and Novotny, A. A.},
     TITLE = {A parallel {F}ree{FEM} framework for topology optimization of
              structures into three spatial dimensions},
   JOURNAL = {Finite Elem. Anal. Des.},
  FJOURNAL = {Finite Elements in Analysis and Design},
    VOLUME = {250},
      YEAR = {2025},
     PAGES = {Paper No. 104384},
      ISSN = {0168-874X,1872-6925},
   MRCLASS = {65M60 (65Y05)},
  MRNUMBER = {4916810},
       DOI = {10.1016/j.finel.2025.104384},
       URL = {https://doi.org/10.1016/j.finel.2025.104384},
}

@article {MR1640142,
    AUTHOR = {Barth, Timothy J. and Sethian, James A.},
     TITLE = {Numerical schemes for the {H}amilton-{J}acobi and level set
              equations on triangulated domains},
   JOURNAL = {J. Comput. Phys.},
  FJOURNAL = {Journal of Computational Physics},
    VOLUME = {145},
      YEAR = {1998},
    NUMBER = {1},
     PAGES = {1--40},
      ISSN = {0021-9991,1090-2716},
   MRCLASS = {65M60},
  MRNUMBER = {1640142},
       DOI = {10.1006/jcph.1998.6007},
       URL = {https://doi.org/10.1006/jcph.1998.6007},
}

@article {MR3843884,
    AUTHOR = {Laurain, A.},
     TITLE = {A level set-based structural optimization code using
              {FE}ni{CS}},
   JOURNAL = {Structural and Multidisciplinary Optimization},
    VOLUME = {58},
      YEAR = {2018},
    NUMBER = {3},
     PAGES = {1311--1334},
      ISSN = {1615-147X,1615-1488},
   MRCLASS = {74P05},
  MRNUMBER = {3843884},
       DOI = {10.1007/s00158-018-1950-2},
       URL = {https://doi.org/10.1007/s00158-018-1950-2},
}

@article{Jia2024,
  title = {FEniTop: a simple FEniCSx implementation for 2D and 3D topology optimization supporting parallel computing},
  volume = {67},
  ISSN = {1615-1488},
  url = {http://dx.doi.org/10.1007/s00158-024-03818-7},
  DOI = {10.1007/s00158-024-03818-7},
  number = {8},
  journal = {Structural and Multidisciplinary Optimization},
  publisher = {Springer Science and Business Media LLC},
  author = {Jia,  Yingqi and Wang,  Chao and Zhang,  Xiaojia Shelly},
  year = {2024},
  month = aug 
}

@article{Wegert2025,
  title = {GridapTopOpt.jl: a scalable Julia toolbox for level set-based topology optimisation},
  volume = {68},
  ISSN = {1615-1488},
  url = {http://dx.doi.org/10.1007/s00158-024-03927-3},
  DOI = {10.1007/s00158-024-03927-3},
  number = {1},
  journal = {Structural and Multidisciplinary Optimization},
  publisher = {Springer Science and Business Media LLC},
  author = {Wegert,  Zachary J. and Manyer,  Jordi and Mallon,  Connor N. and Badia,  Santiago and Challis,  Vivien J.},
  year = {2025},
  month = jan 
}

@article {MR4589221,
    AUTHOR = {Kim, Jong Gwang},
     TITLE = {A new {L}agrangian-based first-order method for nonconvex
              constrained optimization},
   JOURNAL = {Oper. Res. Lett.},
  FJOURNAL = {Operations Research Letters},
    VOLUME = {51},
      YEAR = {2023},
    NUMBER = {3},
     PAGES = {357--363},
      ISSN = {0167-6377,1872-7468},
   MRCLASS = {90C26 (49M29 65K10)},
  MRNUMBER = {4589221},
       DOI = {10.1016/j.orl.2023.04.006},
       URL = {https://doi.org/10.1016/j.orl.2023.04.006},
}

@book {MR1908418,
    AUTHOR = {Murray, J. D.},
     TITLE = {Mathematical biology. {I}},
    SERIES = {Interdisciplinary Applied Mathematics},
    VOLUME = {17},
   EDITION = {Third},
      NOTE = {An introduction},
 PUBLISHER = {Springer-Verlag, New York},
      YEAR = {2002},
     PAGES = {xxiv+551},
      ISBN = {0-387-95223-3},
   MRCLASS = {92B05 (92-01)},
  MRNUMBER = {1908418},
MRREVIEWER = {Trachette\ L.\ Jackson},
}

\end{document}